\documentclass[preprint, 11pt]{elsarticle}

\usepackage[T1]{fontenc}
\usepackage[utf8]{inputenc}
\usepackage[hmargin=0.9in]{geometry}
\usepackage[babel,german=quotes]{csquotes}
\usepackage{subcaption}
\usepackage{graphicx}
\usepackage[colorlinks=true,citecolor=blue,urlcolor=blue]{hyperref}
\usepackage{caption}
\usepackage{amsmath}
\usepackage{amsthm}
\usepackage{amssymb}
\usepackage{bm, amsfonts}
\usepackage{xcolor,soul,framed}
\usepackage{fixmath}
\usepackage{tikz}
\usepackage{bm}
\usepackage{makecell}
\usepackage{titlesec}
\usepackage{multirow}
\usepackage{indentfirst}
\usepackage{mathtools}
\usepackage{float}
\usepackage[normalem]{ulem}
\useunder{\uline}{\ul}{}

\newcommand{\diff}{\mathrm{d} }
\newcommand{\R}{\mathbb{R}}
\newcommand{\N}{\mathbb{N}}

\newcommand{\x}{\mathbf{x}}
\newcommand{\s}{\mathbf{s}}
\newcommand{\vel}{\mathbf{v}}
\newcommand{\n}{\mathbf{n}}
\newcommand{\f}{\mathbf{f}}

\newcommand{\GG}{\mathcal{G}}

\newcommand{\E}{\mathcal{E}}
\newcommand{\con}{\mathbf{c}}
\newcommand{\A}{\mathbf{A}}

\newcommand{\J}{J_L}

\DeclareMathOperator*{\argmin}{arg\,min}

\colorlet{shadecolor}{green}

\titleformat{\paragraph}
{\normalfont\normalsize\itshape}{\theparagraph}{1em}{}
\titlespacing*{\paragraph}
{0pt}{3.25ex plus 1ex minus .2ex}{1.5ex plus .2ex}

\newdefinition{remark}{Remark}
\newtheorem{theorem}{Theorem}

\makeatletter
\def\ps@pprintTitle{%
\let\@oddhead\@empty
\let\@evenhead\@empty
\def\@oddfoot{
\footnotesize\itshape
\ifx\@journal\@empty Elsevier
\else\@journal\fi
\hfill\today
}%
\let\@evenfoot\@oddfoot}
\makeatother

\newcommand\red[1]{#1}
\newcommand\blue[1]{#1}
\newcommand\cyan[1]{#1}
\newcommand\magenta[1]{#1}

  \newcommand{\pd}[2]{\frac{\partial #1}{\partial #2}}
\newcommand{\dx}{\mathrm{d}\mathbf{x}}

\begin{document}

\begin{frontmatter}
\title{Monolithic convex limiting and implicit pseudo-time stepping for calculating steady-state solutions of the Euler equations}
\author{Paul Moujaes, Dmitri Kuzmin}
\ead{paul.moujaes@math.tu-dortmund.de, kuzmin@math.uni-dortmund.de}

\address{Institute of Applied Mathematics (LS III), TU Dortmund University\\ Vogelpothsweg 87,
	D-44227 Dortmund, Germany}

\journal{Journal of Computational Physics}

\begin{abstract}
  In this work, we use the monolithic convex limiting (MCL) methodology to enforce relevant inequality constraints in implicit finite element discretizations of the compressible Euler equations. In this context, preservation of invariant domains follows from positivity preservation for intermediate states of the density and internal energy. To avoid spurious oscillations, we additionally impose local maximum principles on intermediate states of the density, velocity components, and specific total energy. For the backward Euler time stepping, we show the invariant domain preserving (IDP) property of the fully discrete MCL scheme by constructing a fixed-point iteration \red{that meets the requirements of a Krasnoselskii-type theorem.} Our iterative solver for the nonlinear discrete problem employs a more efficient fixed-point iteration. The matrix of the associated linear system is a robust low-order Jacobian approximation that exploits the homogeneity property of the flux function. The limited antidiffusive terms are treated explicitly. We use positivity preservation as a stopping criterion for nonlinear iterations. The first iteration yields the solution of a linearized semi-implicit problem. This solution possesses the discrete conservation property but is generally not IDP. Further iterations are performed if any non-IDP states are detected. The existence of an IDP limit is guaranteed by our analysis. To facilitate convergence to steady-state solutions, we perform adaptive explicit underrelaxation at the end of each time step. The calculation of appropriate relaxation factors is based on an approximate minimization of nodal entropy residuals. The performance of proposed algorithms and alternative solution strategies is illustrated by the convergence history for standard two-dimensional test problems. 
\end{abstract}


\end{frontmatter}

\section{Introduction}
\label{sec:intro}
Recent years have witnessed significant advances in the development of bound-preserving
schemes for the Euler equations of gas dynamics. In particular, flux-corrected transport (FCT)
algorithms that ensure nonnegativity of continuous finite element approximations to the density
and pressure (internal energy) were proposed in \cite{dobrev2018,guermond2018,kuzmin2010a,lohmann2016}. A potential drawback of these explicit predictor-corrector approaches is the fact that they require the use of small time steps and do not converge to steady-state solutions.
The implicit algebraic flux correction schemes employed by Gurris et al. \cite{gurris2012} do not have this limitation. However, 
  positivity preservation is not guaranteed even for converged solutions.

The monolithic convex limiting (MCL) methodology introduced in \cite{kuzmin2020} supports
the use of general time integrators. The residuals of nonlinear discrete problems are well
defined both for individual time steps and for a steady-state limit. If time marching is performed using an explicit strong stability (SSP)
preserving Runge--Kutta method \cite{gottlieb2001}, the invariant
domain preservation (IDP) property\footnote{Preservation of invariant domains, as defined
  in \cite{guermond2018,guermond2016}, is a formal synonym for positivity preservation.} of Shu--Osher
stages can easily be shown following
the analysis of FCT-type convex limiting approaches in \cite{guermond2018}. The
analysis of implicit MCL schemes for linear advection equations can be performed
using the theoretical framework developed in \cite{barrenechea2016,lohmann2019}. So far, no formal proof of the IDP property was provided for implicit MCL
discretizations of the Euler equations. \red{We fill this gap in the present paper by applying a Krasnoselskii-type theorem to a fixed-point iteration that we design to be IDP.}

Having established the existence of an IDP solution to our nonlinear discrete
problem, we discuss ways to calculate it in practice. Building on our
previous experience with the design of iterative solvers for flux-limited
finite element discretizations of the Euler system \cite{gurris2012,kuzmin2012b},
we update intermediate solutions using a low-order approximation to
the Jacobian. The underlying linearization of nodal fluxes was proposed by
Dolej{\v{s}}{\i} and Feistauer \cite{dolejvsi2004}, who performed a single
iteration per time step. In our solver, we perform as many iterations as it
takes to obtain a positivity-preserving result. That is, we use the IDP
property as a natural stopping criterion. Another highlight of the
proposed approach is a new kind of explicit underrelaxation for solution changes.
We select relaxation factors from a discrete set and choose the value that
corresponds to the minimum entropy residual. This selection criterion was
inspired by the work of Ranocha et al. \cite{ranocha2020} on relaxation
Runge--Kutta methods that ensure fully discrete entropy stability for
time-dependent nonlinear problems. In our numerical studies, convergence
is achieved for CFL numbers as high as $10^5$ and Mach numbers as high as 20.

The remainder of this article is organized as follows. We present the low-order
component of our implicit MCL scheme in Section \ref{sec:lodisc} and prove its IDP
property in Section \ref{sec:IDP}. The proof admits a straightforward
generalization to the sequential MCL algorithm that we review in Section \ref{sec:MCL}.
The construction of the low-order Jacobian, our IDP stopping criterion
for fixed-point iterations, and the adaptive choice of relaxation
parameters for steady-state computations are discussed
in Section \ref{sec:Solver}. We present the results of our numerical experiments in 
Section \ref{sec:examples} and draw conclusions in Section \ref{sec:concl}.

\section{Low-order space discretization}
\label{sec:lodisc}

The Euler equations of gas dynamics are a nonlinear
hyperbolic system that can be written as
\begin{equation}\label{eq:hypsys}
	\frac{\partial u}{\partial t}+\nabla\cdot\f(u) = 0.
\end{equation}
We denote by $u = u(\mathbf{x}, t) \in\R^m$, $m=d+2$ the vector of conserved quantities at a space location $\x\in\Omega\subset\R^d,\,d\in\{1,2,3\}$ and time $t\geq0$.
The conservative variables that constitute $$u = (\rho, \rho\magenta{\vel^\top},\rho E)^\top$$ are the density, momentum, and total energy, respectively.
The flux function is defined by
\begin{equation*}
		\f(u)= (f_1(u),\ldots, f_d(u)) = \left( \begin{array}{c}
		\rho\mathbf{v}^\top\\
		\rho \mathbf{v}\otimes \mathbf{v}+pI_d\\
		(\rho E+p)\mathbf{v}^\top
	\end{array}\right)\in\R^{m\times d},
\end{equation*}
where $I_d \in\R^{d\times d}$ is the identity tensor and $p$ is the pressure.
We use the equation of state
\begin{equation}\label{eq:pres}
	p = p(u) = (\gamma - 1)\biggl(\rho E - \frac{|\rho\vel|^2}{2\rho}\biggr)
\end{equation}
of an ideal gas with adiabatic constant $\gamma>1$. \magenta{The work of Tovar et al. \cite{tovar2023}
appears to provide a good stepping stone for extending the proposed methodology to arbitrary equations of state.

Guermond and Popov \cite{guermond2016} call a set $\mathcal G\subset \R^m$ an
  \emph{invariant set} of a hyperbolic system if spatially averaged exact entropy solutions of
  Riemann problems stay in $\mathcal G$ for all $t>0$. For the system of Euler equations with the equation
  of state \eqref{eq:pres}, the convex set
\begin{equation}\label{invdom}
	\GG = \{(\rho,\rho\mathbf v^\top,\rho E)^\top: \rho > 0,\ p > 0 \}
\end{equation}
of physically admissible states
is an invariant set \cite[Sec. 5.3]{batten1997}, and so is the subset of $\GG$ such that the specific 
entropy $s=\log(p/\rho^{\gamma})$ is bounded below by a constant depending on the initial data.

A numerical
method is called invariant domain preserving (IDP) or positivity preserving
if it produces
approximations that stay in an invariant set. The foundations of modern IDP
methods were laid by Hoff \cite{hoff1979}, Frid \cite{frid2001}, Perthame and Shu \cite{perthame1996},
Bouchut \cite{bouchut2004}, and Berthon \cite{berthon2008}. 
}
\smallskip

Let us also define the Jacobians that we need to construct the matrices of our
implicit finite element schemes.
The array of $d$ flux Jacobians associated with the components of $\f(u)$ is given by
\begin{equation*}
	\A(u) = (A_1(u), \ldots, A_d(u)) = \left(\frac{\partial f_1(u)}{\partial u},\ldots, \frac{\partial f_d(u)}{\partial u}\right)  \in\R^{m\times m \times d}.
\end{equation*}
In the two-dimensional case, $\A(u)$ consists of~\cite{dolejvsi2015, gurris2009}
\begin{equation*}
	A_1(u) = \left(\begin{matrix}
		0 & 1 & 0 & 0 \\
		b_2 v_1^2 + b_1 v_2^2 & (3-\gamma)v_1 & (1-\gamma)v_2 & \gamma -1\\
		- v_1v_2 & v_2 & v_1 & 0\\
		b_1(v_1^3 + v_2^2v_1)- (E+ \frac{p}{\rho}) v_1 &  (E+ \frac{p}{\rho})-(\gamma -1)v_1^2 & (\gamma -1)v_1v_2 & \gamma v_1
	\end{matrix}\right)
\end{equation*}
and
\begin{equation*}
	A_2(u) = \left(\begin{matrix}
		0 & 0 & 1 & 0 \\
		- v_1v_2 & v_2 & v_1 & 0\\
		b_2 v_2^2 + b_1 v_1^2 & (1-\gamma)v_1 & (3-\gamma)v_2 & \gamma -1\\
		b_1( v_1^2v_2+  v_2^3)- (E+ \frac{p}{\rho}) v_2 & (\gamma -1)v_1v_2&  (E+ \frac{p}{\rho})-(\gamma -1)v_2^2 & \gamma v_2
	\end{matrix}\right),
\end{equation*}
where 
\begin{equation*}
	b_1 = \frac{\gamma -1}{2}, \quad b_2 = \frac{\gamma-3}{2}.
\end{equation*}

For any vector $\n = (n_1, \ldots, n_d)$, the directional Jacobian
\begin{equation*}
	\n\cdot\A(u) = \sum_{i = 1}^d n_i A_i(u)\in\R^{m\times m}
\end{equation*}
is diagonalizable. Its real eigenvalues are given by~\cite{dolejvsi2015, hajduk2022diss}
\begin{equation*}
	\lambda_1 =\vel\cdot \n - c,\quad\lambda_2 =\ldots =\lambda_{d+1} = \vel\cdot \n,\quad\lambda_m = \vel\cdot\n +c,
\end{equation*}
where $c = c(u) = \sqrt{\frac{\gamma p}{\rho}}$ is the local speed of sound.
The spectral radius
\begin{equation*}
	\lambda(u, \n)  = |\vel\cdot \n| + c=\mathrm{spr}(\n\cdot \A(u))
\end{equation*}
of $\n\cdot \A(u)$ determines the maximum local speed of wave propagation.

Multiplying~\eqref{eq:hypsys} by a test function $w$ and integrating
over $\Omega$, we obtain the weak form
\begin{equation}\label{eq:weak}
	\int_\Omega w\left( \frac{\partial u}{\partial t} +\nabla \cdot \f(u)\right)\,\diff\x = \int_{\Gamma}w[\f(u)\cdot \n - \mathcal F(u,\hat u;\mathbf n)]\diff \s,
\end{equation}
\red{where $\mathcal F(u,\hat u;\mathbf n)$ is a numerical approximation to the normal flux $\f(u)\cdot\n$. The external state $\hat u$ of the approximate Riemann solver is the boundary data that we prescribe on $\Gamma = \partial\Omega$ in a weak sense. In this work, we use
  the property-preserving local Lax--Friedrichs (LLF) flux
\begin{equation}\label{eq:LLF}
    \mathcal F(u_L, u_R;\mathbf n) = \dfrac{\mathbf f(u_R)+\mathbf f(u_L)}{2}\cdot \mathbf n - \dfrac{\lambda_{LR}}{2}(u_R - u_L)
\end{equation}
and define the viscosity parameter $\lambda_{LR} =\lambda_{\max}(u_L,u_R;\mathbf n)$ using the Rusanov approximation
\begin{equation}\label{eq:Rusanov}
\lambda_{\max}(u_L,u_R;\mathbf n)=\max\{\lambda(u_L, \n),\lambda(u_R, \n)\}
\end{equation}
to the maximum wave speed of the one-dimensional Riemann problem. Algorithms for
calculating a guaranteed upper bound $\lambda_{\max}^{\rm GUB}$ for this speed can be found in \cite{guermond2018,guermond2016a,guermond2016,toro2020}.
\medskip

To discretize \eqref{eq:weak} in space using a continuous Galerkin (CG) method,
  we approximate $\bar\Omega$ by the union $\bar\Omega_h=\bigcup_{e=1}^{E_h}K_e$
  of nonoverlapping cells $K_{1},\ldots,K_{E_h}$ that constitute
  a conforming mesh $\mathcal T_h$. We assume that each cell $K\in\mathcal T_h$ can be mapped to the
  same $d$-dimensional simplex or box $\hat K$. The local basis functions of a linear ($\mathbb{P}_1$) or
  multilinear ($\mathbb{Q}_1$) finite element are defined on $\hat K$.
The corresponding global
basis functions $\varphi_1,\ldots,\varphi_{N_h}$ are associated with the vertices $\x_1,\ldots,\x_{N_h}$ of $\mathcal{T}_h$ and satisfy}
\begin{equation*}
	\varphi_i(\x_j)=\delta_{ij}\quad \forall i,j\in\{1,\ldots,N_h\}.
\end{equation*}
\red{We use these Lagrange basis functions to define the continuous finite element approximations
\begin{equation}\label{eq:GFE}
	u_h = \sum_{j= 1}^{N_h} u_j \varphi_j\approx u,\qquad \f_h(u_h) = \sum_{j=1}^{N_h} \f_j\varphi_j\approx\f (u_h),
\end{equation}
where $u_j(t)=u_h(\mathbf x_j,t)$ and $\f_j(t)=\f(u_j(t))$ are the nodal values of $u_h$ and $\f(u_h)$, respectively.

Let $\mathcal{E}_i$ denote the set of indices of mesh cells containing the vertex
$\mathbf{x}_i$ and $\mathcal N^e$ the set of indices of vertices belonging to  the cell $K_e$.
Additionally, we define the nodal stencils
$\mathcal{N}_i =\bigcup_{e\in\mathcal E_i}\mathcal N^e$ and $\mathcal{N}_i^*=\mathcal{N}_i\setminus \{i\}$.
The notation
$\mathcal S_\Gamma(K_e)$ will be used for the set of external boundary faces $S\subset\partial K_e\cap\partial\Omega_h$ of $K_e$ 
and $\mathbf n_S$ for the constant unit outward normal to $S\in\mathcal S_\Gamma(K_e)$.

Substituting \eqref{eq:GFE}
into \eqref{eq:weak} and using $w=\varphi_i$ as a test function, we obtain the CG
discretization~\cite{kuzmin2020, kuzmin2023}}
\begin{equation}\label{eq:standardCG}
  \sum_{j\in \mathcal{N}_i}m_{ij} \frac{\diff u_j}{\diff t}= b_i(u_h,\hat u)-\sum_{j\in \mathcal{N}_i} \f_j\cdot\con_{ij},
\end{equation}
where
\begin{gather*}
	m_{ij} = \sum_{e\in\E_i\cap\E_j} \int_{K_e}\varphi_i\varphi_j\,\diff\x,\qquad
	\con_{ij} =  \sum_{e\in\E_i\cap\E_j} \int_{K_e}\varphi_i\nabla \varphi_j\,\diff\x,\\
	\red{b_i(u_h, \hat u) =  \sum_{e\in \mathcal{E}_i}\sum_{S\in\mathcal S_\Gamma(K_e)}\int_S
        \varphi_i[\f(u_h)\cdot\n -\mathcal F(u,\hat u;\mathbf n)]\,\diff\s.}
\end{gather*}
\red{
\begin{remark}
  The global conservation property of the CG method can be shown as usual
  by using the admissible test function $w\equiv 1$ in \eqref{eq:weak} or
  summing equations \eqref{eq:standardCG} over $i=1,\ldots,N_h$.
\end{remark}  

To derive a property-preserving low-order approximation of local Lax--Friedrichs type
  as in \cite{kuzmin2012b,kuzmin2010a},} we first approximate the left-hand side of \eqref{eq:standardCG}
  by $m_i\frac{\diff u_i}{\diff t}$, where
$$m_i = \sum_{j=1}^{N_h}m_{ij} =\sum_{e\in\E_i} \int_{K_e}\varphi_i\diff\x$$ is a diagonal
entry of the lumped mass matrix. On the right-hand side,
we replace $b_i(u_h, \hat u)$ by\red{
\begin{align}\notag
  \tilde b_i(u_i, \hat u) &=  \sum_{e\in \mathcal{E}_i}\sum_{S\in\mathcal S_\Gamma(K_e)}\int_S
  \varphi_i[\f_i\cdot \n-\mathcal F(u_i,\hat u_{i,S};\mathbf n)]\,\diff\s\\
    &= \sum_{e\in \mathcal{E}_i}\sum_{S\in\mathcal S_\Gamma(K_e)}
    \frac{\sigma_{i,S}}2
    \left[\lambda_{i,S}(\hat u_{i,S}-u_i)-
      (\f(\hat u_{i,S})-\f_i)\cdot\mathbf n_{S}\right],\label{eq:bdrtermLO}
  \end{align}
where
$$
  \sigma_{i,S}=\int_{S}\varphi_i\,\diff\s,\qquad
  \hat u_{i,S}=\frac{1}{\sigma_{i,S}}\int_{S}\varphi_i\hat u\,\diff\s,\qquad
 \lambda_{i,S}=\lambda_{\max}(u_i,\hat u_{i,S};\mathbf n_S).
 $$
 Similarly to mass lumping, this approximation can be interpreted as inexact numerical integration using a nodal quadrature rule.
 We use it because it makes the discretization of boundary terms IDP (see below), while preserving
 the second-order accuracy of the baseline discretization.
}

Owing to the partition of unity property  $\sum_i\varphi_i \equiv 1$ of the Lagrange basis functions, we have
\begin{equation*}
	\sum_{j\in \mathcal{N}_i} \f_j\cdot\con_{ij} = \sum_{j\in\mathcal{N}_i^*} (\f_j-\f_i )\cdot\mathbf{c}_{ij}.
\end{equation*}
\red{Thus the lumped counterpart of the spatial semi-discretization \eqref{eq:standardCG} can be written as
\begin{equation}\label{eq:lumpedCG}
  m_i \frac{\diff u_i}{\diff t}= \tilde b_i(u_h,\hat u)-\sum_{j\in\mathcal{N}_i^*} (\f_j-\f_i )\cdot\mathbf{c}_{ij}.
\end{equation}
This discretization is equivalent to a vertex-centered finite volume scheme with centered numerical
fluxes \cite{selmin1996,selmin1993}. A low-order LLF approximation is defined by~\cite{guermond2018,guermond2016,kuzmin2020,kuzmin2023}
\begin{equation}\label{eq:LLFdef}
		m_i \frac{\diff u_i}{\diff t} 
		=\tilde{b}_i(u_i, \hat u)+ \sum_{j\in\mathcal{N}_i^*} [d_{ij}(u_j-u_i) - (\mathbf{f}_j - \mathbf{f}_i )\cdot\mathbf{c}_{ij}].
\end{equation}
Following Kuzmin et al. \cite{kuzmin2010a}, we define the \emph{graph viscosity} coefficients
\begin{equation}\label{graph}
	d_{ij} = \begin{cases}
		\max\{\lambda_{ij}|\con_{ij}|,\lambda_{ji}|\con_{ji}|\} &\text{if }j\in\mathcal{N}_i^*,\\
		-\sum_{k\in\mathcal{N}_i^*} d_{ik} & \text{if }j = i,\\
		0 & \text{otherwise}
	\end{cases}
\end{equation}
 using the Rusanov wave speed
$\lambda_{ij}=\lambda_{\max}(u_i,u_j;\mathbf n_{ij})$ corresponding to the unit
normal $\n_{ij}=\frac{\con_{ij}}{|\con_{ij}|}$.

The LLF evolution equation \eqref{eq:LLFdef} for $u_i$ can be written in the \emph{bar state} form\footnote{\red{In the literature on
    finite volume methods for conservation laws, such a representation in terms of jumps is known as the \emph{fluctuation form} of an approximate Riemann solver \cite{leveque2002}.}}
  (cf. \cite{guermond2016,kuzmin2020})
\begin{equation}\label{eq:LOdisc}
		m_i \frac{\diff u_i}{\diff t} 
		=\sum_{e\in \mathcal{E}_i}\sum_{S\in\mathcal S_\Gamma(K_e)}
                \sigma_{i,S}\lambda_{i,S}(\bar u_{i,S}-u_i)
                + \sum_{j\in\mathcal{N}_i^*} 2d_{ij}(\overline{u}_{ij}-u_i),
\end{equation}
where
\begin{equation*}
   \bar u_{i,S} = \frac{\hat u_{i,S} + u_i}{2}
  - \frac{(\mathbf f(\hat u_{i,S})- \mathbf f_i)\cdot \mathbf n_{S}}{2\lambda_{i,S}},\qquad
	\overline{u}_{ij} = \frac{u_j + u_i}{2}-\frac{(\f_j-\f_i)\cdot\mathbf{c}_{ij}}{2 d_{ij}}
\end{equation*}
are intermediate states of the HLL Riemann solver \cite{harten1983b}.
\magenta{If $\lambda_{i,S}$ and
  $d_{ij}$ are defined using a guaranteed upper bound (GUB) for the maximum local wave speed (as in \cite{guermond2016}),
    then  $\bar u_{i,S}$ and $\overline{u}_{ij} $
 coincide with spatial averages of exact solutions to the corresponding Riemann problems. It follows that
 \begin{equation}\label{eq:idpuijbar}
 u_i,\hat u_{i,S}\in\GG \quad\Rightarrow\quad\overline{u}_{i,S}\in \GG,\qquad  
 u_i,u_j\in\GG \quad\Rightarrow\quad\overline{u}_{ij}\in \GG
 \end{equation}
 for any convex invariant set $\mathcal G$. The use of the Rusanov wave speed
 \eqref{eq:Rusanov} in our implementation ensures the validity of \eqref{eq:idpuijbar}
 for $\mathcal G$ defined by \eqref{invdom}, i.e., we have positivity preservation
 for the density and pressure. This property was shown by Lin et al. \cite{lin2023}
 using Zhang’s lemma \cite[App. B]{zhang2017} and by Batten et al. \cite[Sec. 5.3]{batten1997}
 for intermediate states of a general HLLC Riemann solver. The GUB version of the
 LLF method additionally guarantees that a minimum principle holds for the specific physical entropy
 \cite{guermond2016} but, in our experience, there is no practical benefit in
 using $\lambda_{\max}^{\rm GUB}$ for the Euler equations.}
\medskip

Using the  forward Euler method to discretize \eqref{eq:LOdisc} in time, we arrive at
$$
u_i^{n+1} = u_i^n +\frac{\Delta t}{m_i}\left[\sum_{e\in \mathcal{E}_i}\sum_{S\in\mathcal S_\Gamma(K_e)}
                \sigma_{i,S}\lambda_{i,S}(\bar u_{i,S}^n-u_i^n)+
\sum_{j\in\mathcal{N}_i^*}2d_{ij}(\overline{u}_{ij}^{n}-u_i^{n})\right].
$$
If conditions \eqref{eq:idpuijbar} hold and 
the time step $\Delta t$ satisfies the CFL-like condition
\begin{equation}\label{eq:cfl}
  \frac{\Delta t}{m_i}\left[\sum_{e\in \mathcal{E}_i}\sum_{S\in\mathcal S_\Gamma(K_e)}
                \sigma_{i,S}\lambda_{i,S}+
    \sum_{j\in\mathcal{N}_i^*}2d_{ij}\right]\le 1,
\end{equation}
then $u_i^{n+1}$ is a convex combination of states belonging
to the invariant set $\mathcal G$. Hence, the combination of
the LLF method \eqref{eq:LOdisc} with forward Euler time stepping
is IDP under condition \eqref{eq:cfl}. 
}

\section{Backward Euler time discretization}
\label{sec:IDP}

\red{ In this work, our main
 interest is in steady-state solutions, which we calculate using
 a flux-corrected high-order extension of \eqref{eq:LOdisc}  and
pseudo-time stepping of backward Euler type. The implicit approach} allows us
to use large pseudo-time steps, but we need to show the existence of
an IDP solution to the resulting nonlinear system and design a
robust iterative solver. We begin with the former task and prove
the IDP property of the fully implicit
low-order scheme in this section.

For simplicity, we assume that periodic boundary conditions are prescribed on $\Gamma$.
In this case, $\tilde{b}_i(u_i, \hat u)=0$ for $i=1,\ldots,N_h$, and the
backward Euler time discretization of~\eqref{eq:LOdisc} yields
\begin{equation}\label{eq:implicitEuler}
  u_i^{n+1} = u_i^n + \frac{\Delta t}{m_i}
  \sum_{j\in\mathcal{N}_i^*}2d_{ij}(\overline{u}_{ij}^{n+1}-u_i^{n+1}),\qquad i=1,\ldots,N_h,
\end{equation}
where $d_{ij} = d_{ij}(u^{n+1})$ \red{are the solution-dependent LLF graph viscosity coefficients defined by \eqref{graph}.}

Let $G:=\{u\in\R^{mN_h}\,:\,\magenta{(u_{(i-1)m+1},\ldots,u_{im})^\top}\in\mathcal G,\
i=1,\ldots,N_h\}$, where $\mathcal G$ is the invariant \cyan{set} defined by \eqref{invdom}. \red{Any conservative approximation $u\in G$
  preserves the mesh-dependent global bounds
$$
  \rho_h^{\max}=\frac{\int_{\Omega}\rho_h^0\dx}{\min_{1\le i\le N_h}m_i},\qquad
(\rho E)_h^{\max}=\frac{\int_{\Omega}(\rho E)_h^0\dx}{\min_{1\le i\le N_h}m_i}
$$
for the nodal values of the density and total energy.
It follows that the internal energy, pressure, and velocity
are bounded as well. Hence, approximations produced by conservative
IDP schemes stay in the bounded subset
$G_h:=\{u\in\R^{mN_h}\,:\,\magenta{(u_{(i-1)m+1},\ldots,u_{im})^\top}
\in \mathcal G_h,\
i=1,\ldots,N_h\}\subset G$, where
\begin{equation}\label{eq:globb}
\mathcal G_h=\{(\rho,\rho\mathbf v^\top,\rho E)^\top\,:\,
0< \rho \le\rho_h^{\max},\ p> 0,\ \rho E\le (\rho E)_h^{\max}\}\subset\mathcal G.
\end{equation}
Moreover, the maximum speed of wave propagation is bounded above by a global constant $\lambda_h^{\max}$.}

\blue{
To prove that the nonlinear system \eqref{eq:implicitEuler} has
an IDP solution $u=u^{n+1}$, we will construct a mapping
$\Psi:G\to G$ to which we can apply the following fixed-point
theorem \cite[Theorem 2]{burton1998}.

\begin{theorem}[A fixed-point theorem of Krasnoselskii]\label{thm1}
  Let $M$ be a closed, convex, and nonempty subset of a Banach space
  $(S,\|\cdot\|)$. Suppose that $A:M\to S$ and $B:S\to S$ are such that
  \begin{itemize}
  \item [(i)] $B$ is a contraction with constant $\alpha<1$,
  \item [(ii)] $A$ is continuous, $AM$ is contained in a compact
    subset of $S$,
  \item [(iii)] $[x=Bx+Ay,\ y\in M] \ \Rightarrow\ x\in M$.
  \end{itemize}
  Then there exists $y\in M$ such that $y=Ay+By$. 
  \end{theorem}

\begin{proof}
See \cite[Theorem 2]{burton1998}.
\end{proof}    
}

\magenta{ 
Following the analysis of implicit schemes in~\cite{chalabi1997, tang2000} and adapting
it to \eqref{eq:implicitEuler}, we introduce a free parameter $s>0$ and define the local mapping
  \begin{equation*}
	\Psi^s_i(u) = \frac{s u_i^n}{1+ s}+ \frac{1}{1+s}\biggl(u_i+ \frac{s\Delta t}{m_i}\sum_{j\in\mathcal{N}^*_i}2d_{ij}(\overline{u}_{ij}-u_i)\biggr)=\frac{s u_i^n}{1+ s}+ \frac{\gamma_i(u)}{1+s},
\end{equation*}
where
\begin{equation*}
	\gamma_i(u) = \biggl(1-\frac{s\Delta t}{m_i}\sum_{j\in\mathcal{N}^*_i}2d_{ij}\biggr)u_i+ \frac{s\Delta t}{m_i}\sum_{j\in\mathcal{N}^*_i}2d_{ij}\overline{u}_{ij}.
\end{equation*}
Using this definition of $\Psi^s_i(u)\in\mathbb{R}^m$, we assemble the global vector} 
$$\magenta{\Psi^s(u)=\{v\in\R^{mN_h}\,:\,(v_{(i-1)m+1},\ldots,v_{im})^\top
  =\Psi_i^s(u),\  i=1,\ldots,N_h\}.}
  $$
Note that $\Psi^s_i(u)$ is a convex combination of $u_i^n$ and $\gamma_i(u)$ \red{and that
  $\Psi^s_i(u^{n+1})=u_i^{n+1}$
  for a solution $u=u_i^{n+1}$ of the nonlinear system \eqref{eq:implicitEuler}.}
The assumption that $u_i^n\in\mathcal G$ implies that
$\Phi(u_i^n)\ge 0$ for $\Phi\in\{\rho,p\}$.
For any concave nonnegative function $\Phi$, we have the estimate
\begin{equation*}
		\Phi(\Psi^s_i(u))\geq \frac{s }{1+ s}\Phi(u_i^n) + \frac{1}{1+s}\Phi(\gamma_i(u))\geq \frac{1}{1+s}\Phi(\gamma_i(u)).
\end{equation*}
Furthermore, $\gamma_i(u)$ is a convex combination of $u_i$ and $\overline{u}_{ij}$ under the CFL condition 
\begin{equation}\label{eq:cflcond}
	\frac{s\Delta t}{m_i}\sum_{j\in\mathcal{N}^*_i}2d_{ij}\leq 1.
\end{equation}
This implies that
\begin{align*}
  \Phi(\gamma_i(u))
  \geq \biggl(1-\frac{s\Delta t}{m_i}\sum_{j\in\mathcal{N}^*_i}2d_{ij}\biggr)\Phi(u_i) + \frac{s\Delta t}{m_i}\sum_{j\in\mathcal{N}^*_i}2d_{ij}\Phi(\overline{u}_{ij}),
\end{align*}
\red{where $\Phi(\overline{u}_{ij})\ge 0$  for $\Phi\in\{\rho,p\}$ because
  of the local IDP property \eqref{eq:idpuijbar} of the LLF bar states
  $\bar u_{ij}$.} From
this we deduce that $\Phi( \Psi^s_i(u))\ge 0$ for $\Phi\in\{\rho,p\}$
and, therefore, $\Psi^s(u) \in G$ for all $u\in G$.
\medskip

Let us now set $u^{(0)}= u^n$ and use the mapping $\Psi^s:G\to G$
to construct the fixed-point iteration 
\begin{equation}\label{eq:fixpIDP}
	u^{(k+1)}= \Psi^s(u^{(k)}), \quad k\in\N_0.
\end{equation}
\blue{
The following application of Theorem \ref{thm1} guarantees the
existence of a bounded IDP solution to \eqref{eq:implicitEuler}.
\begin{theorem}\label{thm2}
  Suppose that $u^n\in G_h$. Then the nonlinear system
  \eqref{eq:implicitEuler} has a solution $u^{n+1}\in G_h$.
\end{theorem}

\begin{proof}
By definition, subvectors $u_i\in\mathbb{R}^m$ of a vector $u\in G_h
\subset G$ are contained in the set $\mathcal G_h\subset \mathcal G$
defined by \eqref{eq:globb}. The set $M:=G_h$ is convex,
while $\bar G_h$ is compact. For $j\in\mathcal N_i^*$, the graph viscosity coefficients $d_{ij}(u)=d_{ji}(u)$ are bounded above by $d_{ij}^{\max}
=\max\{|\mathbf c_{ij}|,|\mathbf c_{ji}|\}\lambda_h^{\max}=d_{ji}^{\max}$.
The global bound $\lambda_h^{\max}$ for the maximum speed depends on
the mesh but is independent of $u\in G_h$.
The mapping $\Psi^s(u^{(k)})$ of the fixed-point iteration
\eqref{eq:fixpIDP} admits the Krasnoselskii splitting 
 $$	\Psi^s_i(u)=(Au)_i+(Bu)_i,$$
 where
 \begin{align*}
 (Au)_i&=\frac{s\Delta t}{(1+s)m_i}
 \sum_{j\in\mathcal{N}^*_i}2[d_{ij}\overline{u}_{ij}+(d_{ij}^{\max}-d_{ij})u_i],\\
 (Bu)_i&=\frac{s u_i^n}{1+s}+\frac{1}{1+s}\Big(1-
 \frac{s\Delta t}{m_i}\sum_{j\in\mathcal{N}^*_i}2d_{ij}^{\max}\Big)u_i.
 \end{align*}
The mapping $A:M\to \mathbb{R}^{mN_h}$ is continuous. If $s>0$ is chosen as small as necessary to satisfy
 $$
\frac{s\Delta t}{m_i}\sum_{j\in\mathcal{N}^*_i}2d_{ij}^{\max}\le 1,
$$
then $AM$ is contained in the compact set $\bar G_h$ and 
$B:M\to\mathbb{R}^{mN_h}$ is a contraction. Hence, the
assumptions (i) and (ii) of Theorem \ref{thm1} are met.
The validity of the assumption (iii) follows from the
fact that if $u^n,y\in M$ and $x=Bx+Ay$, then for any
$i\in \{1,\ldots,N_h\}$ the subvector
$$
 x_i=\frac{1}{1+
 \frac{\Delta t}{m_i}\sum_{k\in\mathcal{N}^*_i}2d_{ik}^{\max}} 
\Big(u_i^n+\frac{\Delta t}{m_i}
  \sum_{j\in\mathcal{N}^*_i}2\left[d_{ij}\overline{y}_{ij}
+(d_{ij}^{\max}-d_{ij})y_i\right]\Big)
$$
is a convex combination of $u_i^n,y_i,\overline{y}_{ij}\in\mathcal G$. Applying Theorem \ref{thm1},
we complete the proof.
\end{proof}

In principle, a similar fixed-point argument can be used to show the existence of a steady-state IDP solution to problem \eqref{eq:LOdisc} with non-periodic boundary conditions, i.e., of $u\in G$ such that
$$
0=\sum_{e\in \mathcal{E}_i}\sum_{S\in\mathcal S_\Gamma(K_e)}
                \sigma_{i,S}\lambda_{i,S}(\bar u_{i,S}-u_i)
                + \sum_{j\in\mathcal{N}_i^*} 2d_{ij}(\overline{u}_{ij}-u_i)
               ,\qquad i=1,\ldots,N_h.
$$
Indeed, this nonlinear system can be written as
\begin{equation}
  \label{eq:nonlinsyssteady}
u_i= \bar u_i
+\frac{1}{\mu_i}
\sum_{j\in\mathcal{N}_i^*} 2d_{ij}(\overline{u}_{ij}-u_i)
,\qquad i=1,\ldots,N_h,
\end{equation}
where
$$
\mu_i=\frac{m_i}{\Delta t}+
\sum_{e\in \mathcal{E}_i}\sum_{S\in\mathcal S_\Gamma(K_e)}
\sigma_{i,S}\lambda_{i,S},\qquad \bar u_i=\frac{1}{\mu_i}
\left[\frac{m_i}{\Delta t}u_i+
\sum_{e\in \mathcal{E}_i}\sum_{S\in\mathcal S_\Gamma(K_e)}
                \sigma_{i,S}\lambda_{i,S}\bar u_{i,S}\right].
                $$
                Note that system \eqref{eq:nonlinsyssteady} has the same structure as
                \eqref{eq:implicitEuler}. Assuming that the iterates stay in a bounded set $G_h\subset G$, existence of a solution $u\in G$ can be shown following the proof of Theorem \ref{thm2}.

\begin{remark}                          
The difficulty of verifying the assumption (i) of Theorem \ref{thm1} for \eqref{eq:nonlinsyssteady} depends on the particular choice of boundary conditions. The contractivity of the steady-state mapping
                  $$
                (Bu)_i=\frac{s \bar u_i}{1+s}+\frac{1}{1+s}\Big(1-
 \frac{s}{\mu_i}\sum_{j\in\mathcal{N}^*_i}2d_{ij}^{\max}\Big)u_i
 $$
 can easily be shown, e.g., for $\Gamma=\bar\Gamma_{\rm in}\cup \bar\Gamma_{\rm out}$,
 where $\Gamma_{\rm in}$ is a supersonic inlet and $\Gamma_{\rm out}$ is a supersonic outlet. In this case, the external state $\bar u_i$ is independent of $u_i$ on $\Gamma_{\rm in}$ and $\bar u_i=u_i$ on $\Gamma_{\rm out}$. Thus
 $$
 |(Bu)_i-(Bv)_i|\le \frac{1}{1+s}\Big(
1+s-\frac{s}{\mu_i}\sum_{j\in\mathcal{N}^*_i}2d_{ij}^{\max}
 \Big)|u_i-v_i|
 $$
 for any vectors $u,v\in G$ such that $u_i=\bar u_i=v_i$ on $\Gamma_{\rm in}$.
 The mapping $B$ is contractive
 if the parameter $s>0$ is chosen small enough to satisfy the CFL-like condition
 $ \frac{s}{\mu_i}\sum_{j\in\mathcal{N}^*_i}2d_{ij}^{\max}\le 1$.
\end{remark}

\begin{remark}
  In the existence proof for \eqref{eq:nonlinsyssteady}, the global bounds of the sets $\mathcal G_h$ and $G_h$ can no longer be defined in terms of initial states. It is intuitively clear that the bounds for (approximations to)  steady-state solutions should depend on the boundary data. Instead of imposing an upper bound on the total energy, one can define 
$
  \mathcal G_h=\{(\rho,\rho\mathbf v^\top,\rho E)^\top\,:\,
  0< \rho \le\rho_h^{\max},\ p> 0,\ \eta \le \eta_h^{\max}\}\subset\mathcal G
  $
using an upper bound $\eta_h^{\max}$ for the mathematical entropy $\eta(u)=\frac{\rho\log(p/\rho^{\gamma})}{1-\gamma}$. In view of the fact that
  $$
  \pd{\rho}{t}+\nabla\cdot(\rho\mathbf v)=0,\qquad
  \pd{\eta}{t}+\nabla\cdot(\eta\mathbf v)\le 0,
  $$
  the total mass and entropy should depend on the values of $\rho$ and $\eta$ on
  $\Gamma_-=\{\mathbf x\in\Gamma\,:\,\mathbf v(\mathbf x)\cdot
  \mathbf n(\mathbf x)<0\}$. Therefore, we conjecture that
  $ \rho_h^{\max}=\frac{|\Omega|\sup_{\mathbf x\in\Gamma_-}
    \bar\rho(\mathbf x)}{\min_{1\le i\le N_h}m_i}$ and
  $\eta_h^{\max}=\frac{|\Omega|\sup_{\mathbf x\in\Gamma_-}\bar\eta(\mathbf x)}{\min_{1\le i\le N_h}m_i}$
  is an appropriate definition of the global bounds to be preserved by the
  mapping $A$ of the Krasnoselskii theorem\footnote{\red{More detailed analysis is difficult
  because even the behavior of exact steady-state solutions to boundary value problems for nonlinear hyperbolic problems is not well understood.}}. 
\end{remark}
}

\section{Monolithic convex limiting}
\label{sec:MCL}

By the Godunov theorem, the spatial semi-discretization \eqref{eq:implicitEuler} is at most first-order
accurate. Obviously, the accuracy of steady-state solutions is also affected by this
order barrier. The monolithic convex limiting (MCL) strategy proposed
in~\cite{kuzmin2020, kuzmin2023} makes it possible to achieve higher
resolution without losing the IDP property and inhibiting convergence
to steady-state solutions. To explain the MCL design philosophy,
we write the \red{lumped Galerkin} discretization \eqref{eq:standardCG} in the
 bar state form
\begin{equation}\label{eq:HOdiscBar}
  m_i \frac{\diff u_i}{\diff t} 
  =\red{\sum_{e\in \mathcal{E}_i}\sum_{S\in\mathcal S_\Gamma(K_e)}
  \sigma_{i,S}\lambda_{i,S}(\bar u_{i,S}-u_i)}
  + \sum_{j\in\mathcal{N}_i^*} 2d_{ij}(\overline{u}_{ij}^H-u_i),
\end{equation}
which exhibits the same structure as \eqref{eq:LOdisc}. \red{The
LLF component $\overline{u}_{ij}$ of the high-order bar states
\begin{equation*}
	\overline{u}_{ij}^H = \overline{u}_{ij} +\frac{f_{ij}^H}{2d_{ij}}
\end{equation*}
is IDP but violations of global and/or local bounds may be caused by the
built-in \emph{antidiffusive flux}}
\begin{equation}\label{eq:fluxtarget}
	\blue{f_{ij}^H =d_{ij}(u_i-u_j).}
\end{equation}

\blue{
  The mass lumping error on the left-hand side of \eqref{eq:HOdiscBar} can be corrected by decomposing this error into fluxes $m_{ij} \left( \dot{u}_{i} - \dot{u}_{j}\right)$ and incorporating them into 
\begin{equation}
	f_{ij} = m_{ij} \left( \dot{u}_{i} - \dot{u}_{j} \right) + d_{ij}(u_i-u_j).
\end{equation}
Following Kuzmin \cite{kuzmin2020}, we define $f_{ij}$
using the LLF time derivatives (cf. \eqref{eq:LOdisc})
$$ \dot{u}_{i} =
  \frac{1}{m_i}\left[\sum_{e\in \mathcal{E}_i}\sum_{S\in\mathcal S_\Gamma(K_e)}
  \sigma_{i,S}\lambda_{i,S}(\bar u_{i,S}-u_i)
  + \sum_{j\in\mathcal{N}_i^*} 2d_{ij}(\overline{u}_{ij}-u_i)\right].
$$
The inclusion of $m_{ij} \left( \dot{u}_{i} - \dot{u}_{j}\right)$ preserves
second-order accuracy and has a stabilizing effect \cite{kuzmin2020,kuzmin2023}.
}
\medskip

The MCL approach introduced in~\cite{kuzmin2020} approximates an unacceptable \emph{target flux}
$f_{ij}=-f_{ji}$
  by a limited flux $f_{ij}^*=-f_{ji}^*$ that has the local IDP property (cf. \eqref{eq:idpuijbar})
\begin{equation}\label{eq:global}
  \overline{u}_{ij},  \overline{u}_{ji}\in\GG\quad\Rightarrow\quad
	\overline{u}_{ij}^* := \overline{u}_{ij} + \frac{f_{ij}^* }{2d_{ij}} \in\GG,\qquad
  \overline{u}_{ji}^* := \overline{u}_{ji} - \frac{f_{ij}^* }{2d_{ij}} \in\GG.
\end{equation}
\red{The admissible set $\mathcal G$ is defined by \eqref{invdom}. Since positivity
  preservation for the density and pressure (internal energy) does not rule out
  spurious oscillations, we additionally impose local discrete maximum
  principles on some scalar quantities of interest. In the remainder
  of this section, we formulate the corresponding inequality constraints and
  review algorithms for constructing a bound-preserving approximation 
  $f_{ij}^*=(f_{ij}^{*,\rho},f_{ij}^{*,\rho v_1},\ldots,f_{ij}^{*,\rho v_d},f_{ij}^{*,\rho E})^\top$
  to $f_{ij}=(f_{ij}^{\rho},f_{ij}^{\rho v_1},\ldots,f_{ij}^{\rho v_d},f_{ij}^{\rho E})^\top$.}

\medskip

In this work, we use a sequential MCL algorithm that imposes local bounds
on the density, velocity, and
specific total energy (as in \cite{kuzmin2020,kuzmin2023}). The local
bounds of the density constraints
\begin{equation}\label{eq:rhoconstraints}
  \rho_{i}^{\min} \leq \overline{\rho}_{ij}^*=
  \overline{\rho}_{ij}-\frac{f_{ij}^{*,\rho}}{2d_{ij}}
  \leq	\rho_{i}^{\max}\qquad
  \forall j\in\mathcal N_i^*
\end{equation}
are given by
\begin{equation*}
		\rho_{i}^{\min} = \min\left\{  \min_{j\in \mathcal{N}_i} \rho_j   , \min_{j\in \mathcal{N}_i^*} \overline{\rho}_{ij}  \right\},\qquad
		\rho_{i}^{\max} = \max\left\{  \max_{j\in \mathcal{N}_i} \rho_j   , \max_{j\in \mathcal{N}_i^*} \overline{\rho}_{ij}  \right\}.
\end{equation*}
The density constraints \eqref{eq:rhoconstraints} are satisfied by the limited antidiffusive fluxes
\cite{kuzmin2020, kuzmin2023} 
\begin{equation*}
	f_{ij}^{*,\rho} = \begin{cases}
		\min\{   f_{ij}^{\rho},\, 2d_{ij} \min\{ \rho_i^{\max} - \overline{\rho}_{ij},\, \overline{\rho}_{ji} - \rho_j^{\min}     \}     \}\quad & \text{if } f_{ij}^{\rho}>0,\\
		\max\{   f_{ij}^{\rho},\, 2d_{ij} \max\{ \rho_i^{\min} - \overline{\rho}_{ij},\, \overline{\rho}_{ji} - \rho_j^{\max}     \}     \}\quad &\text{otherwise.}
	\end{cases}
\end{equation*}

The components of $\mathbf v=(v_1,\ldots,v_d)$ and the specific total energy
$E$ are derived quantities that cannot be limited directly \red{in continuous finite element methods (using limited
  reconstructions of $\mathbf v$ and $E$ would violate the continuity requirement).
  Therefore, we will limit conserved
products $\rho\phi$, where $\phi\in\{ v_1,\ldots,v_d, E \}$,
using a discrete version of the product rule. Introducing the states~\cite{kuzmin2023, kuzmin2020}
\begin{equation*}
	\phi_i = \frac{(\rho\phi)_i}{\rho_i}, \quad \overline{\phi}_{ij} = \frac{\overline{(\rho\phi)}_{ij} + \overline{(\rho\phi)}_{ji}}  {\overline\rho_{ij} + \overline\rho_{ji}},
\end{equation*}
we impose the inequality constraints
\begin{equation}\label{eq:phiconstraints}
	\overline{\rho}_{ij}^* \phi_i^{\min} \leq \overline{(\rho\phi)}_{ij}^* = \overline{(\rho\phi)}_{ij} + \frac{f_{ij}^{\rho\phi,*}}{2d_{ij}} \leq \overline{\rho}_{ij}^* \phi_i^{\max}
\end{equation}
with local bounds
\begin{equation}\label{eq:philocalbounds}
	\phi_i^{\min} = \min\left\{\min_{j\in\mathcal{N}_i} \phi_j, \min_{j\in\mathcal{N}_i^*} \overline{\phi}_{ij}  \right\}, \quad
	\phi_i^{\max} = \max\left\{\max_{j\in\mathcal{N}_i} \phi_j, \max_{j\in\mathcal{N}_i^*} \overline{\phi}_{ij}  \right\}.
\end{equation}
To show that the constraints \eqref{eq:phiconstraints} are feasible, we define 
\begin{equation}\label{eq:limphiflux}
	f_{ij}^{\rho\phi,*} = g_{ij}^{\rho\phi,*} - 2 d_{ij}\left[ \overline{(\rho\phi)}_{ij}-\overline\rho_{ij}^*   \overline{\phi}_{ij}  \right]
\end{equation}
using a limited counterparts  $g_{ij}^{\rho\phi,*}$ of the auxiliary flux
\begin{equation*}
	g_{ij}^{\rho\phi} = f_{ij}^{\rho\phi} + 2 d_{ij}\left[ \overline{(\rho\phi)}_{ij}-\overline\rho_{ij}^*   \overline{\phi}_{ij}  \right].
\end{equation*}
Note that the target flux $f_{ij}^{\rho\phi}$ can be recovered using \eqref{eq:limphiflux} with
 $g_{ij}^{\rho\phi,*} = g_{ij}^{\rho\phi}$. The general definition
of $g_{ij}^{\rho\phi,*}$ should ensure that $f_{ij}^{\rho\phi,*}$ defined by \eqref{eq:limphiflux}
satisfies the local discrete maximum principle \eqref{eq:phiconstraints}.}

The choice $g_{ij}^{\rho\phi,*} = 0$ yields the low-order bar state $\overline{(\rho\phi)}_{ij} = \overline{\rho}_{ij}^* \overline{\phi}_{ij}$, which satisfies \eqref{eq:phiconstraints} by \eqref{eq:philocalbounds}. It is easy to verify that the constraints are also satisfied for
\begin{equation}\label{eq:glimited}
	g_{ij}^{\rho\phi,*} =  \begin{cases}
		\min\left\{  g_{ij}^{\rho\phi}, g_{ij}^{\rho\phi,\max} \right\} & \text{if }	g_{ij}^{\rho\phi} >0,\\
		\max\left\{  g_{ij}^{\rho\phi}, g_{ij}^{\rho\phi,\min}\right\}  & \text{if }	g_{ij}^{\rho\phi} \leq0,\\
	\end{cases}
\end{equation}
where we use the \emph{bounding fluxes}
\begin{equation*}
	\begin{split}
		g_{ij}^{\rho\phi,\max} &= \min\left\{2d_{ij}\overline\rho_{ij}^* (\phi_i^{\max} - \overline\phi_{ij} ),  2d_{ji}\overline\rho_{ji}^* (\overline\phi_{ji} - \phi_i^{\min}) \right\},\\
		g_{ij}^{\rho\phi,\min} &=\max\left\{2d_{ij}\overline\rho_{ij}^* (\phi_i^{\min} - \overline\phi_{ij} ),  2d_{ji}\overline\rho_{ji}^* (\overline\phi_{ji} - \phi_i^{\max}) \right\}.
	\end{split}
\end{equation*}
The corresponding limited flux $f_{ij}^{\rho\phi,*}$ is obtained by substituting \eqref{eq:glimited} into \eqref{eq:limphiflux}.
\medskip

In the final step of the sequential limiting procedure, we use a scalar correction
factor $\alpha_{ij}\in[0,1]$ to ensure the IDP property of the bar states
$$
\bar u_{ij}^{*,\alpha}=\bar u_{ij}+\frac{\alpha_{ij}f_{ij}^*}{2d_{ij}},
$$
where $f_{ij}^* = \left(f_{ij}^{\rho,*}, \f_{ij}^{\rho\vel,*},f_{ij}^{\rho E,*} \right)^\top$
is the vector of prelimited fluxes. Note that $\bar \rho_{ij}^{*,\alpha}\ge 0$ for
any $\alpha_{ij}\in[0,1]$ because the lower bounds of the density constraints
\eqref{eq:rhoconstraints} are nonnegative.
\red{In other words, positivity preservation
  for $\rho$ is already guaranteed.
  It remains to enforce the} pressure constraint
\begin{equation*}
  p(\overline u_{ij}^{*,\alpha})\geq 0\quad\Leftrightarrow\quad
\overline{(\rho E)}^{*,\alpha}_{ij}\geq \frac{|\overline{(\rho\vel)}^{*,\alpha}_{ij}|^2}{2\overline{\rho}^{*,\alpha}_{ij}}.
\end{equation*}
\red{This constraint is satisfied for all} $\alpha_{ij}\in[0,1]$ such that
\begin{equation}\label{eq:posfix}
	P_{ij}(\alpha_{ij})\leq Q_{ij},
\end{equation}
where
\begin{equation*}
		P_{ij}(\alpha) =\left[ \frac{| \f_{ij}^{\rho\vel,*} |^2}{2} - f_{ij}^{\rho E,*} f_{ij}^{\rho,*} \right]\alpha^2
		+ 2 d_{ij}\left[ \overline{(\rho\vel)}_{ij}\cdot\f_{ij}^{\rho\vel,*} -\overline{\rho}_{ij}  f_{ij}^{\rho E,*} - \overline{(\rho E)}_{ij} f_{ij}^{\rho,*}  \right]\alpha
\end{equation*}
and
\begin{equation*}
	Q_{ij} = (2d_{ij})^2\overline{\rho}_{ij} \left[  \overline{(\rho E)}_{ij}  - \frac{ |\overline{(\rho\vel)}_{ij}|^2}{2\overline{\rho}_{ij}}\right].
\end{equation*}
\red{Note that $Q_{ij}\geq 0$ because $\overline{u}_{ij}\in\GG$ and, therefore, $p(\overline{u}_{ij})\ge 0$. It
  follows that the pressure constraint is feasible and
  condition \eqref{eq:posfix} holds trivially for
  $\alpha_{ij}=0$.}  To find a nontrivial
 $\alpha\in[0,1]$ such that $P_{ij}(\alpha)\leq Q_{ij}$, the pressure fix
proposed in \cite{kuzmin2020} replaces \eqref{eq:posfix} by
the linear sufficient\footnote{The validity of \eqref{eq:posfix}
 is implied by \eqref{eq:posfixsuff}  because
 $P_{ij}(\alpha)\leq \alpha R_{ij}$ for any $\alpha\in[0,1]$. The
 derivation of this estimate exploits the fact that $\alpha^2\le\alpha$
for $\alpha\in[0,1]$.}  condition
\begin{equation}\label{eq:posfixsuff}
\alpha_{ij}R_{ij}\le Q_{ij},
\end{equation}
where \begin{equation*}
	\begin{split}
		R_{ij}  &=  \max\left\{ 0, 2 d_{ij}\left[ \overline{(\rho\vel)}_{ij}\cdot\f_{ij}^{\rho\vel,*} -\overline{\rho}_{ij}  f_{ij}^{\rho E,*} - \overline{(\rho E)}_{ij} f_{ij}^{\rho,*}  \right]\right\}\\
		&+ \max \left\{ 0, \frac{| \f_{ij}^{\rho\vel,*} |^2}{2} - f_{ij}^{\rho E,*} f_{ij}^{\rho,*} \right\}.
	\end{split}
\end{equation*}
To ensure continuous dependence
of the limited flux $\alpha_{ij} f_{ij}^*$ on the data and avoid convergence problems at
steady state, we overestimate the bound $R_{ij}$ by \cite{kuzmin2020, kuzmin2023}
\begin{equation*}
	\begin{split}
		R_{ij}^{\max} &= 2 d_{ij}\left[\max \left\{ |\overline{(\rho\vel)}_{ij}|, |\overline{(\rho\vel)}_{ji}|\right\} | \f_{ij}^{\rho\vel,*}|\right.\\
		&\left.+ \max\left\{ \overline{\rho}_{ij}, \overline{\rho}_{ji}\right\}  |f_{ij}^{\rho E,*}| + \max\left\{ \overline{(\rho E)}_{ij}, \overline{(\rho E)}_{ji}\right\} |f_{ij}^{\rho,*}|\right]\\
		&+ \max \left\{ 0, \frac{| \f_{ij}^{\rho\vel,*} |^2}{2} - f_{ij}^{\rho E,*} f_{ij}^{\rho,*} \right\} = R_{ji}^{\max}
	\end{split}
\end{equation*}
but remark that this overestimation may increase the levels of artificial viscosity \cite{rueda2024}.

Taking the symmetry condition $\alpha_{ij}= \alpha_{ji}$ and the constraint
$\alpha_{ji}R_{ji}\le Q_{ji}$ into account, we apply
\begin{equation*}
	\alpha_{ij} = \begin{cases}
		\frac{\min(Q_{ij}, Q_{ji})}{R_{ij}^{\max}}& \text{if }R_{ij}^{\max}>\min(Q_{ij}, Q_{ji}),\\
		1 & \text{otherwise}
	\end{cases}
\end{equation*}
to all components of the prelimited flux $f_{ij}^*$. This final limiting step ensures
that $\overline{u}_{ij}^{*,\alpha}\in\GG$ whenever $\overline{u}_{ij}\in\GG$. The IDP property of the fully
implicit MCL scheme
\begin{equation}\label{eq:implicitEulerMCL}
  u_i^{n+1} = u_i^n + \frac{\Delta t}{m_i}\left[
  \red{\sum_{e\in \mathcal{E}_i}\sum_{S\in\mathcal S_\Gamma(K_e)}
  \sigma_{i,S}\lambda_{i,S}(\bar u_{i,S}^{n+1}-u_i^{n+1})}+
  \sum_{j\in\mathcal{N}_i^*}2d_{ij}(\overline{u}_{ij}^{*,\alpha,n+1}-u_i^{n+1})\right]
\end{equation}
can now be shown in exactly the same manner as for the low-order method
\eqref{eq:implicitEuler} in Section~\ref{sec:IDP}.

\section{Solution of nonlinear systems}
\label{sec:Solver}

In principle, we could solve \eqref{eq:implicitEulerMCL} using the IDP fixed-point
iteration \eqref{eq:fixpIDP} with the mapping
\begin{align}\label{eq:psisMCL}
  \Psi^s_i(u) &= \frac{s u_i^n}{1+ s}+ \frac{1}{1+s}\biggl(u_i+ \frac{s\Delta t}{m_i}
  \Big[\red{\sum_{e\in \mathcal{E}_i}\sum_{S\in\mathcal S_\Gamma(K_e)}
  \sigma_{i,S}\lambda_{i,S}(\bar u_{i,S}-u_i)}
    +\sum_{j\in\mathcal{N}^*_i}2d_{ij}(\overline{u}_{ij}^{*,\alpha}-u_i)\Big]\biggr).
\end{align}
The analysis performed in Section \ref{sec:IDP} guarantees the
IDP property of intermediate solutions under
condition \eqref{eq:cflcond}, which can always be enforced by
choosing a sufficiently small value of the parameter $s>0$.
However, the convergence
behavior of the simple fixed-point
iteration \eqref{eq:psisMCL}
is inferior to that of the quasi-Newton /
deferred correction methods to be discussed in this section.

\magenta{\begin{remark}
    We use a spatially constant time step $\Delta t$ in this work. Further speedups of steady-state
    computations could be achieved using local time stepping. While the use of different time steps
    for different nodes results in a lack of the discrete conservation property for intermediate numerical
    solutions, steady states are independent of the evolution history and, therefore, conservative.
\end{remark}}  

\subsection{Solver for individual time steps}
\label{sec:LoJac}

The flux-corrected nonlinear system \eqref{eq:implicitEulerMCL} can be
written in the following split matrix form:
\begin{equation}\label{eq:matrixform}
M_Lu^{n+1}=M_Lu^n+\Delta t[R_L(u^{n+1})+F^*(u^{n+1})].
\end{equation}
We denote by $M_L=(\delta_{ij}m_iI_m)_{i,j=1}^{N_h}$ the lumped mass matrix.
The components of the low-order steady-state residual $R_L(u)$ and of
the high-order correction term $F^*(u)$ are defined by
$$(R_L(u))_i=\tilde{b}_i(u_i, \hat u)+ \sum_{j\in\mathcal{N}_i^*} [d_{ij}(u_j-u_i) - (\mathbf{f}_j - \mathbf{f}_i )\cdot\mathbf{c}_{ij}],\qquad
F_i^*(u) = \sum_{j\in \mathcal{N}_i^*}\alpha_{ij} f_{ij}^*.$$
A quasi-Newton method for solving \eqref{eq:matrixform} uses an approximate
Jacobian $J(u)$ in updates of the form 
\begin{equation}\label{eq:quasinewtom}
  J(u^{(k)})(u^{(k+1)}-u^{(k)})=M_L(u^n-u^{(k)}) +\Delta t[R_L(u^{(k)})+F^*(u^{(k)})],\quad
      k\in\N_0,
\end{equation}
where $u^{(k)}$ is a given approximation to $u^{n+1}$. By default, we set
$u^{(0)}=u^n$ and use a low-order Jacobian approximation $J(u)=\J(u)$, the
  derivation of which is explained below.

The homogeneity property $\f(u) = \A(u)u$ of the flux function of the
Euler equations implies that
\begin{equation}\label{eq:homoegenity}
\mathbf{f}_j - \mathbf{f}_i=\A(u_j)u_j-\A(u_i)u_i.
\end{equation}
Dolej{\v{s}}{\i} and Feistauer \cite{dolejvsi2004} used this property to
design a conservative linearized implicit scheme. Further representatives
of such schemes were derived and analyzed by Ku{\v{c}}era et al. \cite{kucera2022}.
The algorithm employed by Gurris et al. \cite{gurris2012} is similar but provides
the option of performing multiple iterations.

In view of \eqref{eq:homoegenity}, the Jacobian associated with the low-order component
of \eqref{eq:matrixform} is given by \cite{gurris2009}
\begin{equation*}
	\J = M_L + \Delta t[A - B-D],
\end{equation*}
where $A=(A_{ij})_{i,j=1}^{N_h}$ and $D=(D_{ij})_{i,j=1}^{N_h}$
are sparse block matrices composed from
\begin{equation*}
	A_{ij} = \con_{ij}\cdot \A(u_j),\quad D_{ij} = d_{ij} I_m.
\end{equation*}
The block-diagonal matrix $B=(\delta_{ij}B_{ii})_{i,j=1}^{N_h}$ is \red{an approximate Jacobian
associated with the implicit part of the boundary terms defined by \eqref{eq:bdrtermLO}. It is} constructed
using an algebraic splitting
\begin{equation}\label{eq:btermsplit}
  \tilde b(u, \hat u)= B(u)u + b(\hat u)
\end{equation}
\red{of $\tilde b(u, \hat u)=(\tilde b_i(u_i, \hat u))_{i=1}^{N_h}$ into a
matrix-vector product $B(u)u$ that depends only on $u$ and a remainder
$b(\hat u)=(b_i(\hat u))_{i=1}^{N_h}$ that depends on the
boundary data $\hat u$. This splitting is generally nonunique.}
\medskip

\red{
Using the representation \eqref{eq:bdrtermLO} of the boundary term $\tilde b(u, \hat u)$, we take
advantage of the fact that
$$\f(u_i)-\f(\hat u_{i,S})=\A(u_i)u_i - \A(\hat u_{i,S})\hat u_{i,S}
$$
by the homogeneity of $\f(u)$. 
The vector $\tilde b(u, \hat u)$ of boundary integrals can
be written in the form  \eqref{eq:btermsplit} if the diagonal blocks of
 $B(u)$ and the components of $\hat b(u)$ are defined as follows:
\begin{align*}
  B_{ii} &= \sum_{e\in \mathcal{E}_i}\sum_{S\in\mathcal S_\Gamma(K_e)}B_{ii,S},\qquad  B_{ii,S}=\frac{\sigma_{i,S}}2
 [\A(u_i)\cdot\mathbf n_S -\lambda_{i,S}I_m],\\
 b_i &= \sum_{e\in \mathcal{E}_i}\sum_{S\in\mathcal S_\Gamma(K_e)}b_{i,S},\qquad
b_{i,S}=\frac{\sigma_{i,S}}2
  [\lambda_{i,S}\hat u_{i,S}-\f(\hat u_{i,S})\cdot \n_S].
\end{align*}
\blue{We use the above definitions of the face contributions $B_{ii,S}$ and $b_{i,S}$ by default.}
If there exists a transformation matrix $\hat{B}_{i,S}\in\R^{m\times m}$ such that the external state
$\hat u_{i,S} = \hat{B}_{i,S}u_i$ depends only on the internal state $u_i$, we 
treat the contribution of $S\in\mathcal S_\Gamma(K_e)$ implicitly by setting 
\begin{equation*}
  B_{ii,S} = \frac{\sigma_{i,S}}2
  [\A(u_i)\cdot\mathbf n_S-\A(\hat u_{i,S})\cdot\mathbf n_S\hat B_{i,S} 
    +\lambda_{i,S}(\hat B_{i,S}-I_m)],\qquad b_{i,S}(\hat u)=0.
\end{equation*}
For example,  $\hat B_{i,S} = I_m$} in the case of
a supersonic outlet boundary condition. A two-dimensional reflecting
wall boundary condition can be implemented using the transformation matrix
\begin{equation*}
	\hat B = \left(\begin{matrix}
		1 & 0 & 0 & 0\\
		0 & 1 - 2n_1n_1 & - 2n_1 n_2 & 0\\
		0 & -2 n_2 n_1 & 1- 2n_2n_2 & 0\\
		0 & 0 & 0& 1
	\end{matrix}\right),
\end{equation*}
where $n_1$ and $n_2$ are the two components of the unit outward normal $\n = (n_1, n_2)$. This implicit treatment of the wall boundary condition preserves the conservation property in each iteration of the quasi-Newton method \eqref{eq:quasinewtom} with $J(u)=J_L(u)$. An explicit treatment of $b(\hat u)$ may give rise to conservation errors in intermediate solutions. Hence, conserved quantities may enter or exit the domain through a solid wall if the iterative process is terminated before full convergence is achieved.
\medskip


Substituting the low-order Jacobian $J(u)=J_L(u)$ into \eqref{eq:quasinewtom}, we find that the resulting linear system for $u^{(k+1)}$ has the structure of the \emph{deferred correction} method
\begin{equation}\label{eq:JacBE}
	\J(u^{(k)}) u^{(k+1)} = M_L u^n +\Delta t( F^*(u^{(k)}) + b(\hat u)).
\end{equation}
Criteria for convergence of such fixed-point iterations can be found, e.g., in \magenta{\cite[Prop. 4.3]{abgrall2017b}} and \cite[Prop. 4.66]{lohmann2019}.
If we always stop after the first iteration, then
 $u^{n+1}=u^{(1)}$ is the solution of ~\cite{gurris2009, kuzmin2012b}
\begin{equation}\label{eq:JacBE1}
	\J(u^n) u^{n+1} = M_L u^n +\Delta t( F^*(u^n) + b(\hat u)).
\end{equation}
This approximation to the nonlinear problem \eqref{eq:matrixform}
falls into the category of linearly implicit schemes analyzed
in \cite{kucera2022}. It may perform very well as long as
the time steps
are sufficiently small. However, the linearized form \eqref{eq:JacBE1} of
 \eqref{eq:matrixform}
is generally not IDP. Moreover, it is difficult (if not impossible)
to derive a CFL-like condition under which \eqref{eq:JacBE1}
would \emph{a priori} guarantee the IDP property
of $u^{n+1}$.

\blue{We use the iterative version \eqref{eq:JacBE} of \eqref{eq:JacBE1}
 in this paper because if the fixed-point iteration \eqref{eq:fixpIDP} 
 converges to a unique IDP solution of the nonlinear system \eqref{eq:implicitEulerMCL},
 then this solution can usually be calculated much faster using
 \eqref{eq:JacBE}. Even for time steps corresponding to very large maximum
 CFL numbers,
a single iteration is often sufficient to obtain an IDP result.}
In general, we \magenta{set $u^{(0)}=u^n$ and
update the successive approximations  $u^{(k)},\ k=0,1,2,\ldots$} as follows:
\begin{equation}\label{eq:fpiter}
	\begin{split}
        \Delta u^{(k)} =& [\J(u^{(k)})]^{-1} R_{\Delta t}^*(u^{(k)}),\\
		u^{(k+1)} =& u^{(k)}+ \Delta u^{(k)},
	\end{split}
\end{equation}
where
\begin{align*}
  \blue{R_{\Delta t}^*(u) =M_L(u^n-u)+ \Delta t [R_L(u) +F^*(u)]}
\end{align*}
is the residual of \eqref{eq:matrixform}. \blue{Note that the unsplit form
of the boundary term $\tilde b(u,\hat u)$ is built into $R_L(u)$. A~particular
choice of the splitting \eqref{eq:btermsplit} that defines the component
$B(u)$ of the approximate Jacobian
$J_L(u)$ may affect the convergence behavior of the solver but not a converged
result\footnote{\blue{Convergence of the deferred correction method \eqref{eq:JacBE}  to a fixed point of
   the IDP mapping $\Psi^s(u)$
  can always be achieved using adaptive (pseudo-)time stepping,
  as in the challenging Mach 20 bow shock test of Section \ref{sec:Mach20}.}}.}  
  
Instead of iterating until a standard
stopping criterion, such as
$$
\frac{\|R_{\Delta t}^*(u^{(k)})\|}{ \|R_{\Delta t}^* (u^{(0)})\|} < \varepsilon_1,
\qquad \frac{\|\Delta u^{(k)}\|}{\|u^{(k)}\|} < \varepsilon_2\qquad \mbox{for}
\quad 0< \varepsilon_1,\varepsilon_2\ll 1,
$$
is met, we use the IDP property as the stopping criterion. That is, we exit and set $u^{n+1}=u^{(k+1)}$ if $u^{(k+1)}\in G$. In practice, a single iteration is usually enough regardless of the time step $\Delta t$.

\begin{remark}
  If multiple iterations are necessary to satisfy the IDP criterion, the approximate Jacobian  $\J(u^{n})$ may be used instead of $\J(u^{(k)}),\ k\ge 1$ to reduce the cost of matrix assembly.
\end{remark}

\subsection{Steady-state solver}
\label{sec:steadystate}

Let us now turn to the computation of steady states $u$, for which our system
\eqref{eq:matrixform} reduces to
\begin{equation}\label{MCLsteady}
R_L(u) +F^*(u)=0.
\end{equation}
Iterations of a quasi-Newton solver for this nonlinear system
can be interpreted as pseudo-time steps. To achieve optimal convergence,
the approximate Jacobian must be close enough to the exact one. Moreover, the
initial guess $u^0$ must be close enough to a steady-state solution. If
the latter requirement is not met, a Newton-like method may fail to
converge completely. 
Globalization and robustness enhancing strategies are commonly employed if
the nonlinearity is very strong. 
A popular approach is the use of \emph{underrelaxation} (also known as \emph{backtracking} or \emph{damping}) techniques \cite{patankar1980, ferziger2002, gurris2009, lohmann2021, knoll2004}.

Underrelaxation can  be performed in an implicit way by adding positive numbers
to the diagonal elements of the approximate Jacobian \cite{patankar1980, ferziger2002, gurris2009}. The use of smaller pseudo-time steps in an implicit time marching procedure
has the same effect of enhancing the diagonal dominance.
An explicit underrelaxation procedure controls the step size in the direction $\Delta u^n$ and produces the update 
\begin{equation}\label{eq:ur}
	u^{n+1} = u^n+\omega\Delta u^n,\quad \omega\in(0,1].
\end{equation}

In this work, we solve \eqref{MCLsteady} using \eqref{eq:JacBE} as a pseudo-time stepping
method. The amount of implicit underrelaxation depends on the choice of the parameter
$\Delta t$, which we fit to a fixed upper bound for CFL numbers. Additionally, we perform
explicit underrelaxation in the following way:
\begin{enumerate}
	\item Given $u^n$ compute an IDP approximation $\tilde u^{n+1}$ using \eqref{eq:fpiter} with IDP stopping criterion.
	\item Set $\Delta u^n = \tilde u^{n+1} -u^n$, choose $\omega\in(0,1]$, and calculate the update~\eqref{eq:ur}. 
\end{enumerate}
Since $\tilde u^{n+1}$ and $u^n$ are in the admissible set $G$, so is $u^{n+1}$ for all $\omega\in (0,1]$. The underrelaxation factor $\omega$ can be chosen adaptively or assigned a constant value that is the same for all pseudo-time steps. In our experience, $\omega = 0.5$ is typically sufficient to ensure that the steady-state residual
  $$
R_\infty^*(u)=R_L(u) +F^*(u)
  $$
  becomes small enough.
  However, the rates of convergence may be unsatisfactory for fixed $\omega$.

The explicit underrelaxation technique proposed by Ranocha et al.~\cite{ranocha2020}
is designed to enforce fully discrete entropy stability in the final stage of a
Runge--Kutta method for an evolutionary problem. Badia et al.~\cite{badia2017a} select underrelaxation
factors that minimize the residuals. Combining these ideas,
we propose an adaptive relaxation strategy based on the values of the entropy
residuals
\begin{equation*}
	(R_\eta^*(u))_i = v(u_i)^\top\left(M_L^{-1}R_\infty(u)\right)_i,
\end{equation*}
where  $v(u_i)=\eta'(u_i)$ is the vector of entropy variables corresponding
to a convex entropy $\eta(u)$. The physical entropy $s = \log(p\rho^{-\gamma})$
 is concave. We define $R_\eta^*(u)$ using the mathematical
entropy 
$$
\eta(u)=-\frac{\rho s}{\gamma-1}.
$$

\begin{remark}
We monitor the entropy residual $R_\eta^*(u)$
rather than the residual
$R_\infty^*(u)$ of the discretized Euler system
 because $(R_\eta^*(u))_i$ is a scalar
 quantity, whereas $(R_\infty^*(u))_i$ is a vector whose  components are
 `apples and oranges' with different
physical dimensions / orders of magnitude.
\end{remark}

Let the discrete $L^2$ norm $\|\cdot\|_{2,h}: \R^{N_h}\rightarrow [0,\infty)$
  be defined by
\begin{equation*}
  \|w\|_{2,h} =\sqrt{w^\top M_Cw}\qquad \forall w \in\R^{N_h},
\end{equation*}
where $M_C=(m_{ij})_{i,j=1}^{N_h}$ is the consistent mass matrix. Note that
$\|w\|_{2,h}=\|w_h\|_{L^2(\Omega)}$ for a finite element function
$w_h\in W_h$ with degrees of freedom corresponding to the components of
$w \in\R^{N_h}$.

The 
underrelaxation factor that minimizes the objective function $\|R_{\eta}^*(u^{n+1)}\|_{2,h}$
  is given by
\begin{equation}\label{eq:omegaopt}
	\omega^n = \argmin_{\tilde\omega\in(0,1]}\|R_{\eta}^*(u^n + \tilde\omega\Delta u^n)\|_{2,h}.
\end{equation}
This minimization problem has no closed-form solution, while numerical computation of $\omega^n$ is costly due to the strong nonlinearity of the problem at hand. It is not worthwhile to invest inordinate effort in solving \eqref{eq:omegaopt}. A rough approximation is sufficient to find a usable scaling factor for adaptive relaxation purposes (see, e.g.,~\cite[Sec. 3.3.1]{turek1999}). Following Lohmann~\cite{lohmann-preprint}, we approximate $\omega^n$ by 
\begin{equation}\label{eq:aur}
	\omega^n_K = \argmin_{\tilde\omega\in\overline{\omega}_{K}}\|R_{\eta}^*(u^n+\tilde\omega\Delta u^n)\|_{2,h}, \quad K\in\N,
\end{equation}  
where
\begin{equation*}
  \overline\omega_K \coloneqq \left\{\omega_0 + \frac{k-1}{K -1}(1-\omega_0),\quad
  k=1,\ldots, K\right\}
\end{equation*}
is a discrete set of candidate parameter values.
The approximation \eqref{eq:aur}
to~\eqref{eq:omegaopt} becomes more accurate as $K$ is increased and $\omega_0$ is decreased.
In our numerical experiments, we use $K = 3$ and $\omega_0 = 0.5$. That is, we adaptively
select relaxation parameters belonging to the set $\overline\omega_3 = \{0.5, 0.75, 1\}$.

\section{Numerical examples}
\label{sec:examples}

To study the steady-state convergence behavior of the presented numerical algorithms, we apply them
to standard two-dimensional test problems in this section. In our
numerical studies, the value of the pseudo-time step $\Delta t$ is
determined using the formula~\cite{guermond2018,kuzmin2020, kuzmin2023}
\begin{equation}\label{eq:CFL}
	\max_{i\in\{1,\ldots,N_h\}}\frac{2\Delta t}{m_i}\sum_{j\in\mathcal{N}^*_i}d_{ij} = \mathrm{CFL}
\end{equation}
with a user-defined threshold $\mathrm{CFL}$. Explicit MCL schemes are
IDP for $\mathrm{CFL}=1$. In the implicit case, it is possible to
achieve the IDP property with $\mathrm{CFL}\gg 1$, as we demonstrate below.

Unless stated otherwise, we prescribe the free stream/inflow boundary values as follows~\cite{gurris2009}:
\begin{equation}\label{eq:freestream}
		\rho_\infty = 1,\qquad
		p_\infty = \frac{1}{\gamma},\qquad
		\vel_\infty = (M_\infty, 0)^\top,
\end{equation}
where $M_\infty$ is the free stream Mach number corresponding to the
free stream speed of sound $c_\infty=1$.

The  methods under investigation were implemented using the open-source \texttt{C++} finite element library MFEM~\cite{anderson2021,andrej2024,mfem}. 
To solve the linear system in each step of the fixed-point iteration~\eqref{eq:fpiter}, we use the BiCGstab solver implemented in MFEM and precondition it with Hypre's ILU~\cite{hypre}.
All computations are performed on unstructured triangular meshes generated using Gmsh~\cite{geuzaine2009}. The numerical results are visualized using
Paraview~\cite{ayachit2015}.

We use the low-order solution as initial condition for MCL
and initialize it by extending the free stream values \eqref{eq:freestream}
into $\Omega$.
A numerical solution is considered to be stationary if
\begin{equation*}
	r(u) = \|M_L^{-1}R_\infty^*(u)\|_{2,h} < 10^{-8}.
\end{equation*}

To avoid repetition, let us preview our findings regarding the typical convergence behavior depending on the parameter $\omega$ of the explicit underrelaxation strategy~\eqref{eq:ur} and on the CFL number~\eqref{eq:CFL}. First and foremost, we found that some kind of explicit underrelaxation is an essential prerequisite for steady-state convergence. Without relaxation, i.e., using $\omega \equiv 1$, we were able to achieve convergence for just two simple scenarios. In our experiments with constant underrelaxation factors, we used  $\omega\in\{0.5, 0.6, 0.7, 0.8, 0.9\}$. Steady-state solutions to all test problems could be  obtained with $\omega = 0.5$. The impact of a manually chosen underrelaxation factor on the convergence rate of the pseudo-time stepping method was different in each test. In most cases, the best performance was achieved with the largest value of $\omega$ for which the method does converge. We also report the convergence history of our adaptive underrelaxation strategy~\eqref{eq:aur} with $\omega = \omega_{3}^n$. No stagnation or divergence was observed for this approach. The convergence rates achieved with the proposed adaptive strategy are at least as high as the best rate obtained with a constant damping parameter. Furthermore, increasing the number of candidate underrelaxation factors $\omega\in\overline{\omega}_K$ did not result in significantly faster convergence. Using $\omega_K^n$ defined by~\eqref{eq:aur} with $K=6$, i.e. $\overline{\omega}_6 = \{0.5, 0.6, 0.7, 0.8, 0.9, 1\}$, we observed essentially the same convergence behavior as for $\omega_3^n$, but the computational cost was drastically increased.

In a second set of experiments, we varied the pseudo-time step in our algorithm that uses the adaptive damping parameter~\eqref{eq:aur} to approximately minimize the value of the entropy residual. Specifically, we ran simulations with $\mathrm{CFL}\in\{10^k, k=0,\ldots,5 \}$. In most cases, the proposed scheme converges faster as the CFL number is increased. In the parameter range $\mathrm{CFL}\leq 10^1$, steady-state residuals decrease slowly but almost monotonically. Slow and  oscillatory convergence behavior is observed for $\mathrm{CFL} = 10^2$. The convergence rates become less sensitive to $\Delta t$ for $\mathrm{CFL} \ge 10^3$. For all test cases under consideration, the residuals evolve in the same manner for $\mathrm{CFL} = 10^4$ and $\mathrm{CFL} = 10^5$. The above response to changes in the CFL number is consistent with the behavior of other implicit schemes that use the framework of algebraic flux correction for continuous finite elements~\cite{gurris2009, kuzmin2012b}.

\subsection{GAMM channel}
A well-known two-dimensional benchmark problem for the stationary  Euler equations is the so-called GAMM channel~\cite{feistauer2003}. In this popular test for numerical schemes, the gas enters a channel at free stream Mach number $M_\infty = 0.67$ through the subsonic inlet
\begin{equation*}
	\Gamma_{\mathrm{in}} = \{(x,y)\in\R^2: x = -1, y\in(0,1)\}.
\end{equation*}
The flow accelerates to supersonic speeds over the bump on the lower wall
\begin{equation*}
	\begin{split}
		\Gamma_{\mathrm{lw}} =& \{(x,y)\in\R^2: x\in [-1, -0.5], y = 0\}\\
					&\cup  \{(x,y)\in\R^2: x\in [-0.5, 0.5], y = \sqrt{1.69 - \magenta{x^2}} -1.2\}\\
					&\cup \{(x,y)\in\R^2: x\in [0.5, 1], y = 0\},
	\end{split}
\end{equation*}
 which results in a shock wave, and exits the domain at subsonic speeds at the outlet
 \begin{equation*}
 	\Gamma_{\mathrm{o}} = \{(x,y)\in\R^2: x = 1, y\in(0,1)\}.
 \end{equation*}
 We prescribe reflecting wall boundary conditions on $\Gamma_{\mathrm{lw}}$
 and on the upper wall 
\begin{equation*}
	\Gamma_{\mathrm{uw}} = \{(x,y)\in\R^2: x\in(-1,1), y =1\}.
\end{equation*}

\begin{figure}[h!]
	\centering
	\begin{subfigure}[c]{0.45\textwidth}
		\includegraphics[trim={1cm 2cm 2cm 2cm},clip, width = 0.99\textwidth]{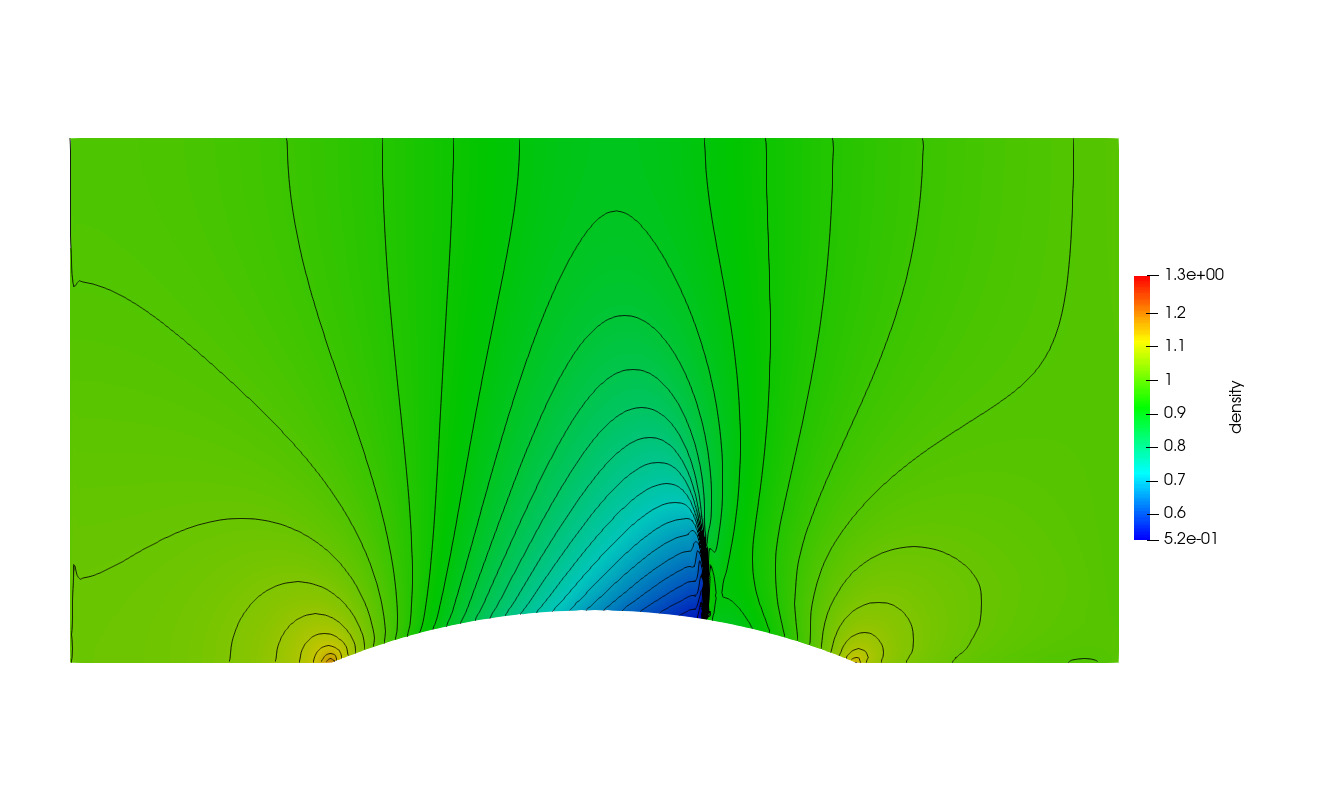}
		\subcaption{Density}
		\label{fig:GAMMden}
	\end{subfigure}
	\begin{subfigure}[c]{0.45\textwidth}
		\includegraphics[trim={1cm 2cm 2cm 2cm},clip,width = 0.99\textwidth]{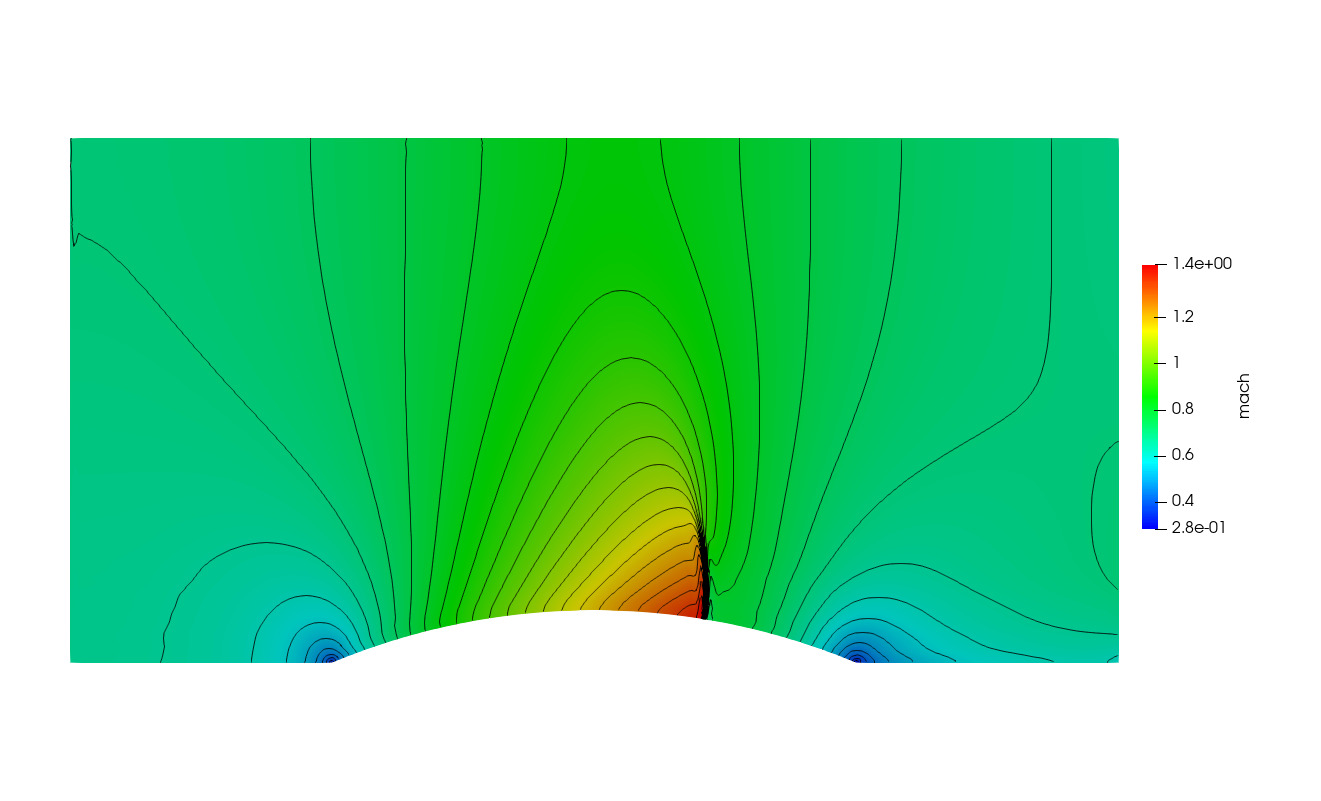}
		\subcaption{Mach number}
		\label{fig:GAMMmach}
	\end{subfigure}
	\begin{subfigure}[c]{0.45\textwidth}
		\includegraphics[trim={1cm 2cm 2cm 2cm},clip, width = 0.99\textwidth]{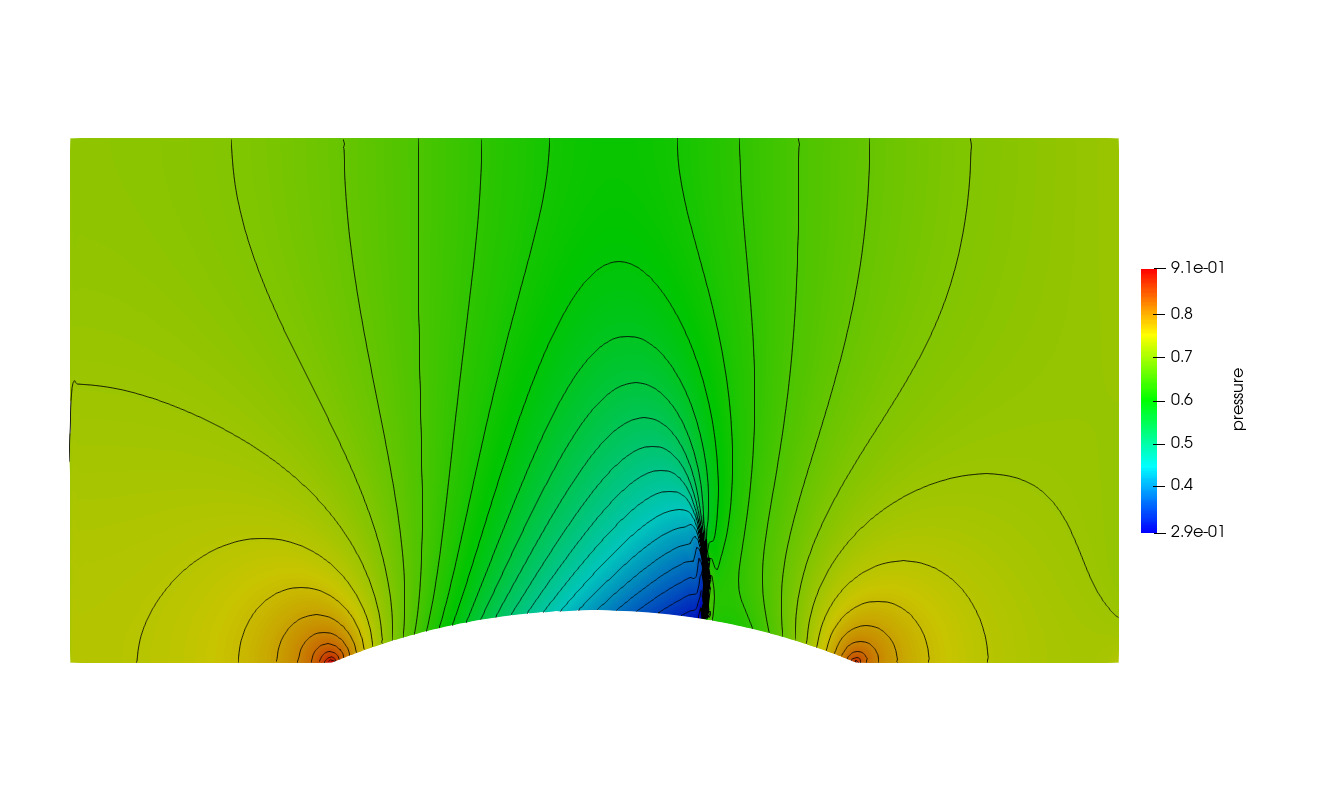}
		\subcaption{Pressure}
		\label{fig:GAMMpres}
	\end{subfigure}
	\caption{GAMM channel ($M_\infty= 0.67$), steady-state MCL results obtained with $N_h = 112,963$ unknowns per component on a mesh consisting of $E_h = 224,576$ triangles.}
	\label{fig:GAMM}
\end{figure}

\begin{figure}[h!]
	\centering
	\begin{subfigure}[c]{0.40\textwidth}
		\includegraphics[trim={0.9cm 0.3cm 2.15cm 1cmm},clip, width = 0.95\textwidth]{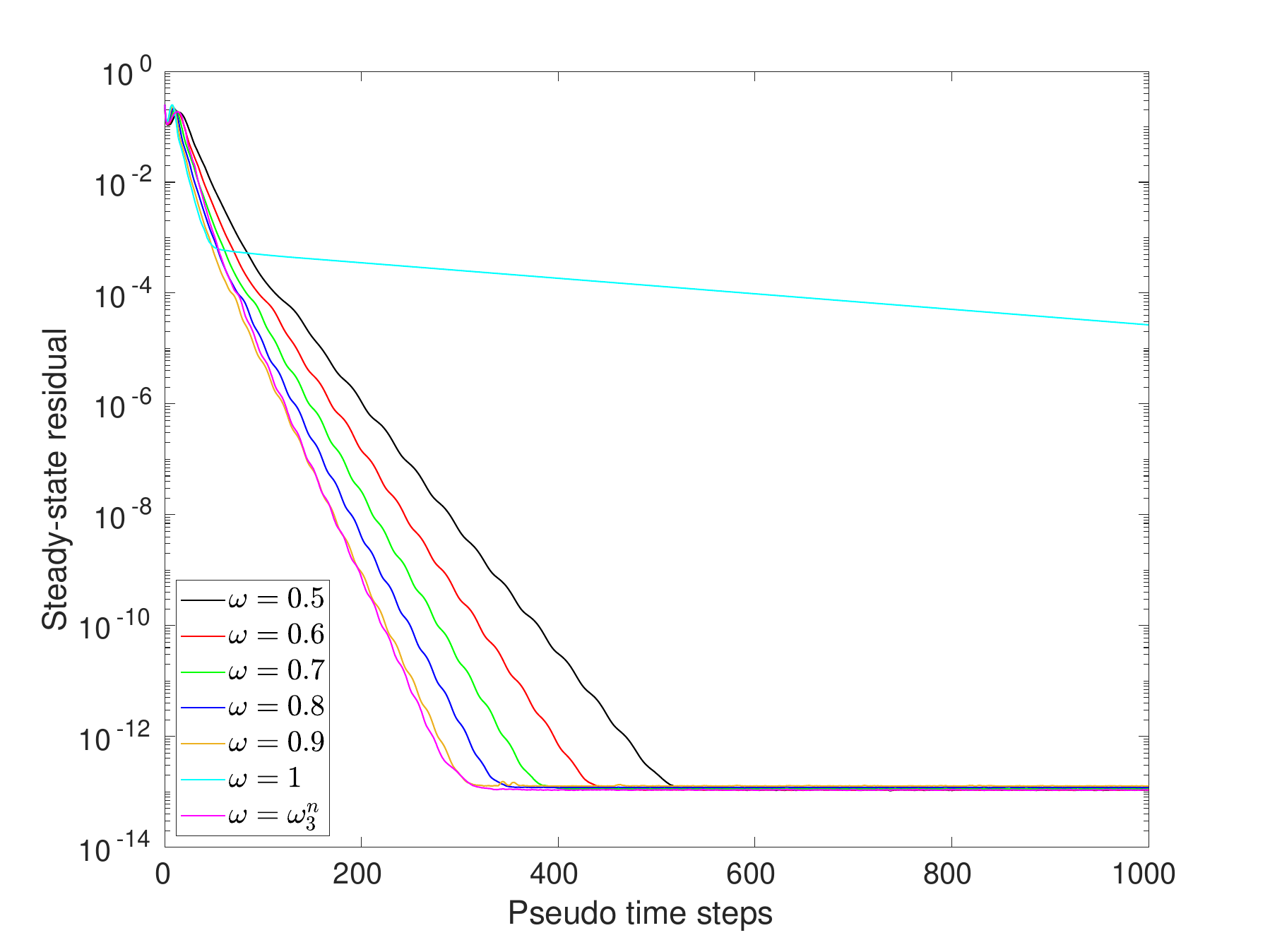}
		\subcaption{Different underrelaxation factors, $\mathrm{CFL} = 10^4$.}
		\label{fig:GAMM_ur}
	\end{subfigure}
	\begin{subfigure}[c]{0.40\textwidth}
		\includegraphics[trim={0.9cm 0.3cm 2.15cm 1cmm},clip, width = 0.95\textwidth]{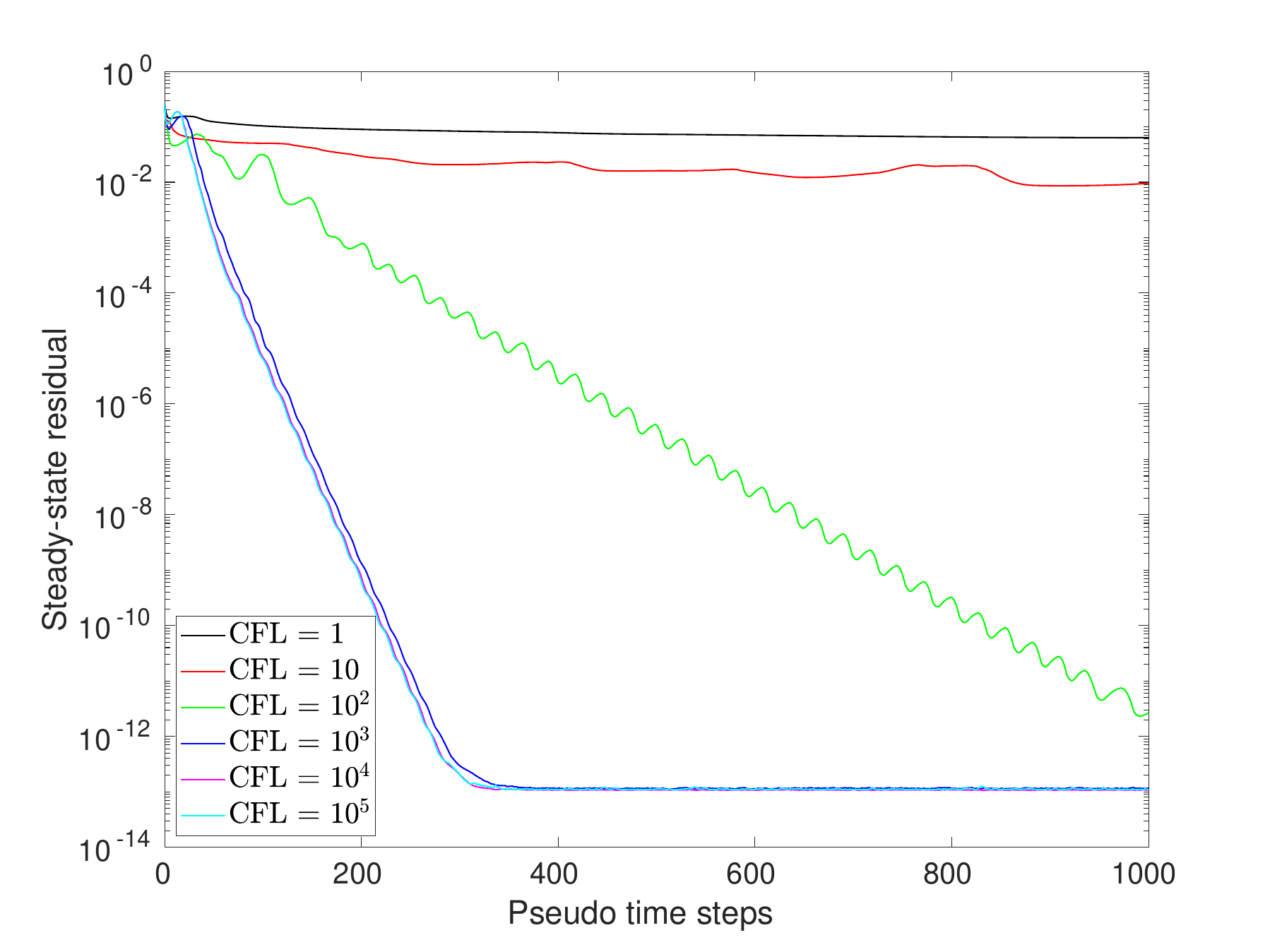}
		\subcaption{Different CFL numbers,
                  $\omega=\omega_3^n$.}
		\label{fig:GAMM_CFL}
	\end{subfigure}
	\caption{GAMM channel ($M_\infty = 0.67$), steady-state convergence history for a long-time MCL simulation on an unstructured triangular mesh with $N_h = 24,037$ nodes and $E_h = 47,256$ cells.}
	\label{fig:GAMM_SSR}
\end{figure}

The numerical results presented in Fig.~\ref{fig:GAMM} are MCL
approximations to the
density, Mach
number, and pressure at steady state. In Fig.~\ref{fig:GAMM_ur},
we plot the evolution history of steady-state residuals for $\mathrm{CFL}=10^4$
and various choices of the underrelaxation factor $\omega$. While convergence is extremely slow for $\omega=1.0$, all kinds of explicit underrelaxation reduce the total number of pseudo-time steps significantly. The best convergence rate is achieved with the adaptive underrelaxation strategy~\eqref{eq:aur}, which meets the stopping criterion $r(u) \le 10^{-13}$ after approximately 300 pseudo-time steps.

Figure~\ref{fig:GAMM_CFL} shows the effect of changing the CFL number on the performance of the fixed-point iteration equipped with the adaptive underrelaxation strategy. We refer to the beginning of Section~\ref{sec:examples} for a general discussion of the way in which $\omega$ and CFL influence the
convergence behavior.

\subsection{Converging-diverging nozzle}

Next, we simulate steady gas flows in a two-dimensional converging-diverging nozzle \cite{gurris2009,hartmann2002}. To evaluate the performance of our method under subsonic and transonic flow conditions, we consider three setups corresponding to different free stream Mach numbers. In all cases, we prescribe a subsonic inflow condition on the left boundary $\Gamma_{\mathrm{in}} =\{ (x,y) \in\R^2: x=-2,\, -1\leq y\leq 1\}$.
The upper and lower reflecting walls of the nozzle are defined by $\Gamma^{\pm}_{\mathrm{w}} = \{(x,g^{\pm}(x))\in\R^2: -2\leq x\leq 8\},$ where \cite{hartmann2002}
\begin{equation*}
	g^{\pm}(x) =\begin{cases}
			\pm 1 &\text{if }-2\leq x\leq 0,\\
			\pm \frac{\cos(\frac{\pi x}{2}) + 3}{4} &\text{if } 0 < x \leq 4,\\
			\pm 1 &\text{if } 4< x\leq 8.
	\end{cases}
\end{equation*}
On the right boundary $\Gamma_{\mathrm{in}} =\{ (x,y) \in\R^2: x=8,\, -1\leq y\leq 1\}$
of the domain, we prescribe a subsonic or supersonic outlet condition depending on the test case.

In the first setup \cite{gurris2009}, we define the free stream values on
$\Gamma_{\mathrm{in}}$ using
~\eqref{eq:freestream} with $M_\infty = 0.2$ and impose
a subsonic outlet condition.
The flow stays subsonic and accelerates to $M\approx 0.48$ at the throat before decelerating in the diverging part of the nozzle.
The characteristic stiffness associated with low Mach numbers makes this problem very challenging despite the fact that the solution is smooth.

\begin{figure}[t!]
	\centering
	\begin{subfigure}[c]{0.7\textwidth}
		\includegraphics[trim={1cm 9cm 0cm 8.5cm},clip, width = \textwidth]{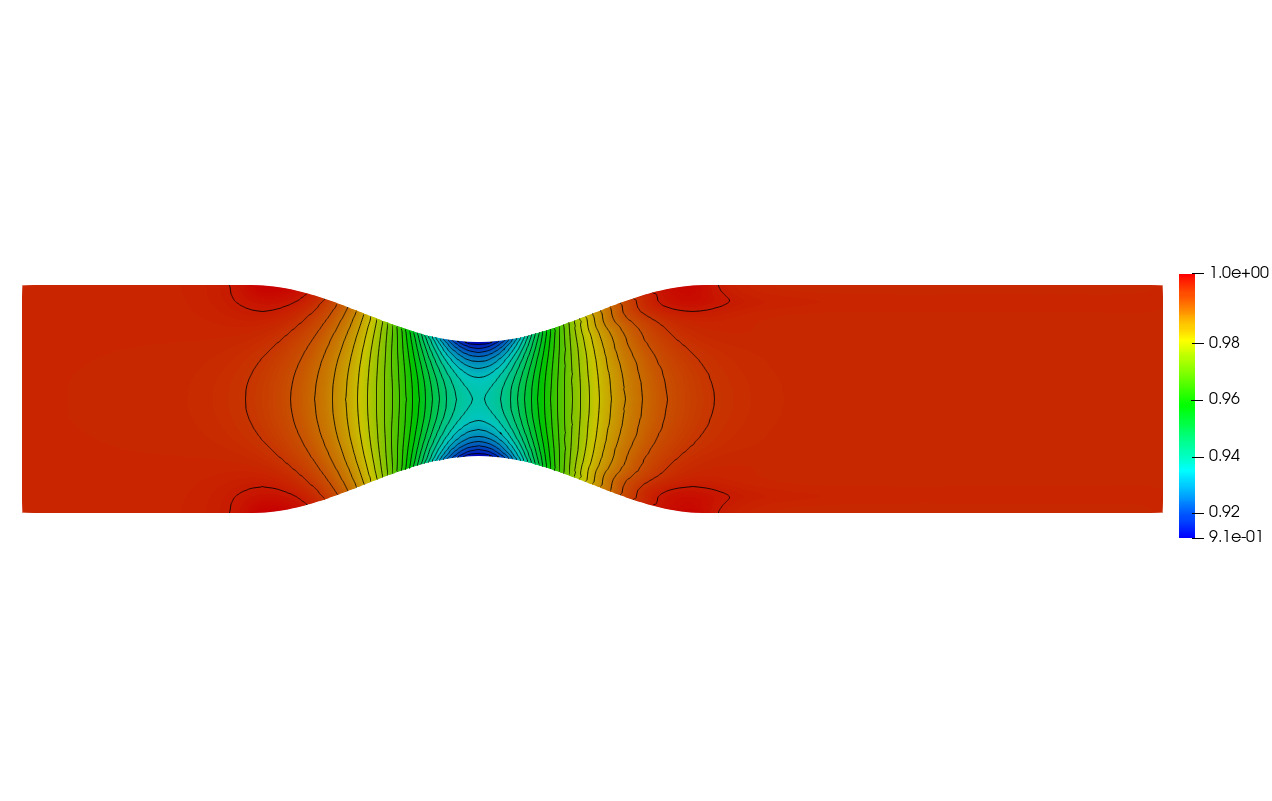}
		\subcaption{Density}
		\label{fig:02Nozzleden}
	\end{subfigure}
	\begin{subfigure}[c]{0.7\textwidth}
		\includegraphics[trim={1cm 9cm 0cm 8.5cm},clip, width = \textwidth]{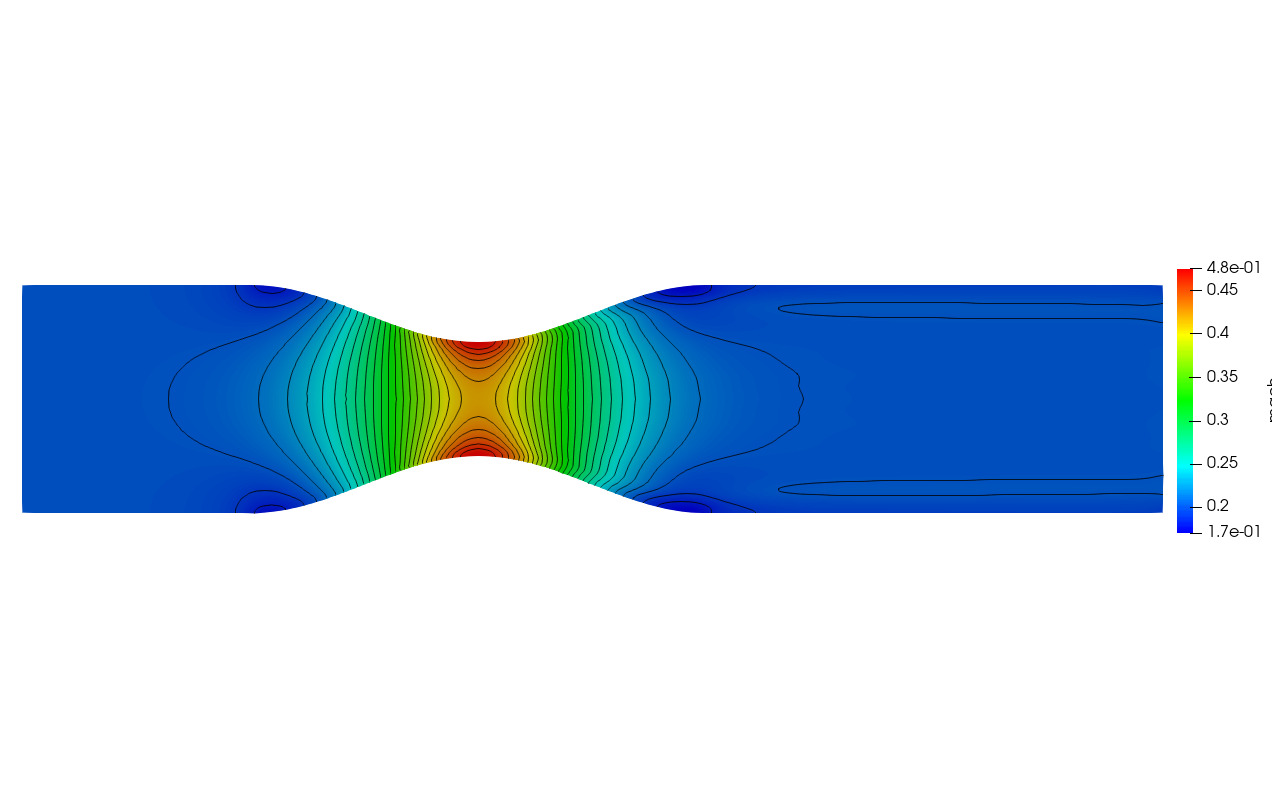}
		\subcaption{Mach number}
		\label{fig:02Nozzlemach}
	\end{subfigure}
	\begin{subfigure}[c]{0.7\textwidth}
		\includegraphics[trim={1cm 9cm 0cm 8.5cm},clip, width = \textwidth]{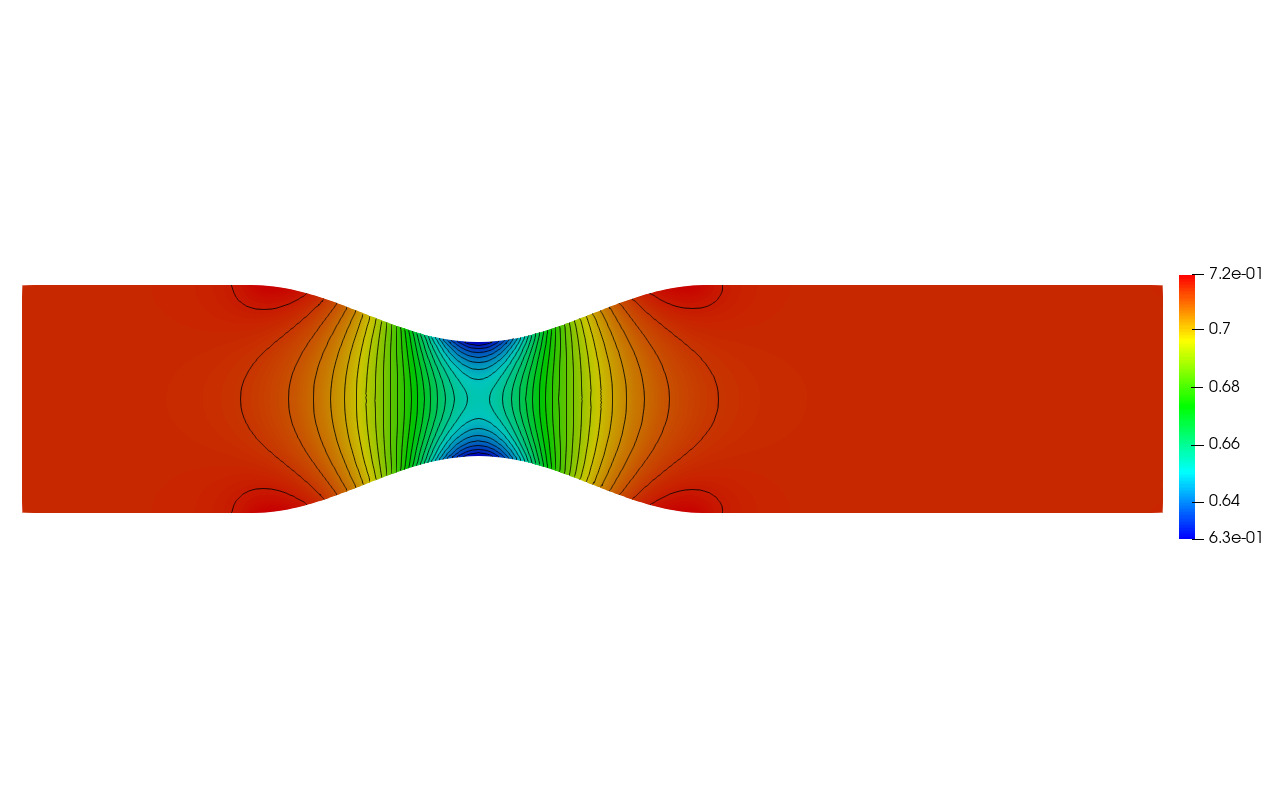}
		\subcaption{Pressure}
		\label{fig:02Nozzlepres}
	\end{subfigure}
	\caption{Subsonic nozzle ($M_\infty= 0.2$), steady-state MCL results obtained with $N_h = 95,329$ unknowns per component on a mesh consisting of $E_h = 189,024$ triangles.}
	\label{fig:02Nozzle}
\medskip   

	\centering
	\begin{subfigure}[c]{0.40\textwidth}
		\includegraphics[trim={0.9cm 0.3cm 2.15cm 1cmm},clip, width = 0.95\textwidth]{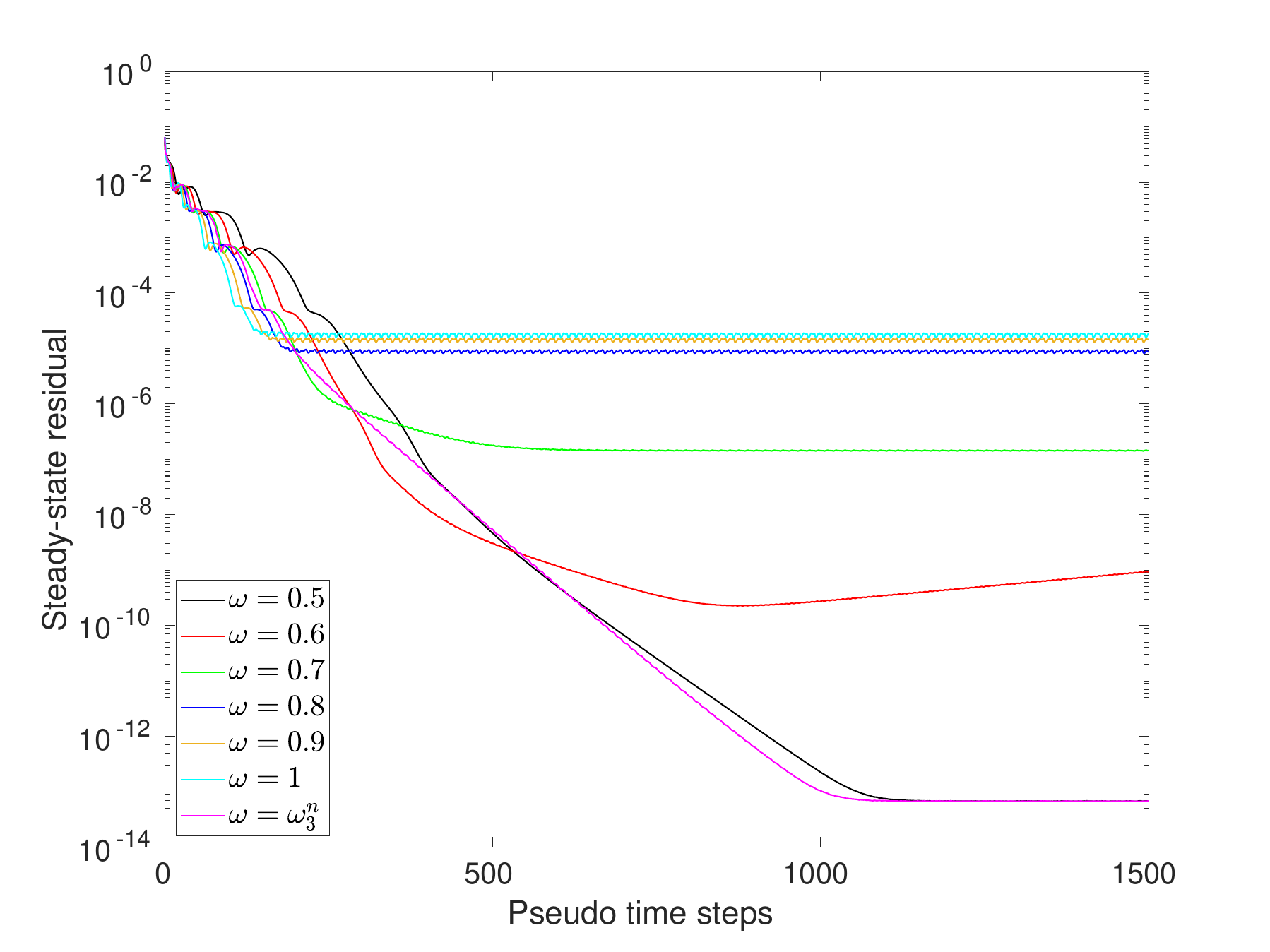}
		\subcaption{Different underrelaxation factors, $\mathrm{CFL} = 10^4$.}
		\label{fig:02Nozzle_ur}
	\end{subfigure}
	\begin{subfigure}[c]{0.40\textwidth}
		\includegraphics[trim={0.9cm 0.3cm 2.15cm 1cmm},clip, width = 0.95\textwidth]{SSR_SubSonicNozzle_CFL.pdf}
		\subcaption{Different CFL numbers, $\omega = \omega_3^n$.}
		\label{fig:02Nozzle_CFL}
	\end{subfigure}
	\caption{Subsonic nozzle ($M_\infty= 0.2$), steady-state convergence history for a long-time MCL simulation on an unstructured triangular mesh with $N_h = 24,037$ nodes and $E_h = 47,256$ cells.}
	\label{fig:02Nozzle_SSR}
\end{figure}

\begin{figure}[h!]
	\centering
	\begin{subfigure}[c]{0.7\textwidth}
		\includegraphics[trim={1cm 9cm 0cm 8.5cm},clip, width = \textwidth]{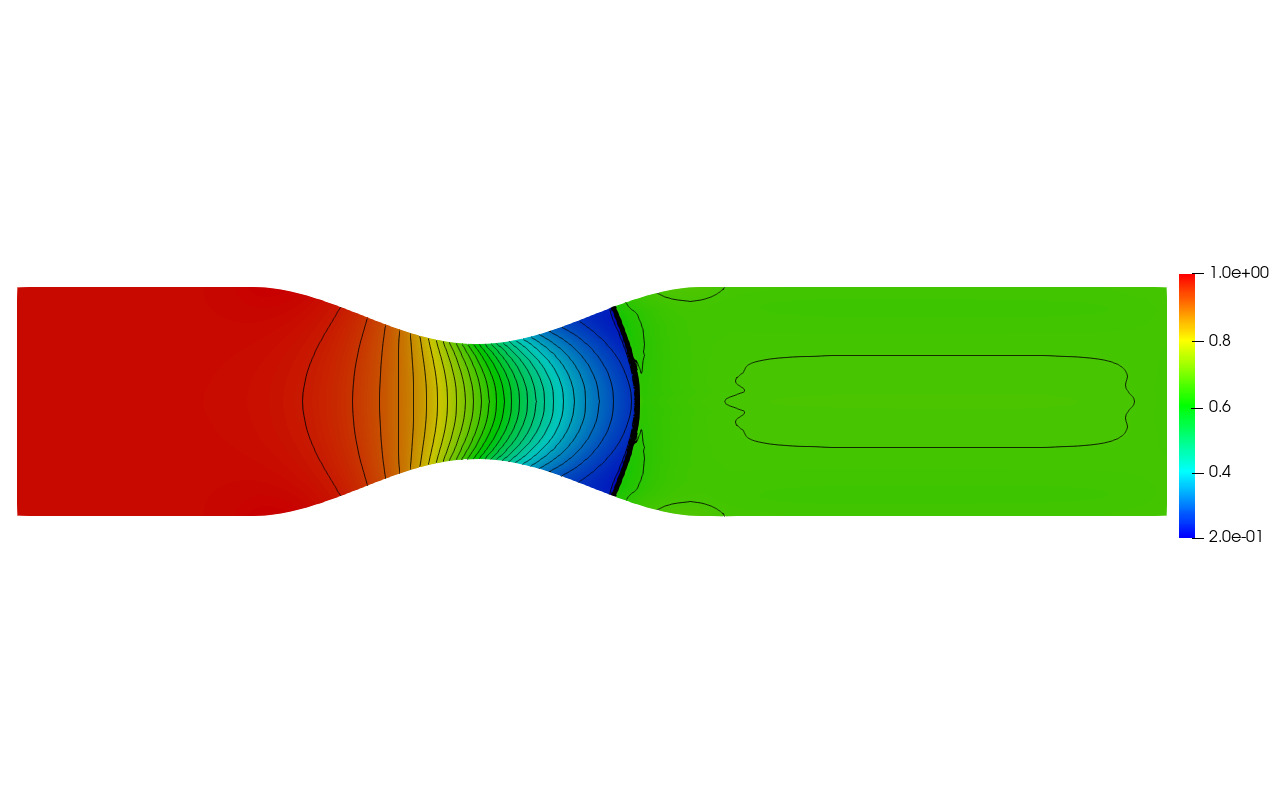}
		\subcaption{Density}
		\label{fig:SingleShockNozzleden}
	\end{subfigure}
	\begin{subfigure}[c]{0.7\textwidth}
		\includegraphics[trim={1cm 9cm 0cm 8.5cm},clip, width = \textwidth]{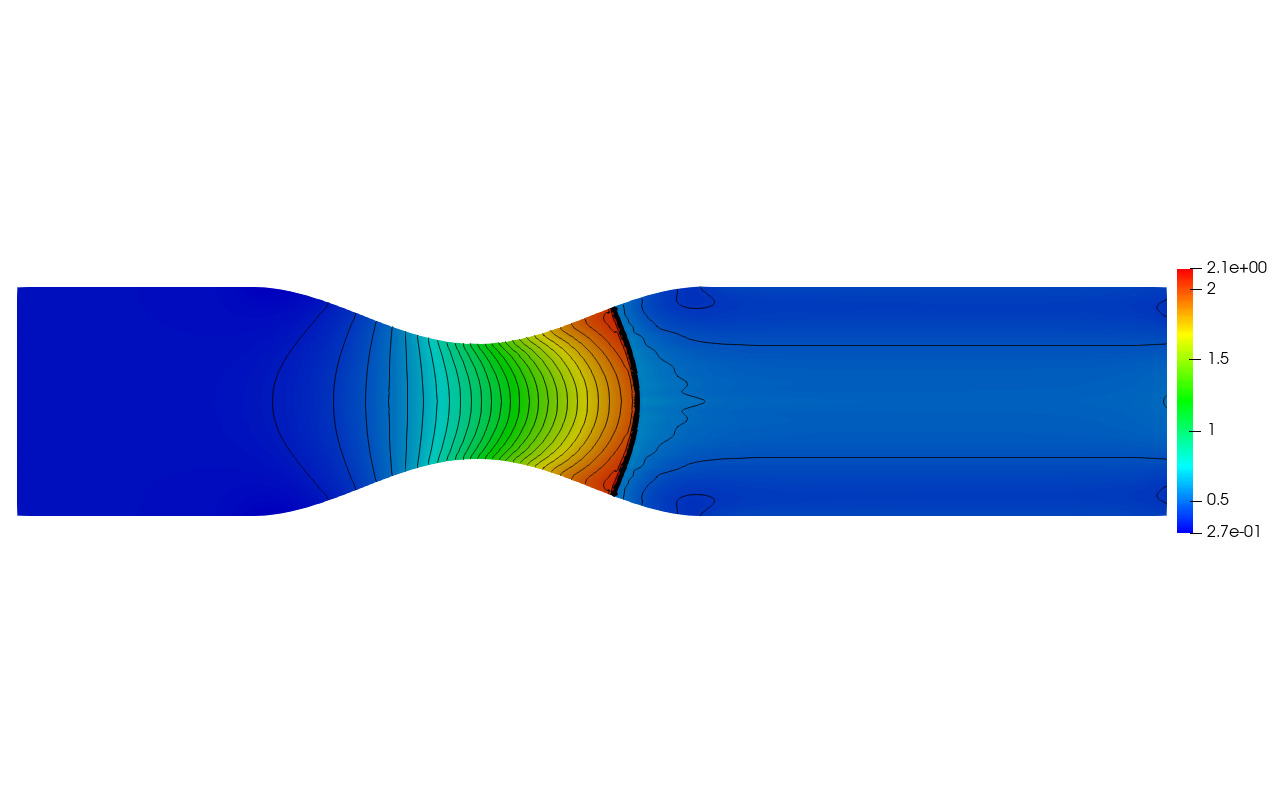}
		\subcaption{Mach number}
		\label{fig:SingleShockNozzlemach}
	\end{subfigure}
	\begin{subfigure}[c]{0.7\textwidth}
		\includegraphics[trim={1cm 9cm 0cm 8.5cm},clip, width = \textwidth]{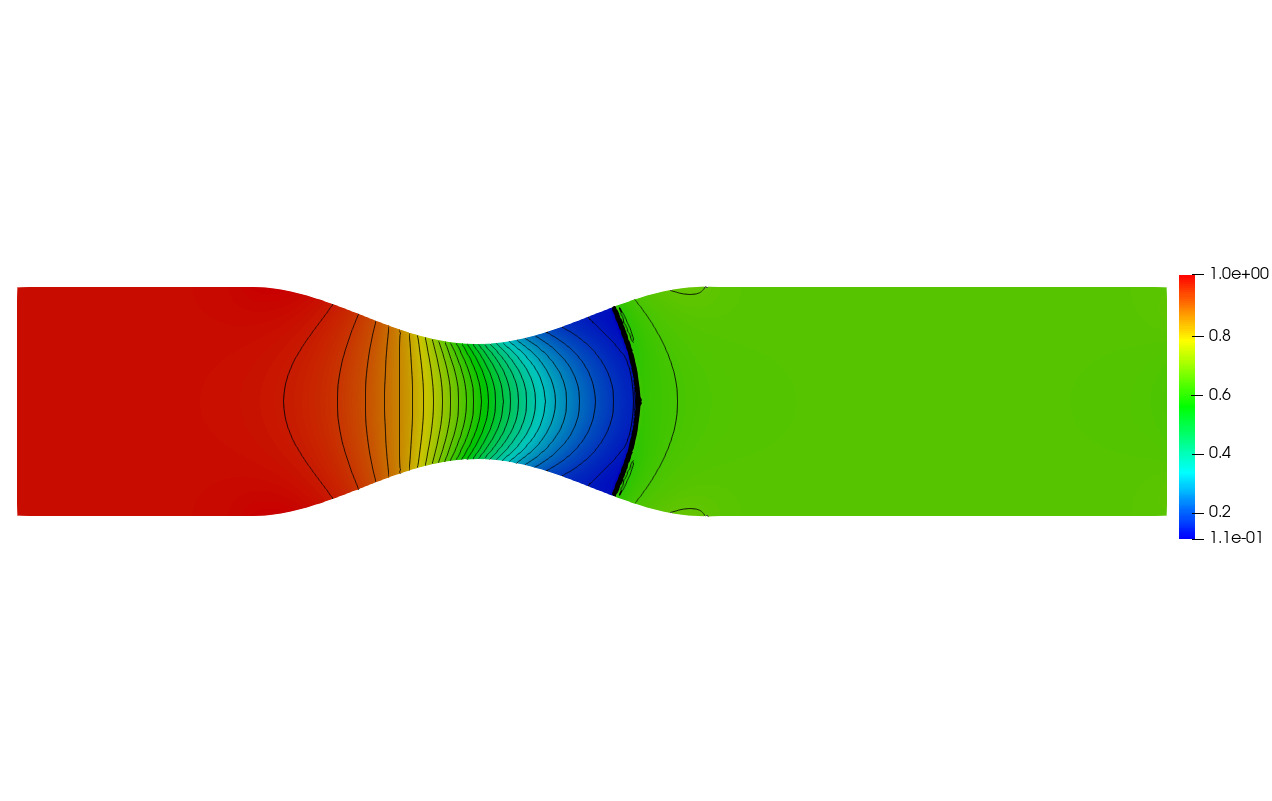}
		\subcaption{Pressure}
		\label{fig:SingleShockNozzlepres}
	\end{subfigure}
	\caption{Transonic nozzle ($M_\infty= 0.3$), steady-state MCL results obtained with $N_h = 95,329$ unknowns per component on a mesh consisting of $E_h = 189,024$ triangles.}
	\label{fig:SingleShockNozzle}
\end{figure}

The steady-state MCL results and the convergence history for this test
are shown in Fig. \ref{fig:02Nozzle} and 
Fig.~\ref{fig:02Nozzle_SSR},
respectively. In this low Mach number test, convergence to a stationary
solution is achieved only with $\omega = 0.5$ and with the adaptive
 strategy~\eqref{eq:aur}. In both cases,
the stopping criterion $r(u) \le 10^{-13}$ is met after $1000$ pseudo-time
steps. All other choices of $\omega$ result in
stagnation or even divergence (see Fig.~\ref{fig:02Nozzle_ur}).
Even for small CFL numbers, which
introduce strong implicit underrelaxation, convergence is not guaranteed
without additional explicit underrelaxation (not shown here).
Similarly to the previous test case, the adaptive underrelaxation
strategy exhibits the fastest convergence rate. The impact of the
CFL number is illustrated by the evolution of steady-state residuals
in Fig.~\ref{fig:02Nozzle_CFL}. The convergence rates for $\mathrm{CFL}=10^4$
and $\mathrm{CFL}=10^5$ are similar and superior to that for
$\mathrm{CFL}=10^3$. 

The second setup induces a transonic flow. 
For better comparison with the literature~\cite{hartmann2002}, we prescribe the $M_\infty = 0.3$ free stream values using the pressure $p_\infty = 1$ instead of $p_\infty = 1/\gamma$. The free stream velocity corresponding to
$\rho_\infty = 1$ is given by
$\vel_\infty = (\sqrt{\gamma}M_\infty , 0)^\top$.
Furthermore, we prescribe $p_{\mathrm{out}} = \frac{2}{3}$ instead of the free stream pressure at the subsonic outlet.
The resulting flow accelerates to $M \approx 2.1$ and develops a sonic shock, behind which it decelerates to subsonic speeds again.

\begin{figure}[H]
	\centering
	\begin{subfigure}[c]{0.40\textwidth}
		\includegraphics[trim={0.9cm 0.3cm 2.15cm 1cmm},clip, width = 0.95\textwidth]{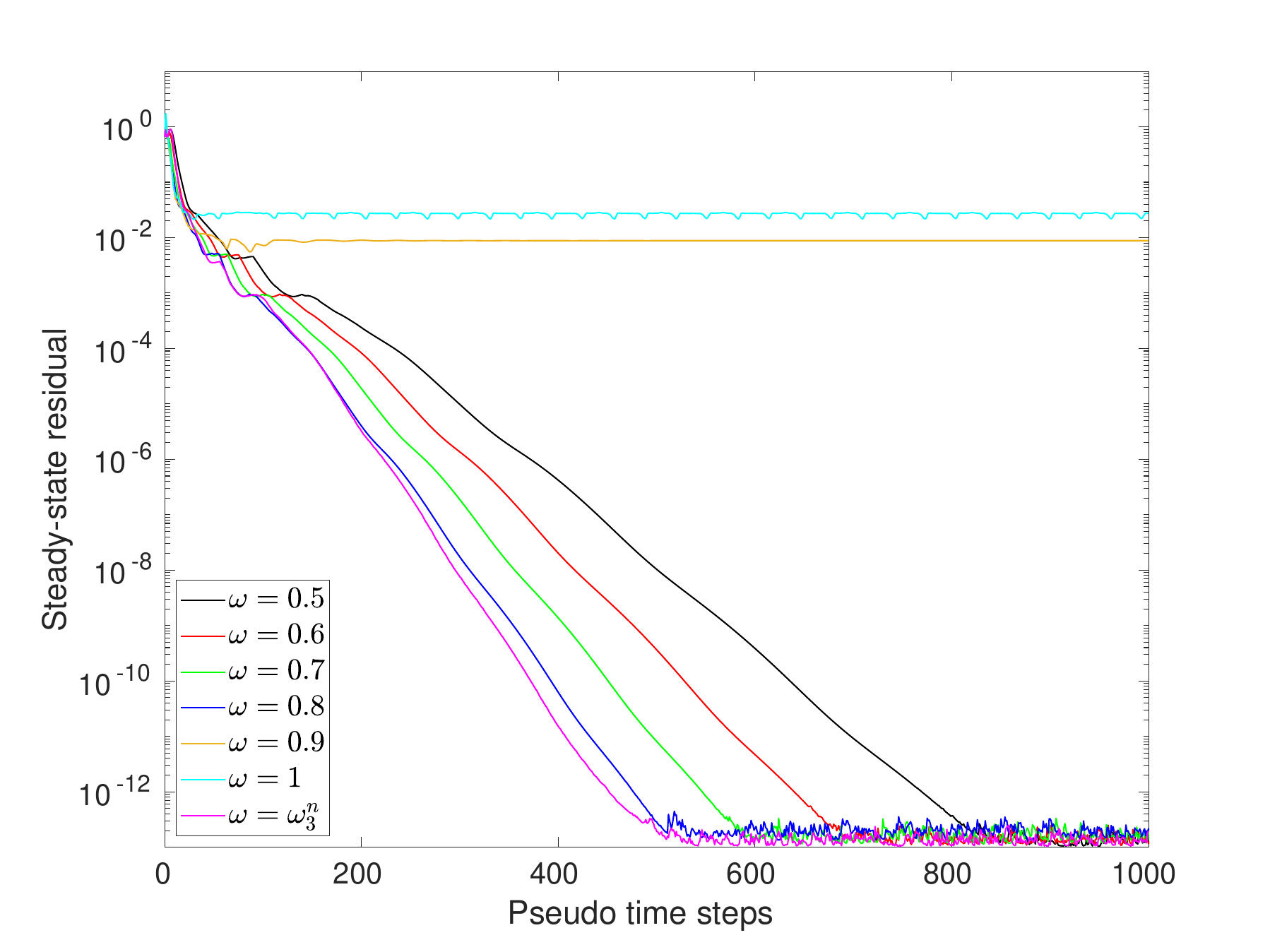}
		\subcaption{Different underrelaxation factors, $\mathrm{CFL} = 10^4$.}
		\label{fig:SingelShockNozzle_ur}
	\end{subfigure}
	\begin{subfigure}[c]{0.40\textwidth}
		\includegraphics[trim={0.9cm 0.3cm 2.15cm 1cmm},clip, width = 0.95\textwidth]{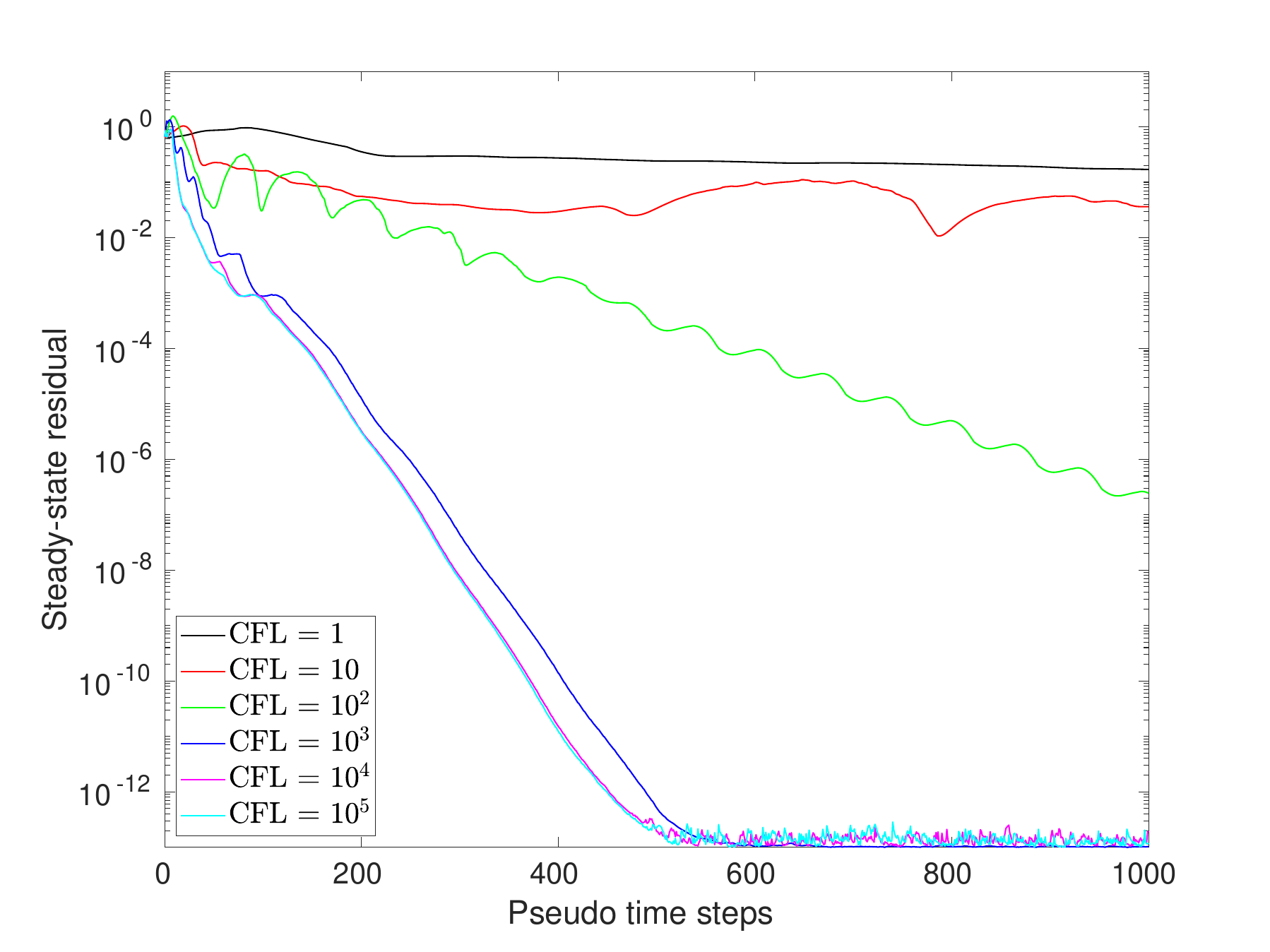}
		\subcaption{Different CFL numbers, $\omega = \omega_3^n$.}
		\label{fig:SingelShockNozzle_CFL}
	\end{subfigure}
	\caption{Transonic nozzle ($M_\infty = 0.3$), steady-state convergence history for a long-time MCL simulation on an unstructured triangular mesh with $N_h = 24,037$ nodes and $E_h = 47,256$ cells.}
	\label{fig:SingelShockNozzle_SSR}
\end{figure}

As we observe in Fig.~\ref{fig:SingelShockNozzle_ur}, solvers that use $\omega\in\{0.5, 0.6, 0.7, 0.8, \omega^n_3\}$ do converge, while the residuals stagnate after a short initial phase for 
$\omega = 1$ and $\omega = 0.9$. The algorithm using the
adaptive strategy~\eqref{eq:aur} requires about 500 pseudo-time steps
to meet the stopping criterion $r(u) \le 10^{-13}$. As in the previous examples,
convergence becomes faster as the CFL number is increased (see Fig.~\ref{fig:SingelShockNozzle_CFL}).

\begin{figure}[t!]
	\centering
	\begin{subfigure}[c]{0.7\textwidth}
		\includegraphics[trim={0cm 9cm 0cm 8.5cm},clip, width = \textwidth]{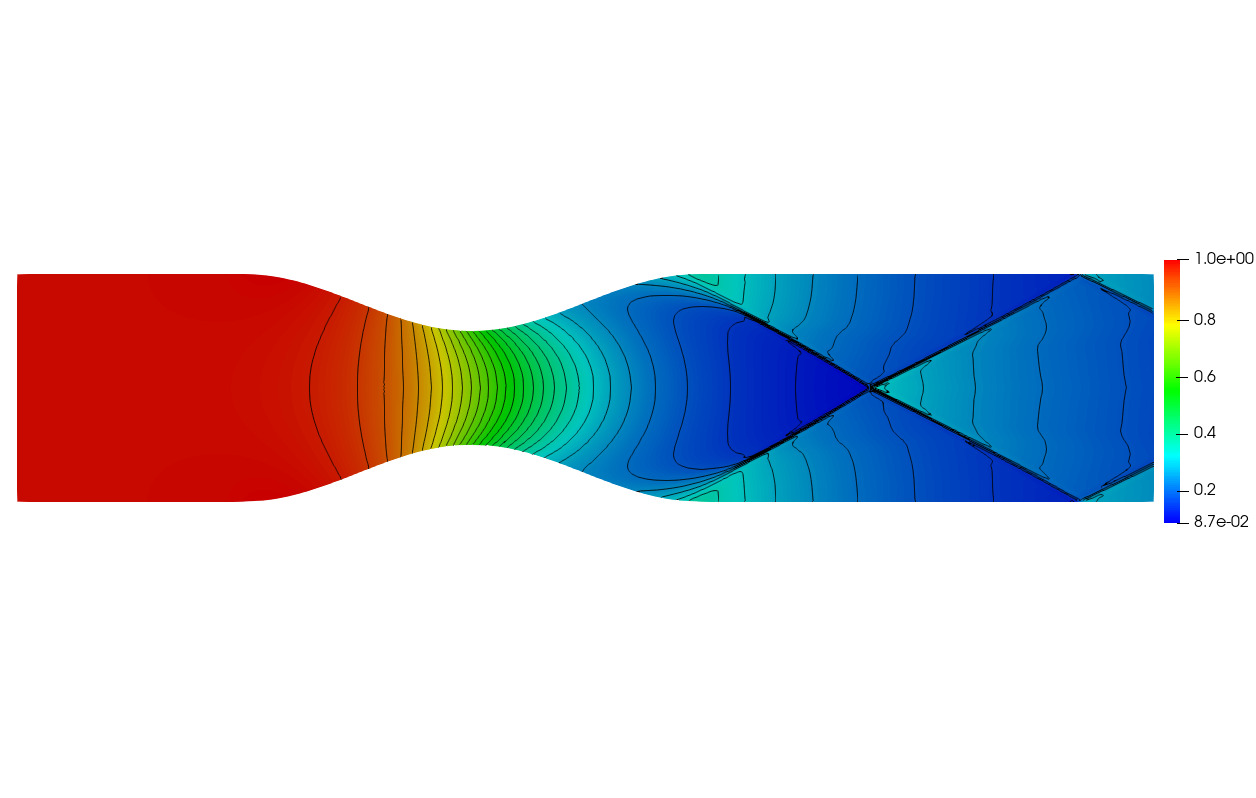}
		\subcaption{Density}
		\label{fig:TransonicNozzleden}
	\end{subfigure}
	\begin{subfigure}[c]{0.7\textwidth}
		\includegraphics[trim={0cm 9cm 0cm 8.5cm},clip, width = \textwidth]{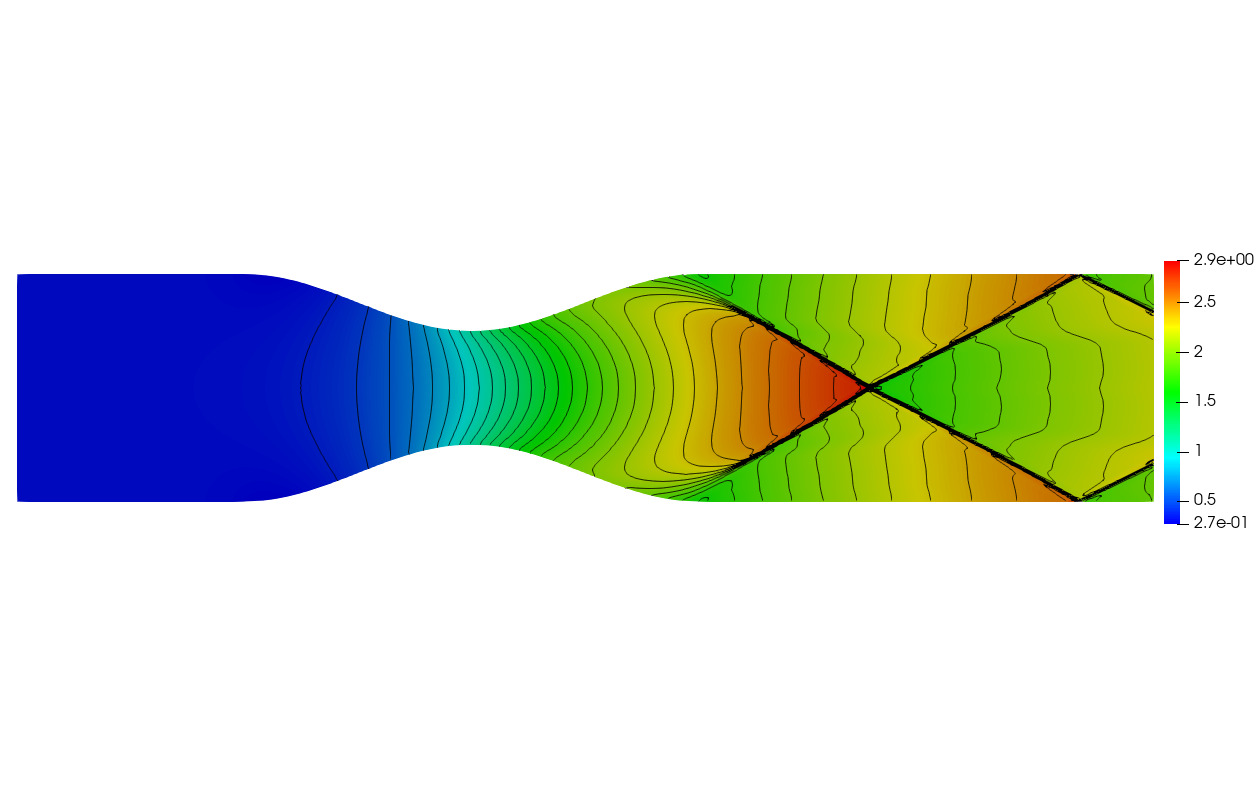}
		\subcaption{Mach number}
		\label{fig:TransonicNozzlemach}
	\end{subfigure}
	\begin{subfigure}[c]{0.7\textwidth}
		\includegraphics[trim={0cm 9cm 0cm 8.5cm},clip, width = \textwidth]{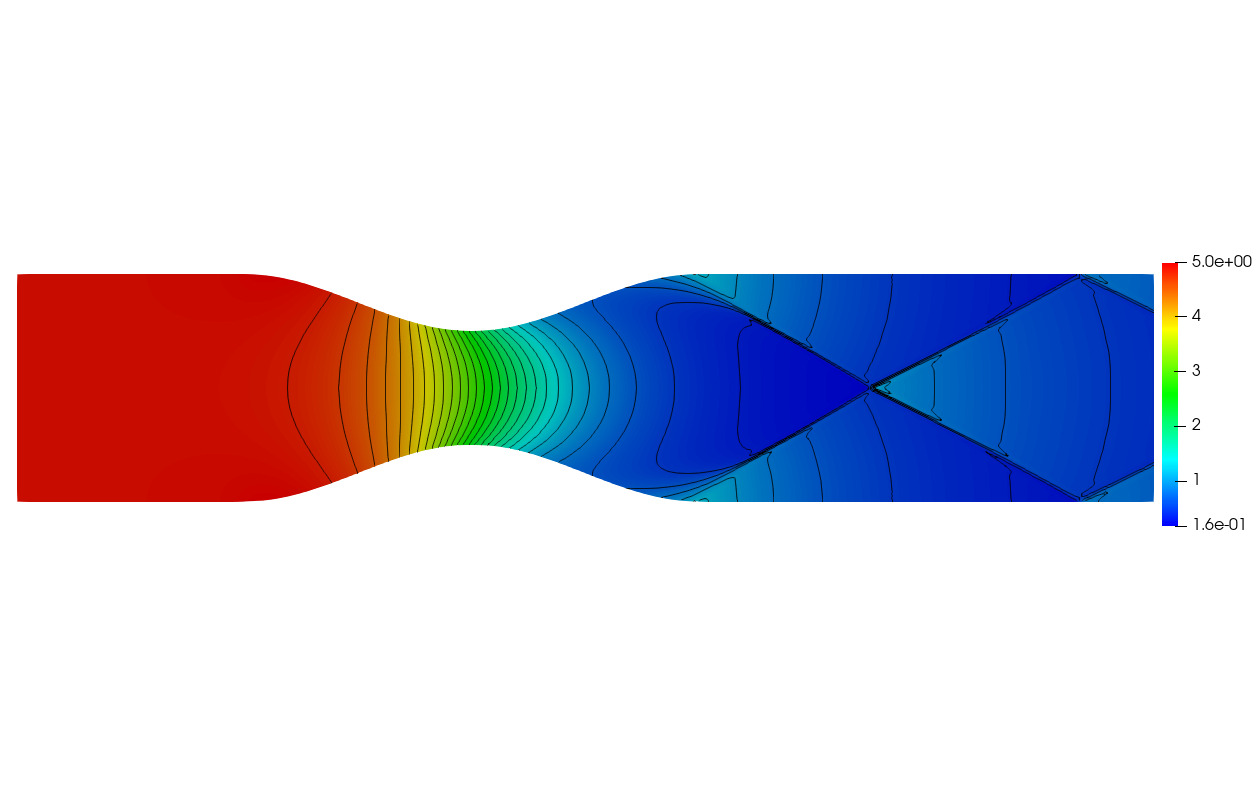}
		\subcaption{Pressure}
		\label{fig:TransonicNozzlepres}
	\end{subfigure}
	\caption{Transonic nozzle ($M_\infty= 0.8$), steady-state MCL results obtained with $N_h = 95,329$ unknowns per component on a mesh consisting of $E_h = 189,024$ triangles.}
	\label{fig:TransonicNozzle}
\end{figure}

\begin{figure}[h!]
	\centering
	\begin{subfigure}[c]{0.40\textwidth}
		\includegraphics[trim={0.9cm 0.3cm 2.15cm 1cmm},clip, width = 0.95\textwidth]{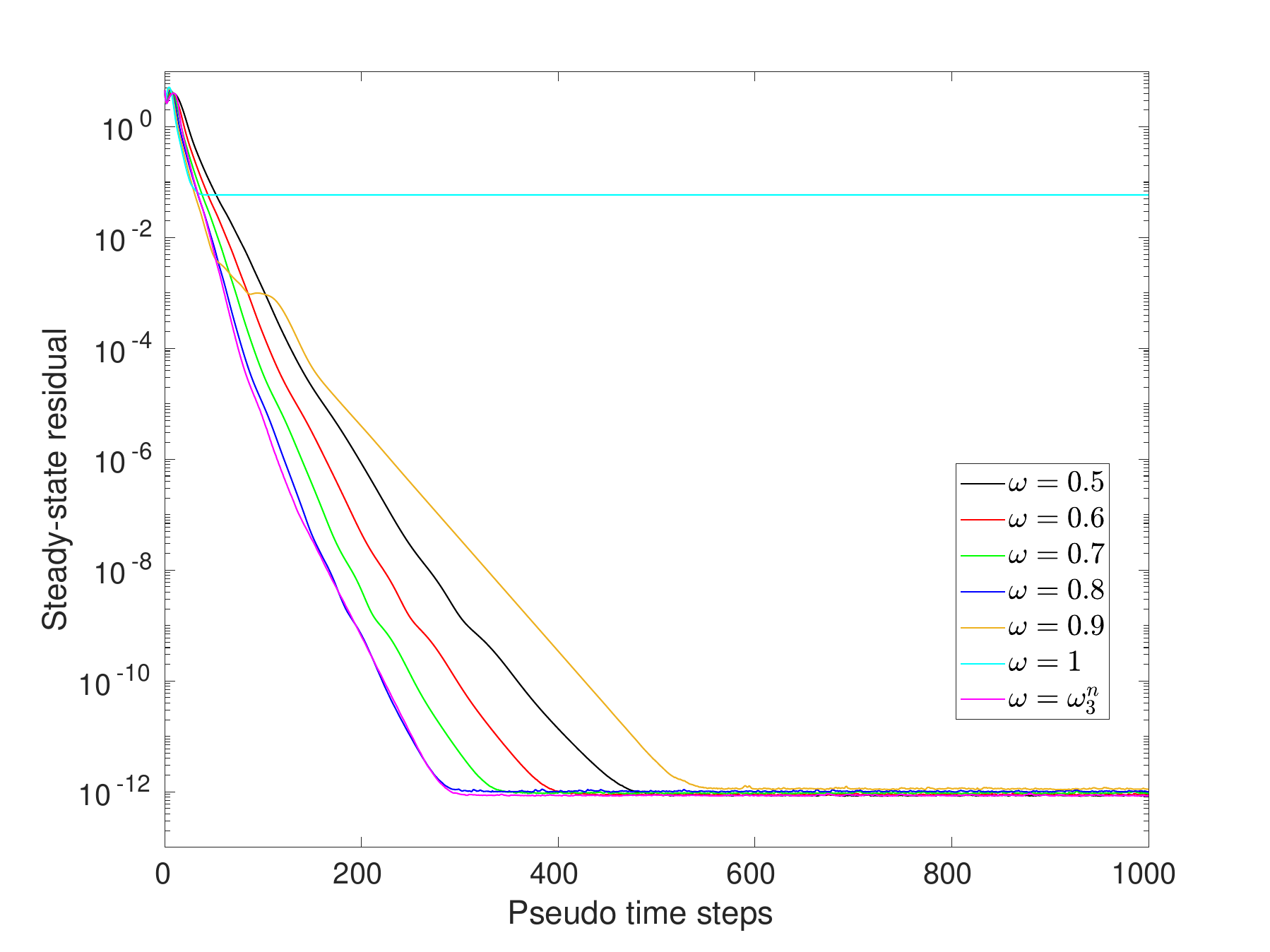}
		\subcaption{Different underrelaxation factors, $\mathrm{CFL} = 10^4$.}
		\label{fig:TransonicNozzle_ur}
	\end{subfigure}
	\begin{subfigure}[c]{0.40\textwidth}
		\includegraphics[trim={0.9cm 0.3cm 2.15cm 1cmm},clip, width = 0.95\textwidth]{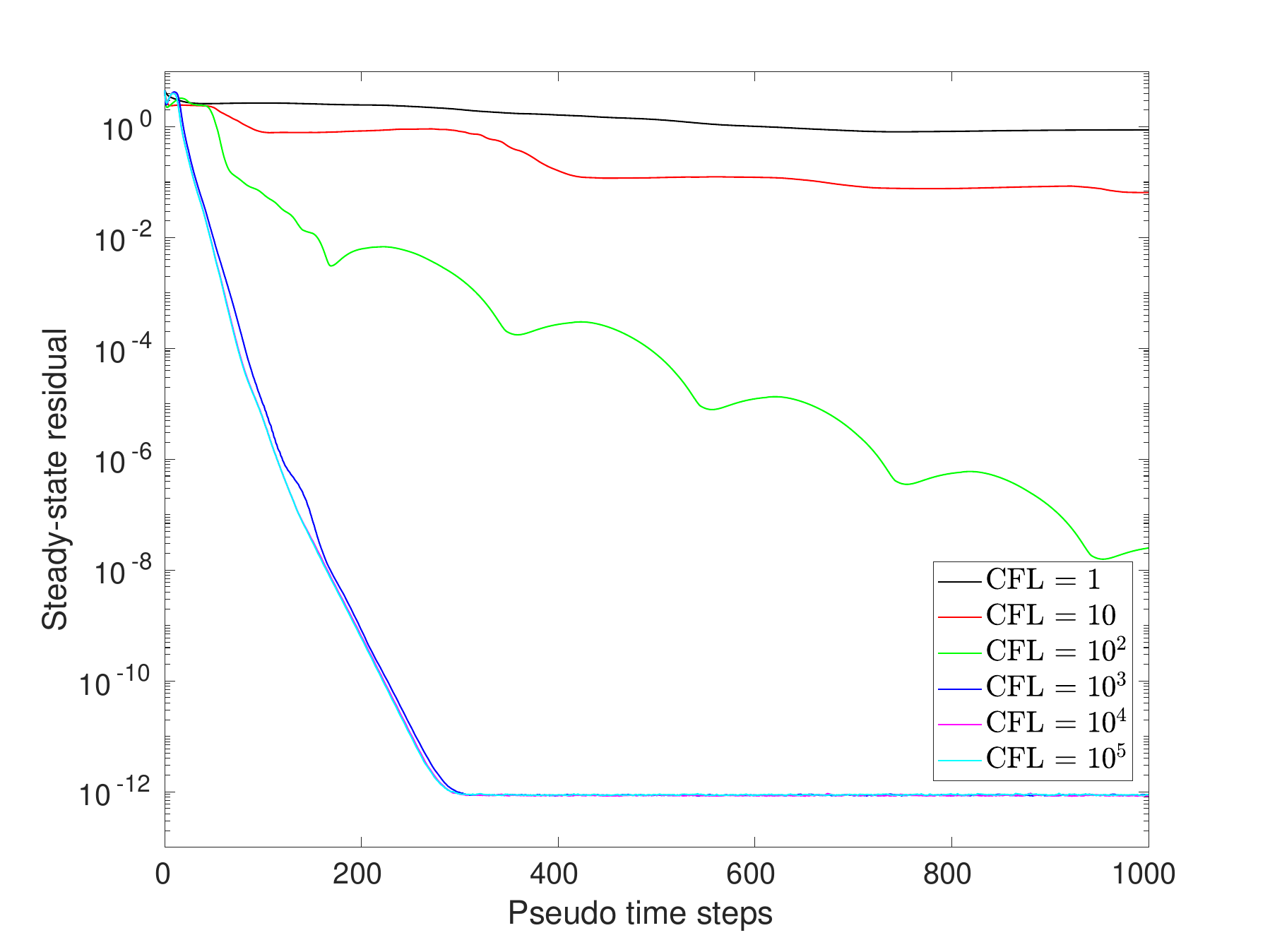}
		\subcaption{Different CFL numbers,  $\omega = \omega_3^n$.}
		\label{fig:TransonicNozzle_CFL}
	\end{subfigure}
	\caption{Transonic nozzle ($M_\infty = 0.8$), steady-state convergence history for a long-time MCL simulation on an unstructured triangular mesh with $N_h = 24,037$ nodes and $E_h = 47,256$ cells.}
	\label{fig:TransonicNozzle_SSR}
\end{figure}

The last setup corresponds to a transonic flow with a supersonic outlet~\cite{gurris2009}. 
We define the free stream values using~\eqref{eq:freestream} with $M_\infty = 0.8$.
The flow is accelerated in the throat of the nozzle and stays supersonic. It features diamond-shaped shocks caused by reflections from the upper and lower walls.
For this transonic flow setup, the steady-state computation fails to converge only if $\omega=1$ is employed. The adaptive strategy~\eqref{eq:aur} meets the stopping criterion $r(u) \le 10^{-13}$ after 300 pseudo-time steps. As we can see in Fig.~\ref{fig:TransonicNozzle}, the use of $\omega = 0.9$ gives rise to small fluctuations around the threshold $10^{-3}$.
Figure~\ref{fig:TransonicNozzle_CFL} shows that
$\mathrm{CFL}\geq10^3$ is needed to achieve satisfactory
convergence behavior.

\subsection{NACA 0012 airfoil}

To further exemplify the typical convergence behavior of our scheme, we simulate external flows over a NACA 0012 airfoil. The upper and lower surface of this airfoil are defined by \cite{gurris2009, kuzmin2012b}
\begin{equation*}
	\Gamma_{\pm} = \left\{(x, f_\pm(x))^\top \in \R^2: x\in[0, 1.00893]\right\},
\end{equation*}
where
\begin{equation*}
	f_\pm(x) = \pm 0.6(0.2969\sqrt{x} - 0.126x - 0.3516x^2- 0.1015x^4).
\end{equation*}
The outer boundary of the computational domain is a circle of radius 10 centered at the tip of the airfoil. The steady-state flow pattern depends on the Mach number and on the inclination angle.

\begin{figure}[h!]
	\centering
	\begin{subfigure}[c]{0.49\textwidth}
		\includegraphics[trim={1cm 2cm 0.1cm 2cm},clip, width = 0.99\textwidth]{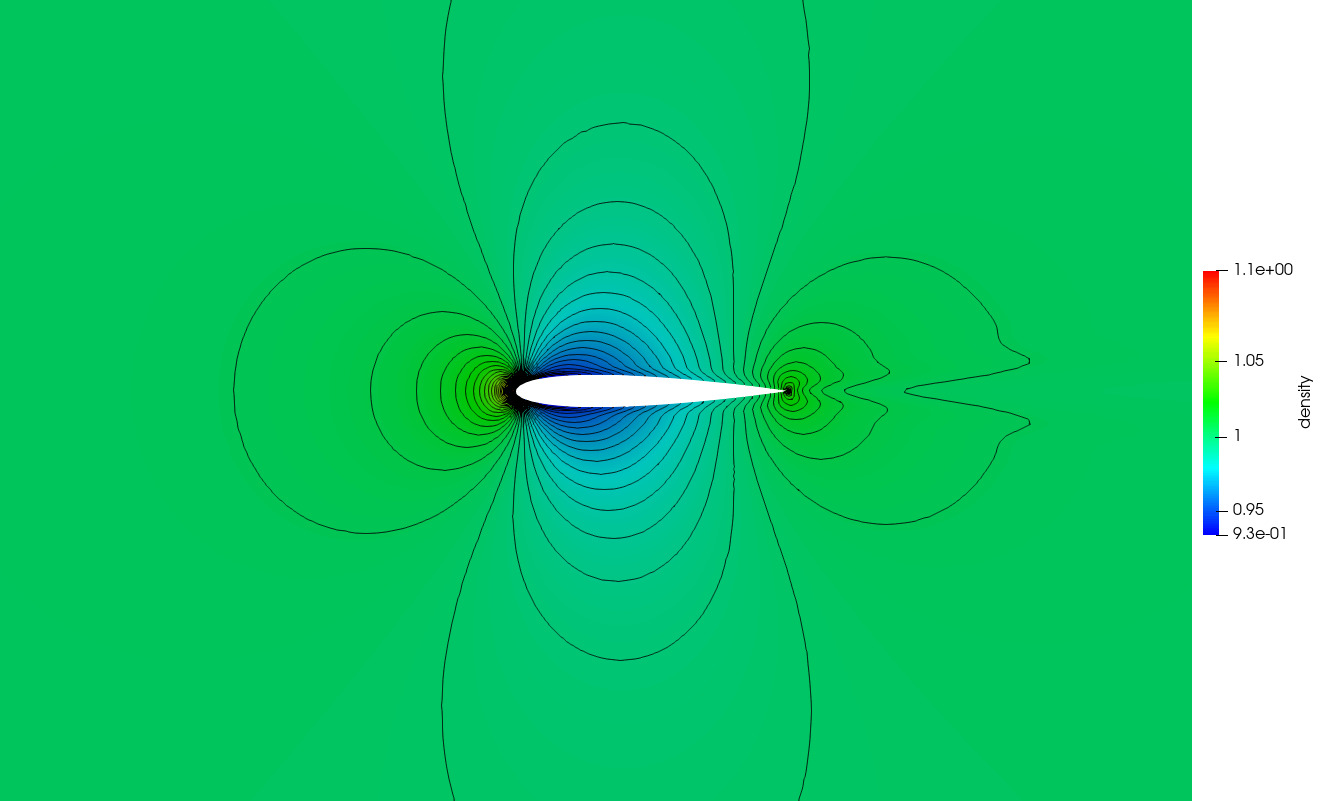}
		\subcaption{Density}
		\label{fig:SubSonicNacaden}
	\end{subfigure}
	\begin{subfigure}[c]{0.49\textwidth}
		\includegraphics[trim={1cm 2cm 0.1cm 2cm},clip,width = 0.99\textwidth]{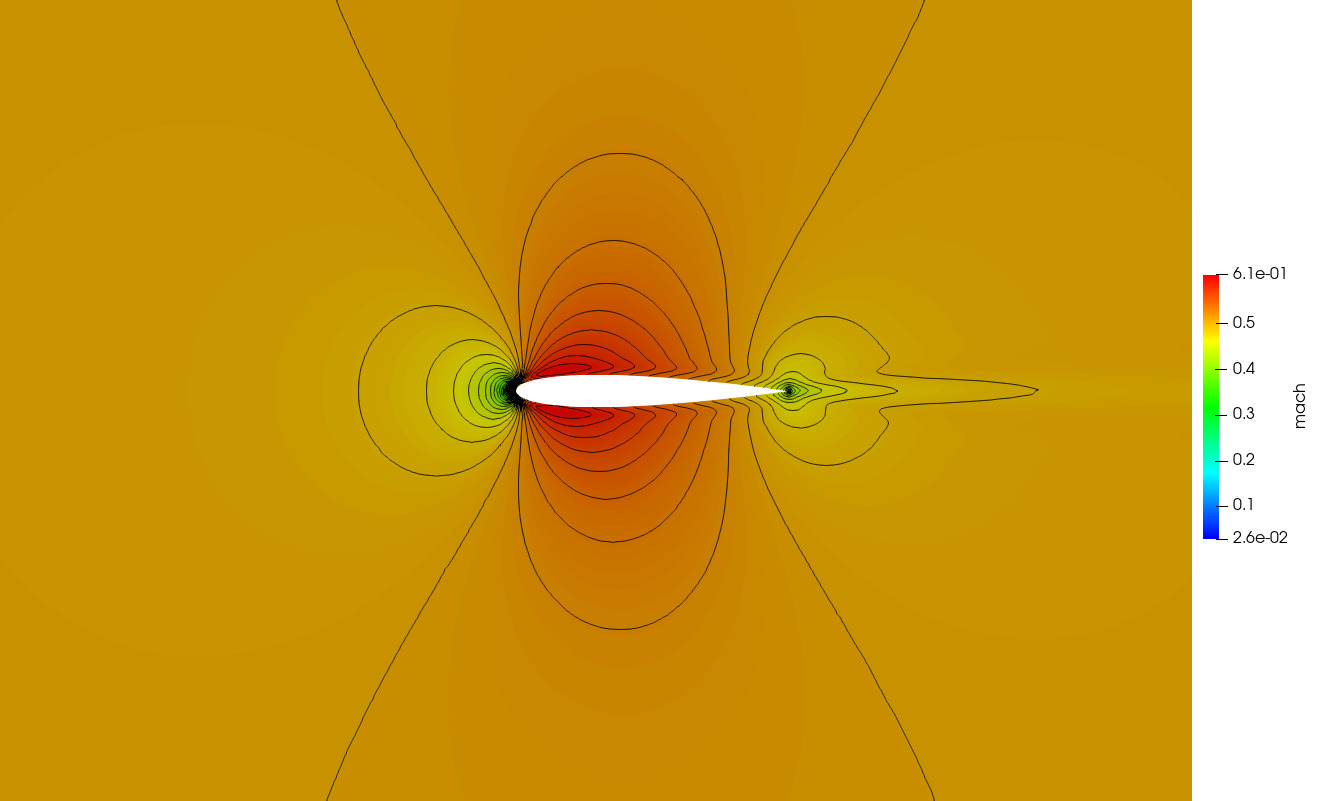}
		\subcaption{Mach number}
		\label{fig:SubSonicNacamach}
	\end{subfigure}
	\begin{subfigure}[c]{0.49\textwidth}
		\includegraphics[trim={1cm 2cm 0.1cm 2cm},clip, width = 0.99\textwidth]{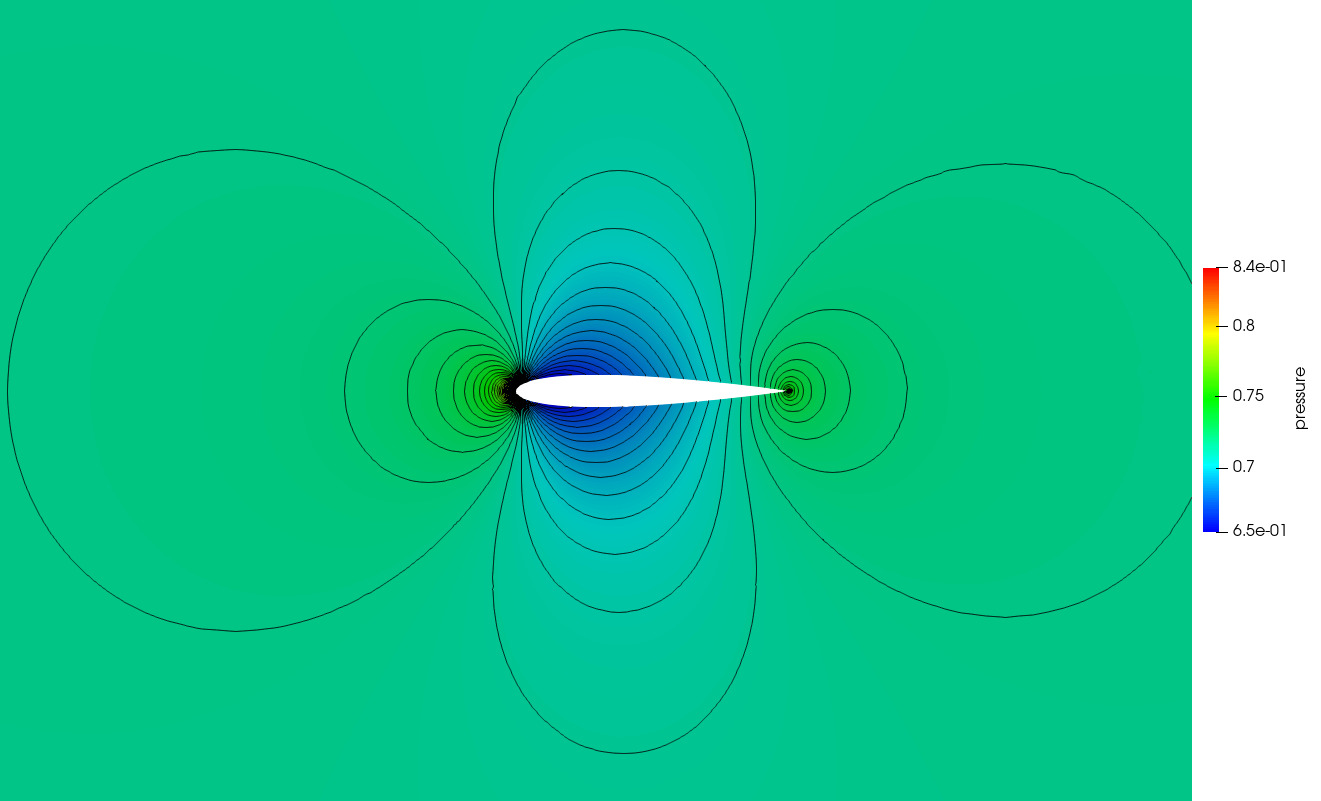}
		\subcaption{Pressure}
		\label{fig:SubSonicNacapres}
	\end{subfigure}
	\caption{Subsonic NACA 0012 airfoil ($M_\infty= 0.5$, $\alpha = 0^\circ$), steady-state MCL results obtained with $N_h = 112,242$ unknowns per component on a mesh consisting of $E_h = 223,424$ triangles.}
	\label{fig:SubsonicNaca}
\end{figure}

\begin{figure}[h!]
	\centering
	\begin{subfigure}[c]{0.40\textwidth}
		\includegraphics[trim={0.9cm 0.3cm 2.15cm 1cmm},clip, width = 0.95\textwidth]{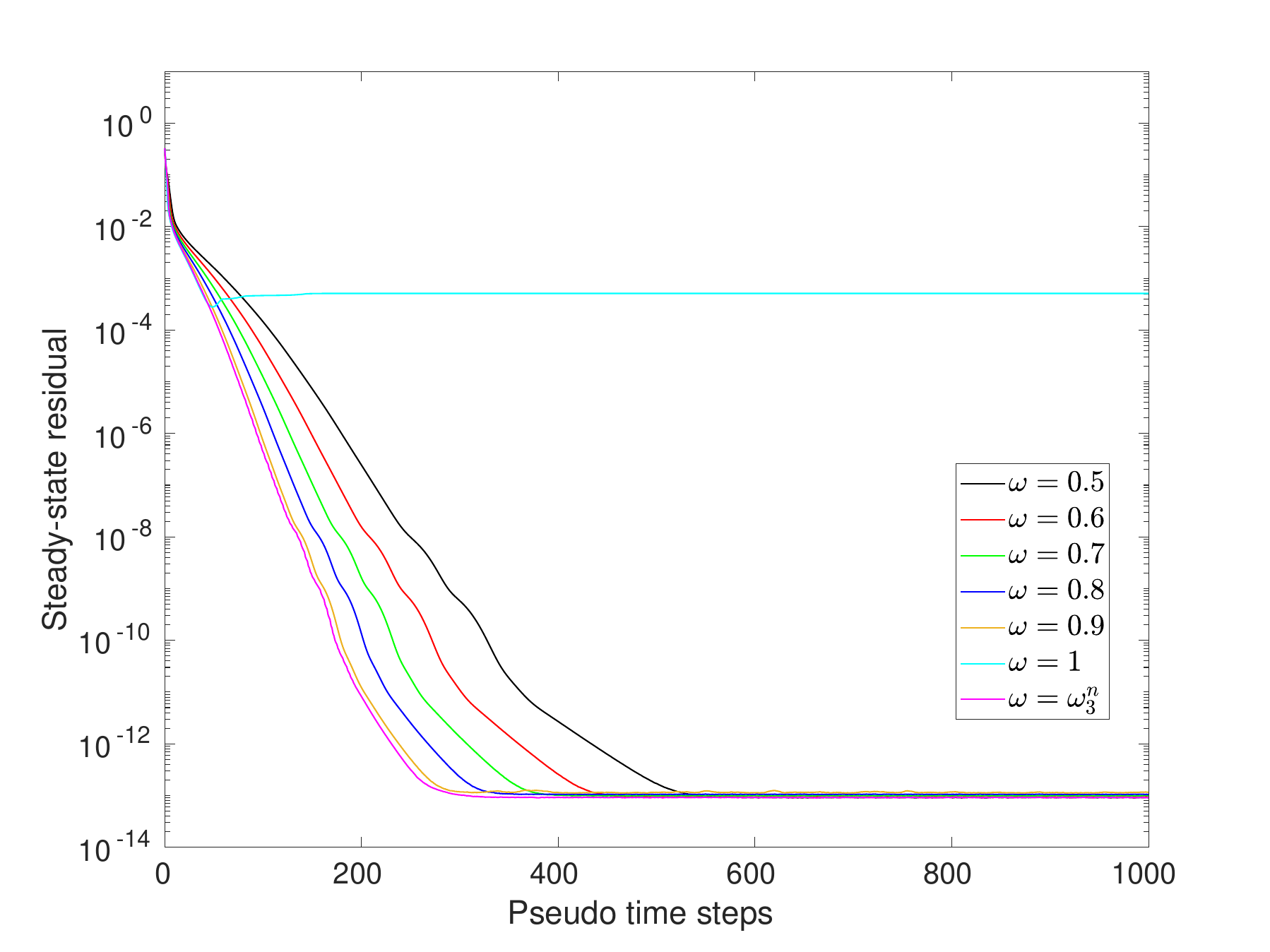}
		\subcaption{Different underrelaxation factors, $\mathrm{CFL} = 10^4$.}
		\label{fig:SubsonicNaca_ur}
	\end{subfigure}
	\begin{subfigure}[c]{0.40\textwidth}
		\includegraphics[trim={0.9cm 0.3cm 2.15cm 1cmm},clip, width = 0.95\textwidth]{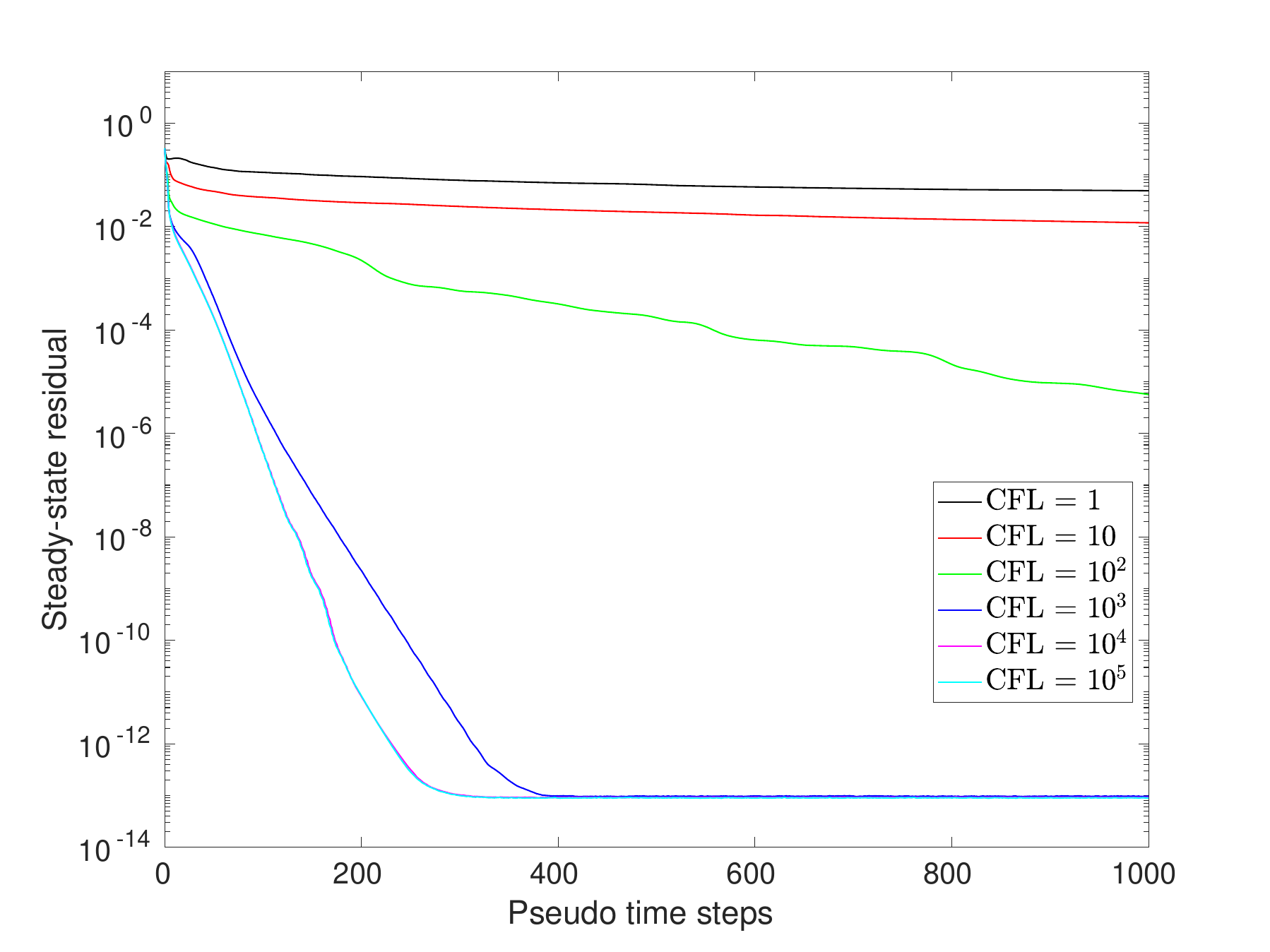}
		\subcaption{Different CFL numbers, $\omega = \omega_3^n$.}
		\label{fig:SubsonicNaca_CFL}
	\end{subfigure}
	\caption{Subsonic NACA 0012 airfoil ($M_\infty= 0.5$, $\alpha = 0^\circ$), steady-state convergence history for a long-time MCL simulation on an unstructured triangular mesh with $N_h = 28,193$ nodes and $E_h = 55,856$ cells.}
	\label{fig:SubsonicNaca_SSR}
\end{figure}

\begin{figure}[h!]
	\centering
	\begin{subfigure}[c]{0.49\textwidth}
		\includegraphics[trim={1cm 2cm 0.1cm 2cm},clip, width = 0.99\textwidth]{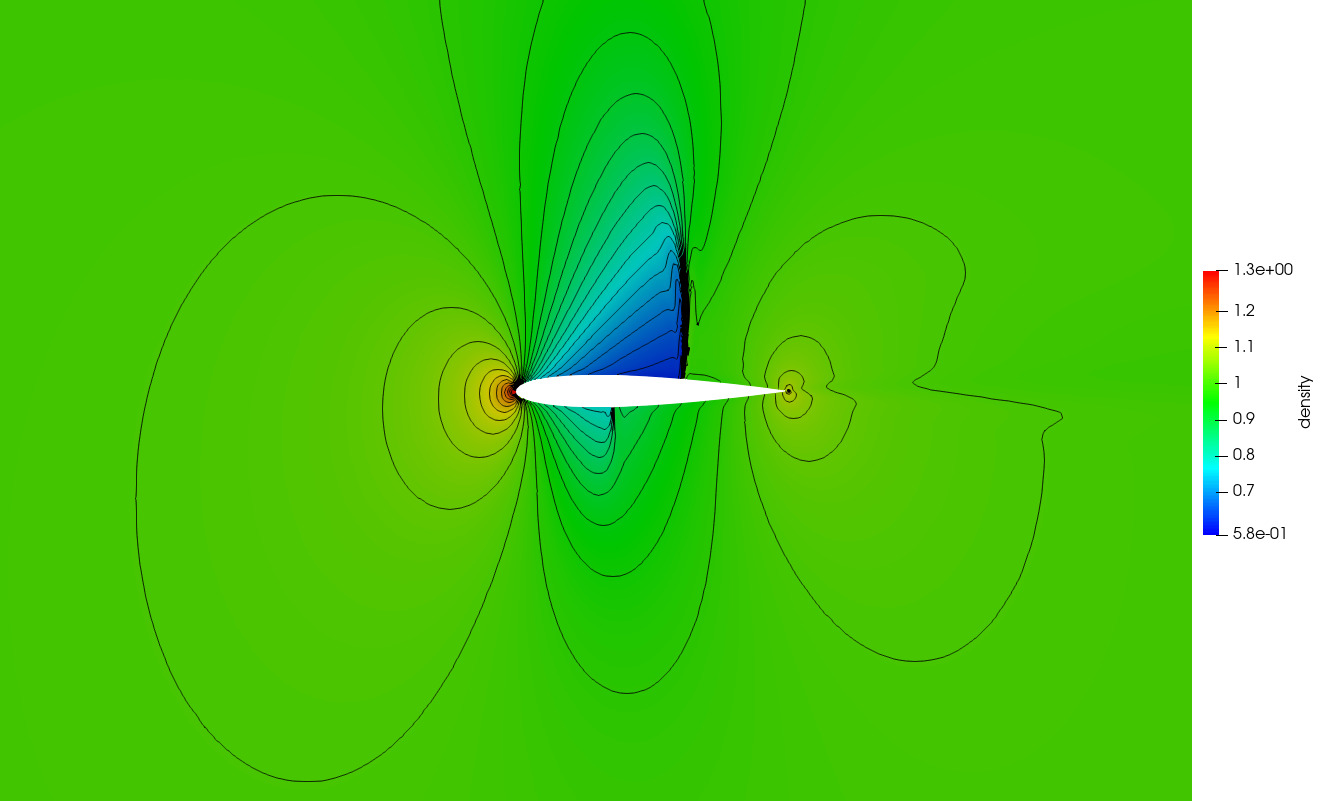}
		\subcaption{Density}
		\label{fig:125Nacaden}
	\end{subfigure}
	\begin{subfigure}[c]{0.49\textwidth}
		\includegraphics[trim={1cm 2cm 0.1cm 2cm},clip, width = 0.99\textwidth]{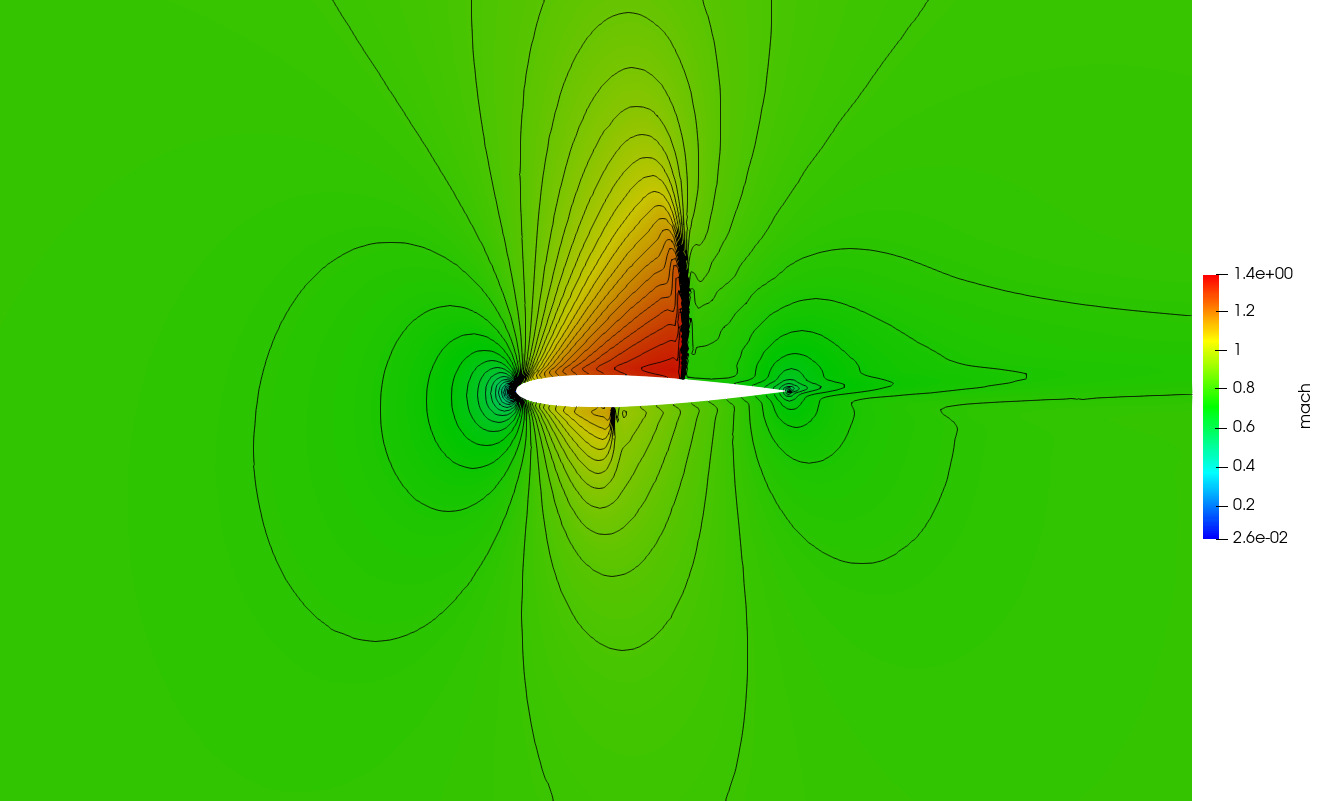}
		\subcaption{Mach number}
		\label{fig:125Nacamach}
	\end{subfigure}
	\begin{subfigure}[c]{0.49\textwidth}
		\includegraphics[trim={1cm 2cm 0.1cm 2cm},clip, width = 0.99\textwidth]{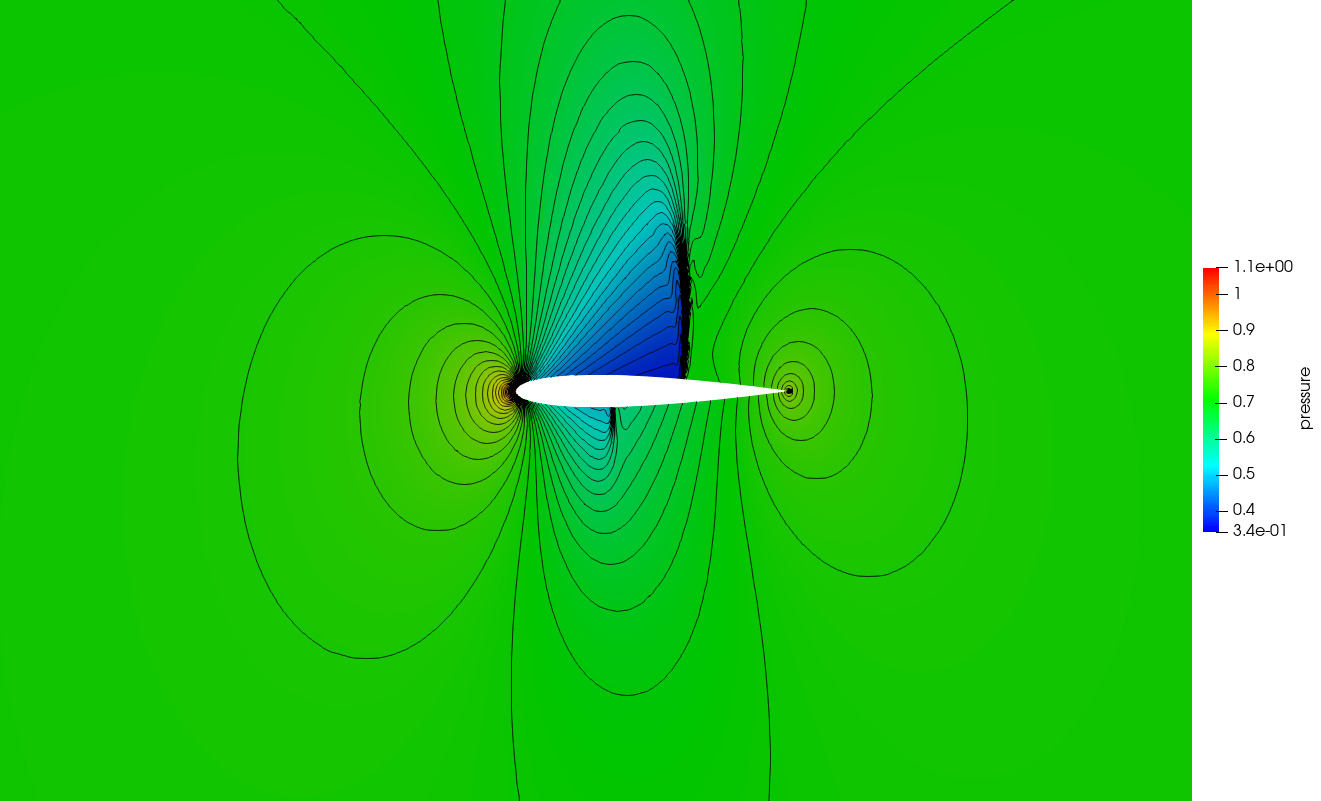}
		\subcaption{Pressure}
		\label{fig:125Nacapres}
	\end{subfigure}
	\caption{Transonic NACA 0012 airfoil at $\alpha = 1.25^\circ$ angle ($M_\infty= 0.8$), steady-state MCL results obtained with $N_h = 112,242$ unknowns per component on a mesh consisting of $E_h = 223,424$ triangles.}
	\label{fig:125Naca}
\end{figure}

\begin{figure}[h!]
	\centering
	\begin{subfigure}[c]{0.40\textwidth}
		\includegraphics[trim={0.9cm 0.3cm 2.15cm 1cmm},clip, width = 0.95\textwidth]{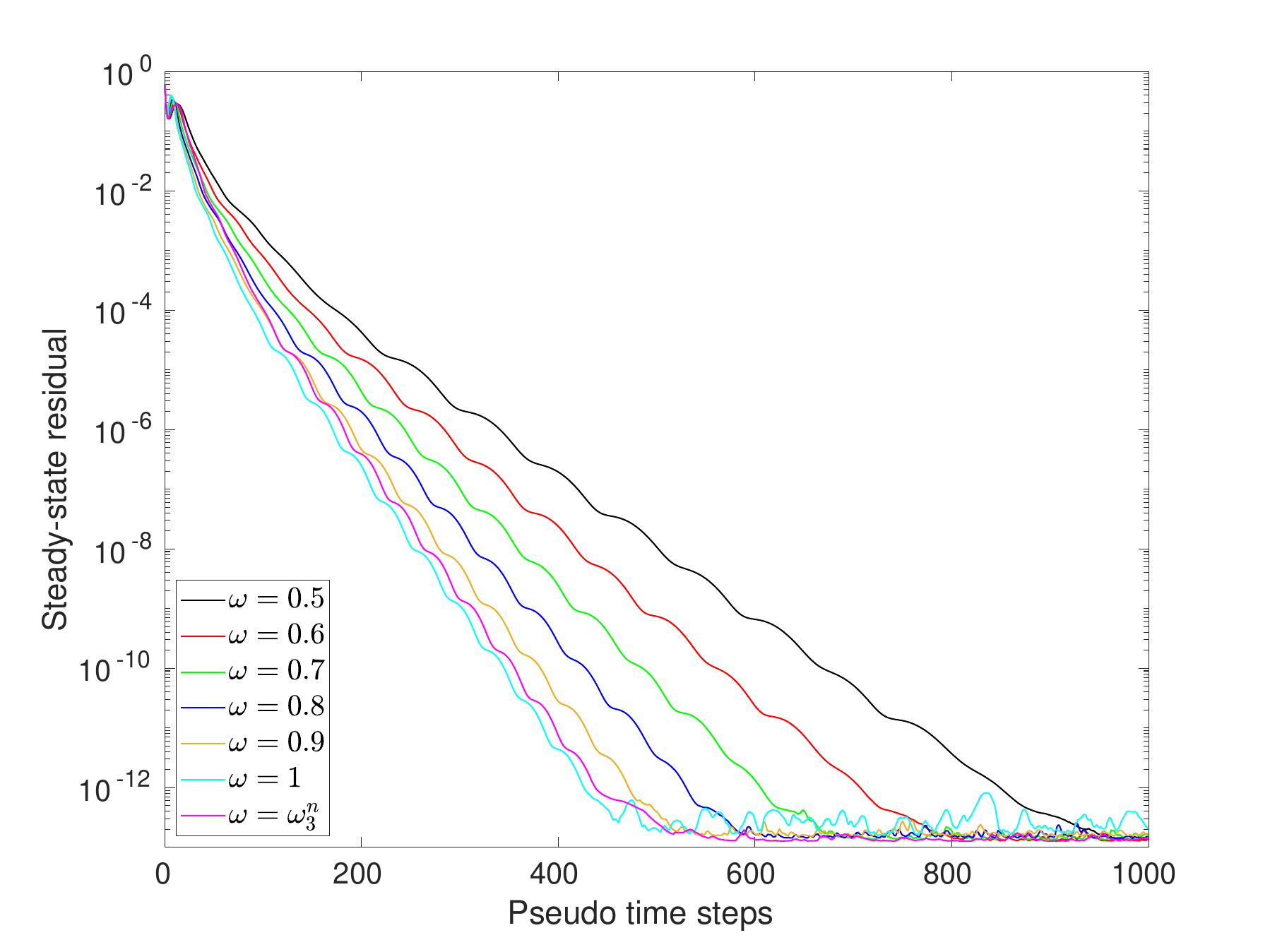}
		\subcaption{Different underrelaxation factors, $\mathrm{CFL} = 10^4$.}
		\label{fig:125Naca_ur}
	\end{subfigure}
	\begin{subfigure}[c]{0.40\textwidth}
		\includegraphics[trim={0.9cm 0.3cm 2.15cm 1cmm},clip, width = 0.95\textwidth]{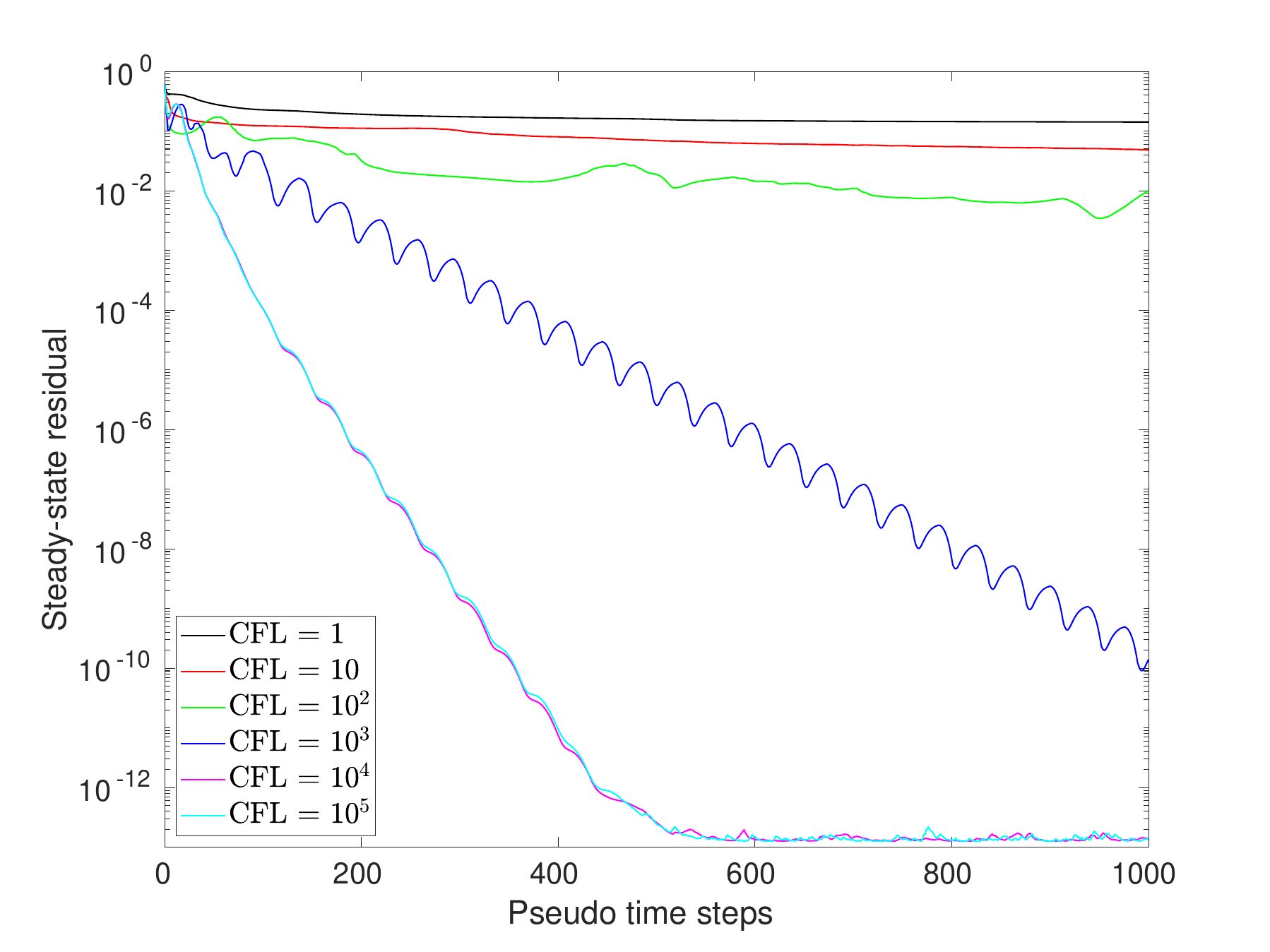}
		\subcaption{Different CFL numbers, $\omega = \omega_3^n$.}
		\label{fig:125Naca_CFL}
	\end{subfigure}
	\caption{Transonic NACA 0012 airfoil ($M_\infty = 0.8$, $\alpha = 1.25^\circ$), steady-state convergence history for a long-time MCL simulation on an unstructured triangular mesh with $N_h = 28,193$ nodes and $E_h = 55,856$ cells.}
	\label{fig:125Naca_SSR}
\end{figure}

In the first setup, we prescribe a $M_\infty= 0.5$ free stream flow with the angle of attack $\alpha = 0^\circ$. The steady state solutions and the convergence history are shown in Fig.~\ref{fig:SubsonicNaca} and Fig.~\ref{fig:SubsonicNaca_SSR}, respectively. In the second series of experiments, we set $M_\infty= 0.8$ and $\alpha =1.25^\circ$. The results are shown in Figs.~\ref{fig:125Naca} and~\ref{fig:125Naca_SSR}. As expected, relaxation is needed to achieve convergence in the stiff subsonic test. Our adaptive relaxation strategy performs slightly better than the best manual choice of a constant relaxation factor ($\omega=0.9$). In the transonic test, convergence is achieved for all values of $\omega$. The best convergence rates are obtained with $\omega=1.0$ and $\omega=\omega_3^n$. For both choices of the Mach number $M_\infty$, steady-state computations with $\omega=\omega_3^n$ become faster as the value of the CFL number is increased. 

\subsection{Mach 20 bow shock}
\label{sec:Mach20}

\begin{figure}[h!]
	\centering
	\begin{subfigure}[c]{0.25\textwidth}
		\includegraphics[trim={17cm 0cm 13.7cm 0cm},clip, width = 0.95\textwidth]{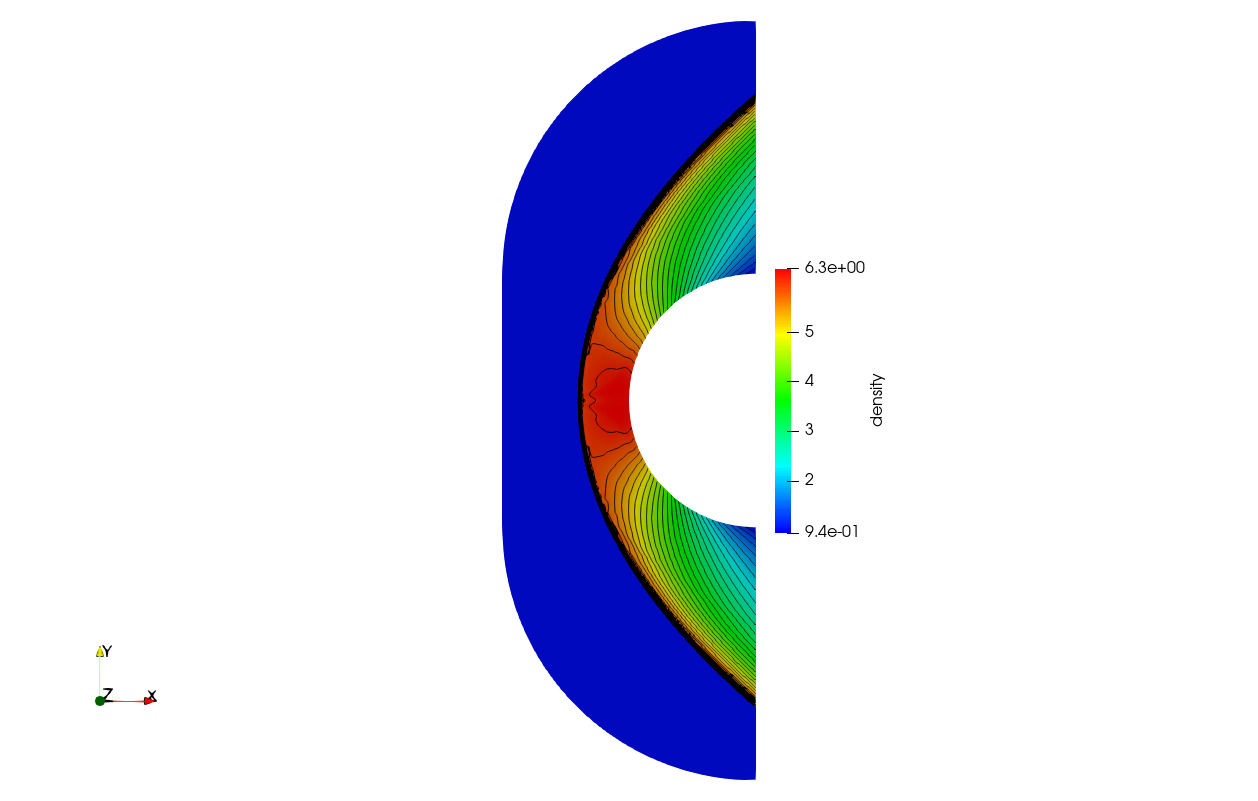}
		\subcaption{Density}
		\label{fig:BowShockden}
	\end{subfigure}
	\begin{subfigure}[c]{0.25\textwidth}
		\includegraphics[trim={17cm 0cm 13.7cm 0cmm},clip, width = 0.95\textwidth]{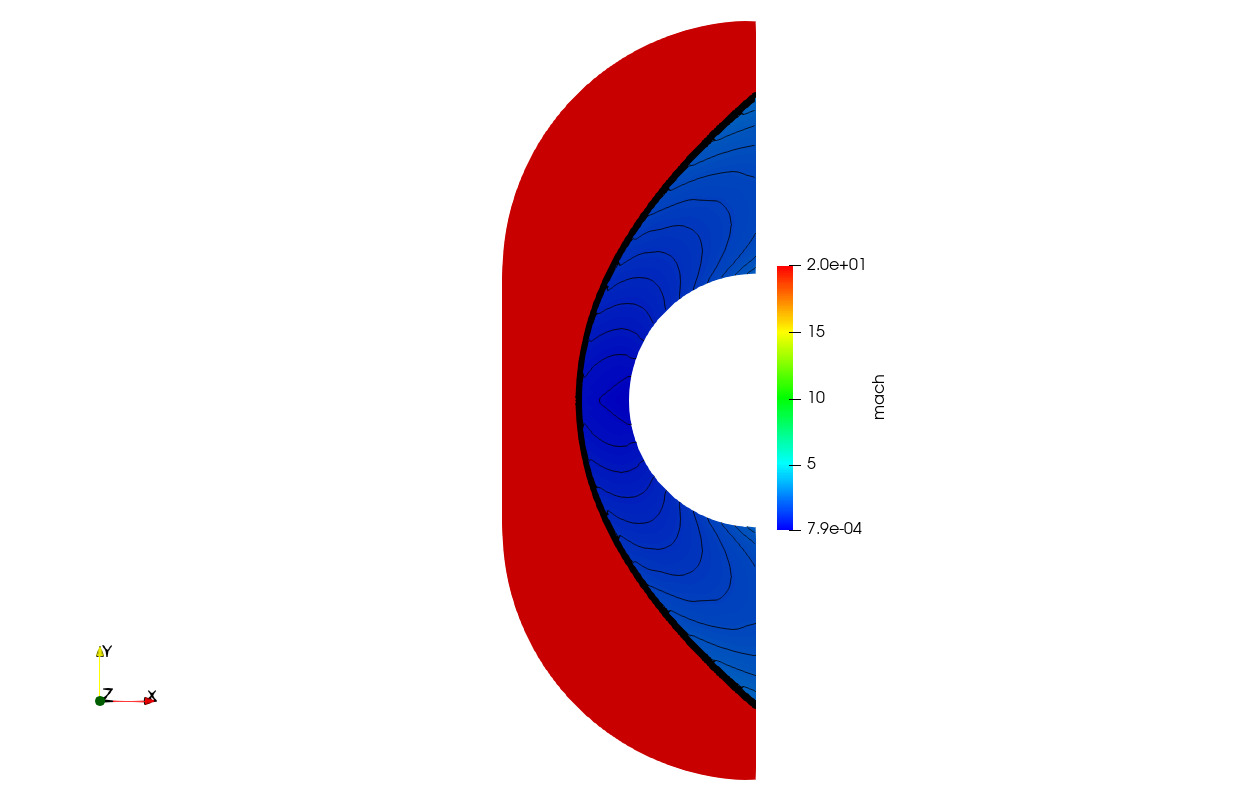}
		\subcaption{Mach number}
		\label{fig:BowShockmach}
	\end{subfigure}
	\begin{subfigure}[c]{0.25\textwidth}
		\includegraphics[trim={17cm 0cm 13.7cm 0cmm},clip, width = 0.95\textwidth]{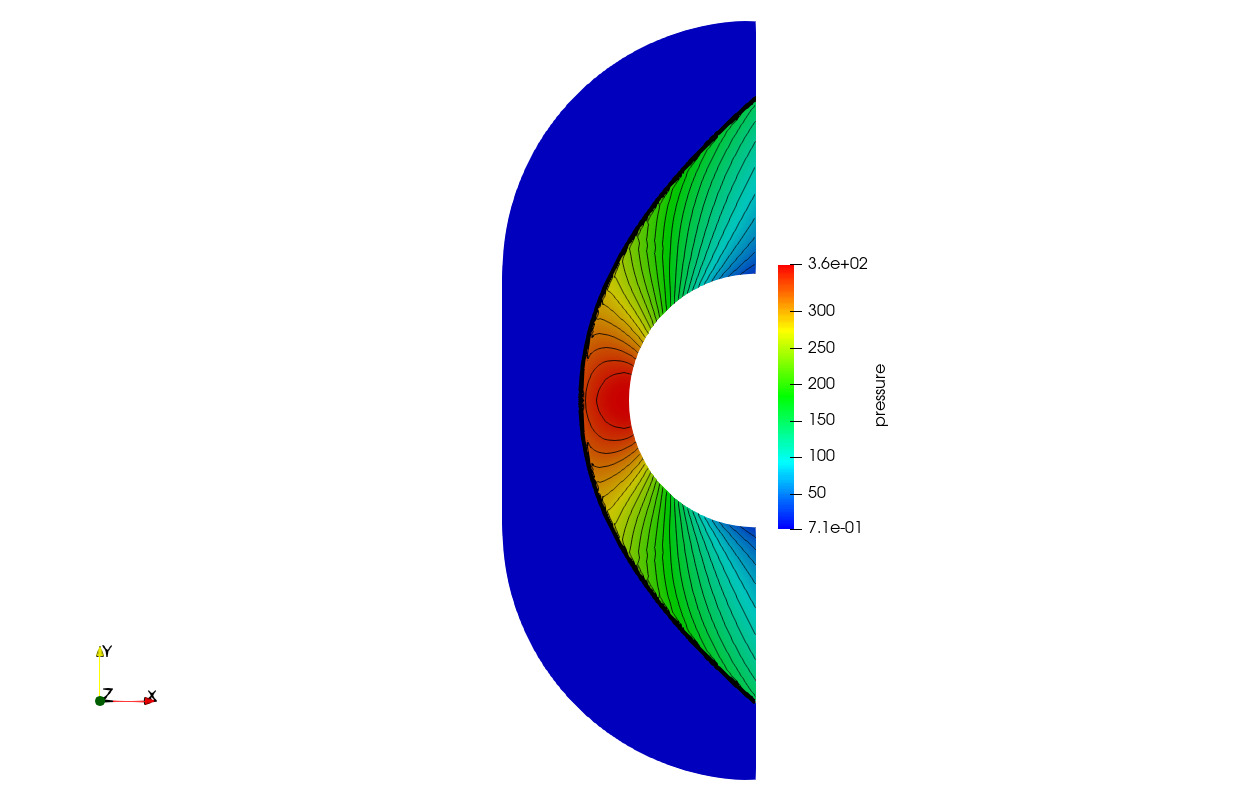}
		\subcaption{Pressure}
		\label{fig:BowShockpres}
	\end{subfigure}
	
	\caption{Bow shock ($M_\infty = 20$), steady-state MCL results obtained with $N_h = 105,361$ unknowns per component on a mesh consisting of $E_h = 209,152$ triangles.}
	\label{fig:BowShock}
\end{figure}

\begin{figure}[h!]
	\centering
	\begin{subfigure}[c]{0.40\textwidth}
		\includegraphics[trim={0.9cm 0.3cm 2.15cm 1cmm},clip, width = 0.95\textwidth]{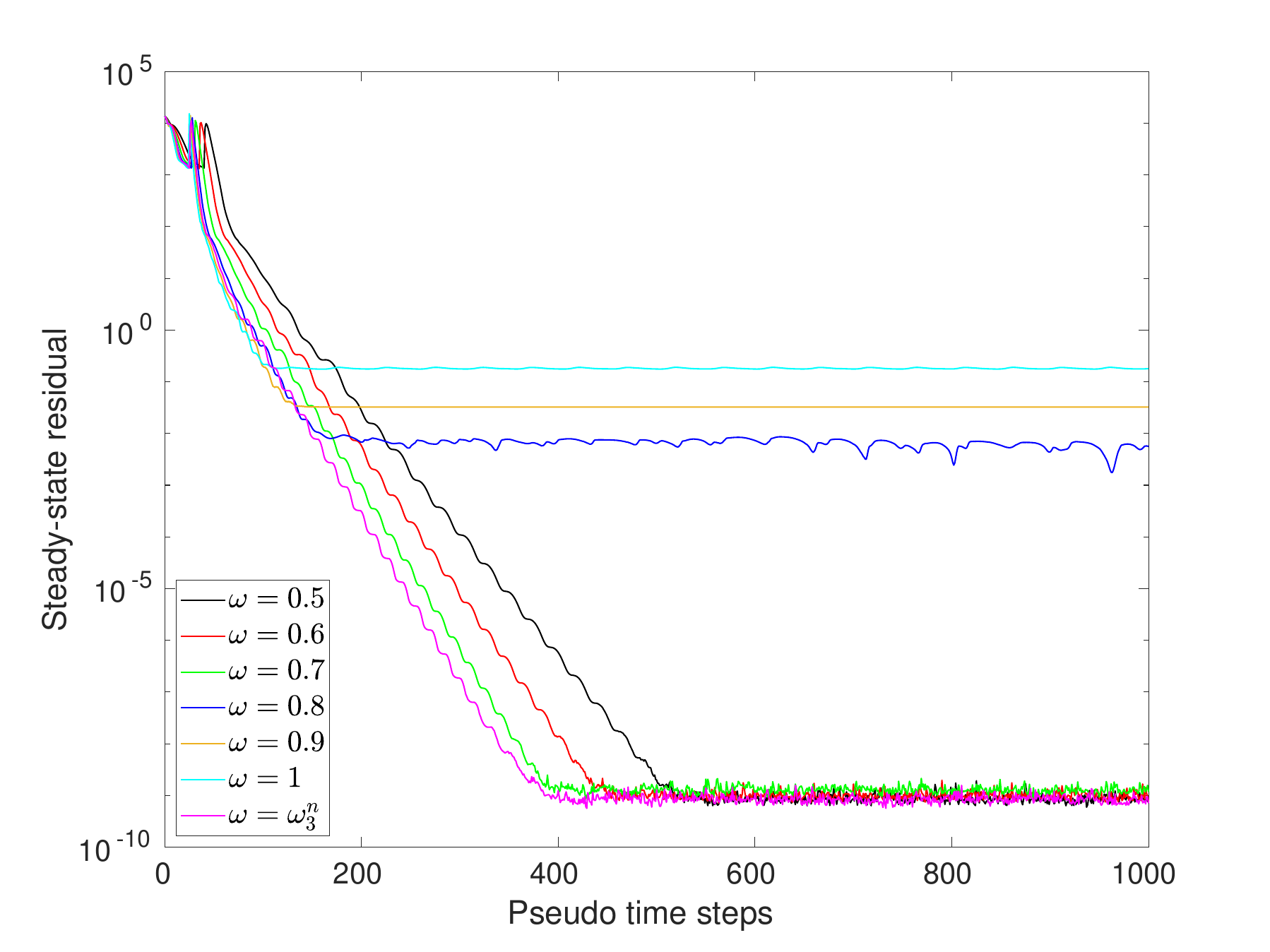}
		\subcaption{Different underrelaxation factors, $\mathrm{CFL} = 10^4$.}
		\label{fig:BowShock_ur}
	\end{subfigure}
	\begin{subfigure}[c]{0.40\textwidth}
		\includegraphics[trim={0.9cm 0.3cm 2.15cm 1cmm},clip, width = 0.95\textwidth]{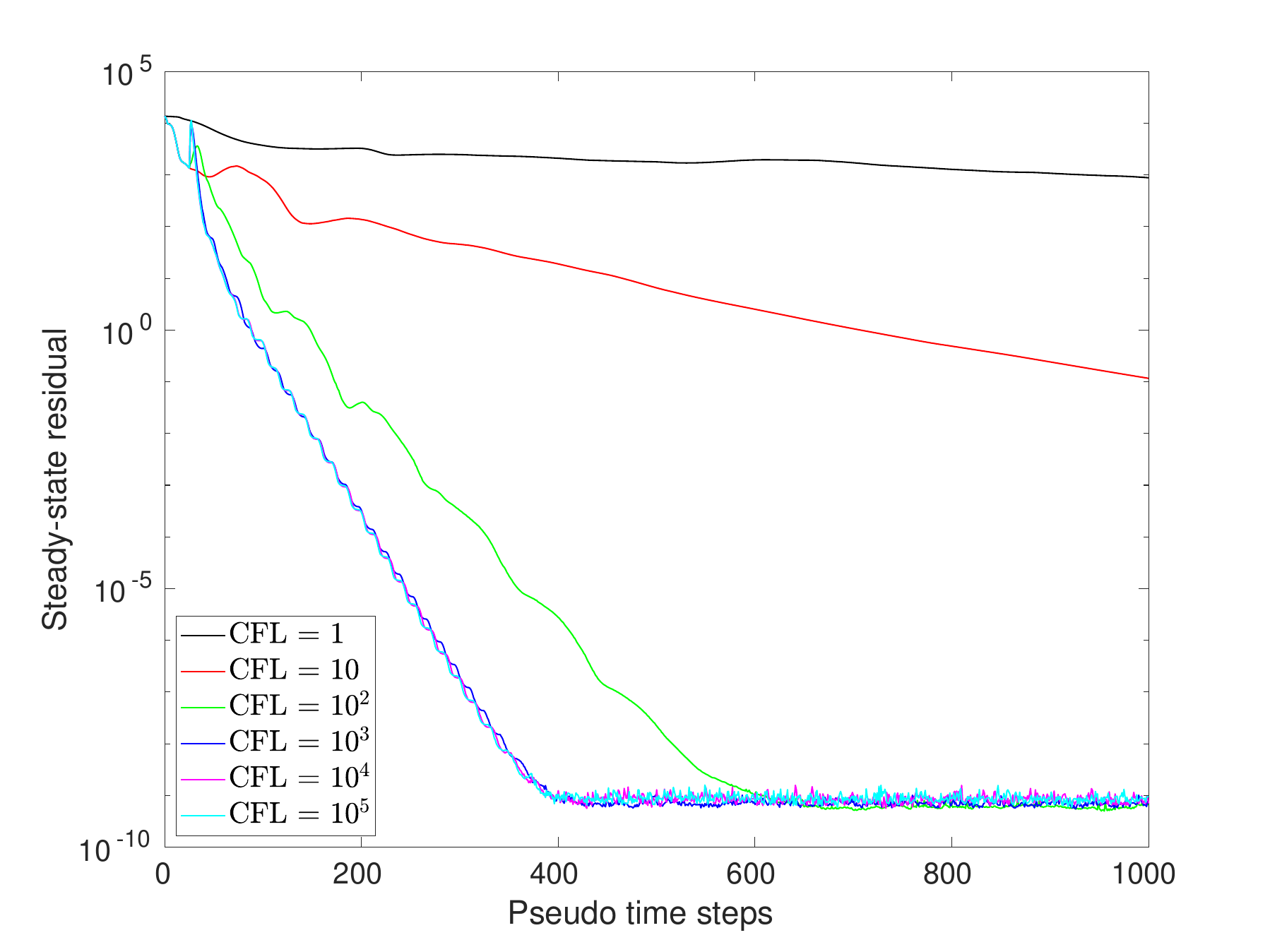}
		\subcaption{Different CFL numbers, $\omega = \omega_3^n$.}
		\label{fig:BowShock_CFL}
	\end{subfigure}
	\caption{Bow Shock ($M_\infty= 20$), steady-state convergence history for a long-time MCL simulation on an unstructured triangular mesh with $N_h = 26,537$ nodes and $E_h = 52,288$ cells.}
	\label{fig:BowShock_SSR}
\end{figure}

We finish with an example in which we simulate a steady two-dimensional hypersonic flow around a half-cylinder of unit radius \cite{cossart2024}. The free stream values at the supersonic inlet correspond to $M_\infty = 20$. At the supersonic outlet, the external state of the weakly imposed boundary condition coincides with the internal limit. The abrupt transition from hypersonic to subsonic flow conditions at the bow of the cylinder makes this benchmark problem very challenging for linearized implicit solvers. 
To avoid divergence of nonlinear and linear iterations in the early stage of computations\footnote{\blue{If we directly start with $\Delta t$ corresponding to $\mathrm{CFL}\gg 10$, then (linear or nonlinear) iterations diverge or too many iterations per pseudo-time step are needed to obtain a positivity-preserving result.}}, \red{we always start with $\mathrm{CFL} \le 10$ in this test. If the target CFL number exceeds 10, we use the pseudo-time step $\Delta t_0$ corresponding to $\mathrm{CFL}=10$
until the relative residual falls below $10^{-1}$, that is, until
\begin{equation*}
	\|M_L^{-1}R_\infty^*(u^n)\|_{2,h}< 10^{-1}\|M_L^{-1}R_\infty^*(u^0)\|_{2,h},
\end{equation*}
where $u^0$ is the low-order solution. Then we increase $\Delta t$ to match
the target CFL number.
Underrelaxation is performed at the end of each time step for any value of the CFL parameter.}
The results presented in Figs.~\ref{fig:BowShock} and
\ref{fig:BowShock_SSR} illustrate the ability of our scheme to handle
arbitrary Mach numbers. 

\section{Conclusions}
\label{sec:concl}
The analysis and algorithms presented in this paper illustrate the potential of
monolithic convex limiting as an algebraic flux correction tool for implicit
finite element discretizations of nonlinear hyperbolic systems. In addition
to proving the existence of a positivity-preserving solution to the nonlinear
discrete problem, we designed an iterative solver that exhibits robust
convergence behavior and is well suited for steady-state computations. Two
essential features of this solver are the stopping criterion based on
positivity preservation and the use of relaxation parameters that minimize
entropy residuals at steady state. The robustness of our method with respect 
to CFL and Mach numbers was verified numerically. Further improvements can
presumably be achieved by using adaptive/local time stepping and local
preconditioning aimed at reducing the characteristic condition number.

\cyan{
\section*{Acknowledgments}
The authors thank three anonymous reviewers for many insightful remarks and helpful suggestions that have greatly improved this manuscript.}

\bibliographystyle{plain}
\bibliography{impcg_paper}

\begin{thebibliography}{10}

\bibitem{abgrall2017b}
R{\'e}mi Abgrall.
\newblock High order schemes for hyperbolic problems using globally continuous
  approximation and avoiding mass matrices.
\newblock {\em J. Sci. Comput.}, 73(2):461--494, 2017.

\bibitem{anderson2021}
Robert Anderson, Julian Andrej, Andrew Barker, Jamie Bramwell, Jean-Sylvain
  Camier, Jakub Cerveny, Veselin Dobrev, Yohann Dudouit, Aaron Fisher, Tzanio
  Kolev, Will Pazner, Mark Stowell, Vladimir Tomov, Ido Akkerman, Johann Dahm,
  David Medina, and Stefano Zampini.
\newblock {MFEM: A} modular finite element methods library.
\newblock {\em Comput. Math. Appl.}, 81:42--74, 2021.

\bibitem{andrej2024}
Julian Andrej, Nabil Atallah, Jan-Phillip B{\"a}cker, Jean-Sylvain Camier,
  Dylan Copeland, Veselin Dobrev, Yohann Dudouit, Tobias Duswald, Brendan
  Keith, Dohyun Kim, et~al.
\newblock High-performance finite elements with mfem.
\newblock {\em The International Journal of High Performance Computing
  Applications}, page 10943420241261981, 2024.

\bibitem{ayachit2015}
Utkarsh Ayachit.
\newblock {\em The ParaView Guide: A Parallel Visualization Application}.
\newblock Kitware, 2015.

\bibitem{badia2017a}
Santiago Badia and Jes{\'u}s Bonilla.
\newblock Monotonicity-preserving finite element schemes based on
  differentiable nonlinear stabilization.
\newblock {\em Computer Methods Appl. Meth. Engrg.}, 313:133--158, 2017.

\bibitem{barrenechea2016}
Gabriel~R Barrenechea, Volker John, and Petr Knobloch.
\newblock Analysis of algebraic flux correction schemes.
\newblock {\em SIAM J. Numer. Anal.}, 54(4):2427--2451, 2016.

\bibitem{batten1997}
Paul Batten, Nicholas Clarke, Claire Lambert, and Derek~M Causon.
\newblock On the choice of wavespeeds for the {HLLC} {R}iemann solver.
\newblock {\em SIAM Journal on Scientific Computing}, 18(6):1553--1570, 1997.

\bibitem{berthon2008}
Christophe Berthon.
\newblock An invariant domain preserving {MUSCL} scheme.
\newblock In Luis~L. Bonilla, Miguel Moscoso, Gloria Platero, and Jose~M. Vega,
  editors, {\em Progress in Industrial Mathematics at ECMI 2006}, pages
  933--938, Berlin, Heidelberg, 2008. Springer Berlin Heidelberg.

\bibitem{bouchut2004}
Fran{\c{c}}ois Bouchut.
\newblock {\em Nonlinear Stability of Finite Volume Methods for Hyperbolic
  Conservation Laws {a}nd Well-Balanced Schemes for Sources}.
\newblock Springer Science \& Business Media, 2004.

\bibitem{burton1998}
T.~A. Burton.
\newblock A fixed-point theorem of {K}rasnoselskii.
\newblock {\em Appl. Math. Letters}, 11(1):85--88, 1998.

\bibitem{chalabi1997}
Abdallah Chalabi.
\newblock On convergence of numerical schemes for hyperbolic conservation laws
  with stiff source terms.
\newblock {\em Mathematics of computation}, 66(218):527--545, 1997.

\bibitem{cossart2024}
Benoît Cossart, Jean-Philippe Braeunig, and Rapha\"el Loub\`ere.
\newblock Toward robust linear implicit schemes for steady state hypersonic
  flow.
\newblock {\em {P}reprint, available at SSRN:
  \url{http://dx.doi.org/10.2139/ssrn.4820055}}, 2015.

\bibitem{dobrev2018}
Veselin Dobrev, Tzanio Kolev, Dmitri Kuzmin, Robert Rieben, and Vladimir Tomov.
\newblock Sequential limiting in continuous and discontinuous {G}alerkin
  methods for the {E}uler equations.
\newblock {\em J. Comput. Phys.}, 356:372--390, 2018.

\bibitem{dolejvsi2004}
V~Dolej{\v{s}}{\i} and M~Feistauer.
\newblock A semi-implicit discontinuous {G}alerkin finite element method for
  the numerical solution of inviscid compressible flow.
\newblock {\em Journal of Computational Physics}, 198(2):727--746, 2004.

\bibitem{dolejvsi2015}
V{\'\i}t Dolej{\v{s}}{\'\i} and Miloslav Feistauer.
\newblock Discontinuous {G}alerkin {M}ethod.
\newblock {\em Analysis and Applications to Compressible Flow. Springer Series
  in Computational Mathematics}, 48:234, 2015.

\bibitem{feistauer2003}
Miloslav Feistauer, Ji{\v{r}}{\'i} Felcman, and Ivan Stra{\v{s}}kraba.
\newblock {\em Mathematical and Computational Methods for Compressible Flow}.
\newblock Oxford University Press, USA, 2003.

\bibitem{ferziger2002}
Joel~H. Ferziger and Milovan Peric.
\newblock {\em Computational Methods for Fluid Dynamics}.
\newblock Springer, 3 edition, 2002.

\bibitem{frid2001}
Hermano Frid.
\newblock Maps of convex sets and invariant regions for finite-difference
  systems of conservation laws.
\newblock {\em Archive for Rational Mechanics and Analysis}, 160(3):245--269,
  2001.

\bibitem{geuzaine2009}
Christophe Geuzaine and Jean-Fran{\c{c}}ois Remacle.
\newblock Gmsh: A 3-d finite element mesh generator with built-in pre-and
  post-processing facilities.
\newblock {\em Int. J. Numer. Methods Eng.}, 79(11):1309--1331, 2009.

\bibitem{gottlieb2001}
Sigal Gottlieb, Chi-Wang Shu, and Eitan Tadmor.
\newblock Strong stability-preserving high-order time discretization methods.
\newblock {\em SIAM Rev.}, 43(1):89--112, 2001.

\bibitem{guermond2018}
Jean-Luc Guermond, Murtazo Nazarov, Bojan Popov, and Ignacio Tomas.
\newblock Second-order invariant domain preserving approximation of the {E}uler
  equations using convex limiting.
\newblock {\em SIAM J. Sci. Comput.}, 40(5):A3211--A3239, 2018.

\bibitem{guermond2016a}
Jean-Luc Guermond and Bojan Popov.
\newblock Fast estimation from above of the maximum wave speed in the {R}iemann
  problem for the {E}uler equations.
\newblock {\em J. Comput. Phys.}, 321:908--926, 2016.

\bibitem{guermond2016}
Jean-Luc Guermond and Bojan Popov.
\newblock Invariant domains and first-order continuous finite element
  approximation for hyperbolic systems.
\newblock {\em SIAM J. Numer. Anal.}, 54(4):2466--2489, 2016.

\bibitem{gurris2009}
Marcel Gurris.
\newblock {\em Implicit Finite Element Schemes for Compressible Gas and
  Particle-Laden Gas Flows}.
\newblock PhD thesis, TU Dortmund Univerity, 2009.

\bibitem{gurris2012}
Marcel Gurris, Dmitri Kuzmin, and Stefan Turek.
\newblock Implicit finite element schemes for the stationary compressible
  {E}uler equations.
\newblock {\em Int. J. Numer. Meth. Fluids}, 69(1):1--28, 2012.

\bibitem{hajduk2022diss}
Hennes Hajduk.
\newblock {\em Algebraically Constrained Finite Element Methods for Hyperbolic
  Problems with Applications in Geophysics and Gas Dynamics}.
\newblock PhD thesis, TU Dortmund University, 2022.

\bibitem{harten1983b}
Ami Harten, Peter~D. Lax, and Bram~van Leer.
\newblock On upstream differencing and {G}odunov-type schemes for hyperbolic
  conservation laws.
\newblock {\em SIAM review}, 25(1):35--61, 1983.

\bibitem{hartmann2002}
Ralf Hartmann and Paul Houston.
\newblock Adaptive discontinuous {G}alerkin finite element methods for the
  compressible {E}uler equations.
\newblock {\em Journal of Computational Physics}, 183(2):508--532, 2002.

\bibitem{hoff1979}
David Hoff.
\newblock A finite difference scheme for a system of two conservation laws with
  artificial viscosity.
\newblock {\em Math. Comp.}, 33(148):1171--1193, 1979.

\bibitem{hypre}
{\sl hypre}: High performance preconditioners.
\newblock \url{https://llnl.gov/casc/hypre},
  \url{https://github.com/hypre-space/hypre}.

\bibitem{knoll2004}
Dana~A Knoll and David~E Keyes.
\newblock Jacobian-free {N}ewton--{K}rylov methods: a survey of approaches and
  applications.
\newblock {\em Journal of Computational Physics}, 193(2):357--397, 2004.

\bibitem{kucera2022}
V{\'a}clav Ku{\v{c}}era, M{\'a}ria Luk{\'a}{\v{c}}ov{\'a}-Medvid’ov{\'a},
  Sebastian Noelle, and Jochen Sch{\"u}tz.
\newblock Asymptotic properties of a class of linearly implicit schemes for
  weakly compressible {E}uler equations.
\newblock {\em Numerische Mathematik}, pages 1--25, 2022.

\bibitem{kuzmin2020}
Dmitri Kuzmin.
\newblock Monolithic convex limiting for continuous finite element
  discretizations of hyperbolic conservation laws.
\newblock {\em Computer Methods Appl. Meth. Engrg.}, 361:112804, 2020.

\bibitem{kuzmin2023}
Dmitri Kuzmin and Hennes Hajduk.
\newblock {\em Property-Preserving Numerical Schemes for Conservation Laws}.
\newblock World Scientific, 2023.

\bibitem{kuzmin2012b}
Dmitri Kuzmin, Matthias M{\"o}ller, and Marcel Gurris.
\newblock Algebraic flux correction {II}. {C}ompressible flow problems.
\newblock In Dmitri Kuzmin, Rainald L{\"o}hner, and Stefan Turek, editors, {\em
  Flux-Corrected Transport: Principles, Algorithms, and Applications}, pages
  193--238. Springer, 2 edition, 2012.

\bibitem{kuzmin2010a}
Dmitri Kuzmin, Matthias M{\"o}ller, John~N Shadid, and Mikhail Shashkov.
\newblock Failsafe flux limiting and constrained data projections for equations
  of gas dynamics.
\newblock {\em J. Comput. Phys.}, 229(23):8766--8779, 2010.

\bibitem{leveque2002}
Randall~J. LeVeque.
\newblock {\em Finite-Volume Methods for Hyperbolic Problems}.
\newblock Cambridge University Press, 2002.

\bibitem{lin2023}
Yimin Lin, Jesse Chan, and Ignacio Thomas.
\newblock A positivity preserving strategy for entropy stable discontinuous
  {G}alerkin discretizations of the compressible {E}uler and {N}avier--{S}tokes
  equations.
\newblock {\em Journal of Computational Physics}, 475:111850, 2023.

\bibitem{lohmann-preprint}
Christoph Lohmann.
\newblock On the solvability and iterative solution of algebraic flux
  correction problems for convection-reaction equations.
\newblock {\em {P}reprint, Ergebnisber. Angew. Mathematik 612, TU Dortmund
  University}, 2019.

\bibitem{lohmann2019}
Christoph Lohmann.
\newblock {\em Physics-Compatible Finite Element Methods for Scalar and
  Tensorial Advection Problems}.
\newblock Springer Spektrum, 2019.

\bibitem{lohmann2021}
Christoph Lohmann.
\newblock An algebraic flux correction scheme facilitating the use of
  {N}ewton-like solution strategies.
\newblock {\em Comput. Math. Appl.}, 84:56--76, 2021.

\bibitem{lohmann2016}
Christoph Lohmann and Dmitri Kuzmin.
\newblock Synchronized flux limiting for gas dynamics variables.
\newblock {\em J. Comput. Phys.}, 326:973--990, 2016.

\bibitem{mfem}
{MFEM:} {M}odular {F}inite {E}lement {M}ethods [{S}oftware].
\newblock \url{https://mfem.org}.

\bibitem{patankar1980}
Suhas Patankar.
\newblock {\em Numerical Heat Transfer and Fluid Flow}.
\newblock CRC press, 1980.

\bibitem{perthame1996}
Benoit Perthame and Chi-Wang Shu.
\newblock On positivity preserving finite volume schemes for {E}uler equations.
\newblock {\em Numer. Math.}, 73(1):119--130, 1996.

\bibitem{ranocha2020}
Hendrik Ranocha, Mohammed Sayyari, Lisandro~D. Dalcin, Matteo Parsani, and
  David~I. Ketcheson.
\newblock Relaxation {R}unge--{K}utta methods: Fully discrete explicit
  entropy-stable schemes for the compressible {E}uler and {N}avier--{S}tokes
  equations.
\newblock {\em SIAM J. Sci. Comput.}, 42(2):A612--A638, 2020.

\bibitem{rueda2024}
Andr{\'e}s~M. Rueda-Ram{\'\i}rez, Benjamin Bolm, Dmitri Kuzmin, and Gregor~J.
  Gassner.
\newblock Monolithic convex limiting for {L}egendre-{G}auss-{L}obatto
  discontinuous {G}alerkin spectral-element methods.
\newblock {\em Commun. Appl. Math. \& Computation}, pages 1--39, 2024.

\bibitem{selmin1996}
V.~Selmin and L.~Formaggia.
\newblock Unified construction of finite element and finite volume
  discretizations for compressible flows.
\newblock {\em Int. J. Numer. Methods Eng.}, 39:1--32, 1996.

\bibitem{selmin1993}
Vittorio Selmin.
\newblock The node-centred finite volume approach: Bridge between finite
  differences and finite elements.
\newblock {\em Comput. Methods Appl. Mech. Engrg.}, 102(1):107--138, 1993.

\bibitem{tang2000}
Hua-Zhong Tang and Kun Xu.
\newblock Positivity-preserving analysis of explicit and implicit
  {L}ax--{F}riedrichs schemes for compressible {E}uler equations.
\newblock {\em J. Sci. Computing}, 15:19--28, 2000.

\bibitem{toro2020}
Eleuterio~F Toro, Lucas~O M{\"u}ller, and Annunziato Siviglia.
\newblock Bounds for wave speeds in the {R}iemann problem: {D}irect theoretical
  estimates.
\newblock {\em Computers \& Fluids}, 209:104640, 2020.

\bibitem{tovar2023}
Eric~Joseph Tovar, Jean-Luc Guermond, Bennett Clayton, Matthias Maier, and
  Bojan Popov.
\newblock Robust second-order approximation of the compressible euler equations
  with an arbitrary equation of state.
\newblock {\em Journal of Computational Physics}, 478:111926, 2023.

\bibitem{turek1999}
Stefan Turek.
\newblock {\em Efficient Solvers for Incompressible Flow Problems: An
  Algorithmic and Computational Approach}.
\newblock Springer, 1999.

\bibitem{zhang2017}
Xiangxiong Zhang.
\newblock On positivity-preserving high order discontinuous {G}alerkin schemes
  for compressible {N}avier--{S}tokes equations.
\newblock {\em J. Comput. Phys.}, 328:301--343, 2017.

\end{thebibliography}

\end{document}